\documentclass[10pt]{smfart}

\RequirePackage[T1]{fontenc}
\RequirePackage{amsfonts,latexsym,amssymb}
\RequirePackage[frenchb]{babel}
\addto\extrasfrenchb{\bbl@nonfrenchitemize
\bbl@nonfrenchspacing}

\RequirePackage{mathrsfs}
\let\mathcal\mathscr

\RequirePackage{smfenum}

\theoremstyle{remark}
\newtheorem{defi}{\definame}[section]
\newtheorem{rema}[defi]{\remaname}

\newtheorem{ques}[defi]{Question}
\theoremstyle{plain}
\newtheorem{prop}[defi]{\propname}
\newtheorem{theo}[defi]{\theoname}

\newtheorem{coro}[defi]{\coroname}

\newtheorem{lemm}[defi]{\lemmname}
\theoremstyle{definition}

\makeatother

\usepackage[matrix,arrow]{xy}
\usepackage{url}

\let\cal\mathcal

\def\wct{\widehat{\cal T}}

\def\demo{\noindent{\itshape D\'emonstration}.\ ---\ \ignorespaces}
\def\Q{{\bf Q}} \def\Z{{\bf Z}}

\def\N{{\bf N}}

\def\O{{\cal O}}
\def\Qbar{{\overline{{\bf Q}}}}
\def\dual{{\boldsymbol *}}

\def\zp{{\Z_p}}
\def\zpet{{\Z_p^\dual}}
\def\qp{{\Q_p}}
\def\qpet{{\Q_p^\dual}}
\def\p1{{\bf P}^1}

\def\matrice#1#2#3#4{{\big(\begin{smallmatrix}#1&#2\\ #3&#4\end{smallmatrix}\big)}}

\def\linv{{\lim\limits_{{\longleftarrow}}}}
\def\epsilon{\varepsilon}

\def\fget{\Phi\Gamma^{\rm et}}
\def\G{{\cal G}} 

\def\edag{{\cal E}^\dagger}
\def\oe{\O_{\cal E}}

\def\piqp{{{\bf P}^1}}

\def\oed#1{{\oe^{\dagger,#1}}}

\def\reso{{\text{\rm r\'es}_0}}

\begin{document}
\title[Compl\'et\'es universels]
{Compl\'et\'es universels de repr\'esentations de ${\bf GL}_2(\Q_p)$}
\author{Pierre Colmez}
\address{C.N.R.S., Institut de math\'ematiques de Jussieu, 4 place Jussieu,
75005 Paris, France}
\email{colmez@math.jussieu.fr}
\author{Gabriel Dospinescu}
\address{UMPA, \'Ecole Normale Sup\'erieure de Lyon, 46 all\'ee d'Italie, 69007 Lyon, France}
\email{gabriel.dospinescu@ens-lyon.fr}
\begin{abstract}
Soit $\Pi$ une repr\'esentation unitaire
de ${\bf GL}_2(\Q_p)$, topologiquement de longueur finie.
Nous d\'ecrivons la sous-repr\'esentation $\Pi^{\rm an}$
de ses vecteurs localement analytiques, et sa filtration par rayon d'analyticit\'e,
 en termes du $(\varphi,\Gamma)$-module
qui lui est associ\'e via la correspondance de Langlands locale $p$-adique, et
nous en d\'eduisons
que le compl\'et\'e universel de $\Pi^{\rm an}$ n'est autre que $\Pi$.
\end{abstract}
\begin{altabstract}
Let $\Pi$ be a unitary
representation of ${\bf GL}_2(\Q_p)$, topologically of finite length.
We describe the sub-representation $\Pi^{\rm an}$ 
made of its locally analytic vectors, and its filtration by radius of
analyticity, in terms of the
$(\varphi,\Gamma)$-module attached to $\Pi$ via the
$p$-adic local Langlands correspondence, and we deduce
that the universal completion of $\Pi^{\rm an}$ is
$\Pi$ itself.
\end{altabstract}
\setcounter{tocdepth}{3}

\maketitle

\stepcounter{tocdepth}
{\Small
\tableofcontents
}

\section*{Introduction}

\Subsection{Notations}
Soit $p$ un nombre premier.
On fixe une cl\^oture alg\'ebrique $\Qbar_p$ de $\Q_p$,
et on note $\G_{\Q_p}={\rm Gal}(\Qbar_p/\Q_p)$ le groupe de Galois
absolu de $\Q_p$. On note $\Gamma$ le groupe
de Galois ${\rm Gal}(\qp(\mu_{p^{\infty}})/\Q_p)$ de l'extension cyclotomique.
 Le 
caract\`ere cyclotomique 
\linebreak
$\chi:\G_{\Q_p}\to \Z_p^\dual$ induit un isomorphisme
de groupes topologiques $ \Gamma\simeq\zpet$, dont l'inverse 
$a\mapsto\sigma_a$ est caract\'eris\'e par $\sigma_a(\zeta)=\zeta^a$ pour
$a\in \zpet$ et $\zeta\in\mu_{p^{\infty}}$.

On fixe une extension finie $L$ de $\Q_p$, et on note $\O_L$ l'anneau
de ses entiers et $k_L$ son corps r\'esiduel.
Soit $\wct(L)$ l'ensemble des caract\`eres continus
$\delta:\Q_p^\dual\to L^\dual$, et, pour $\delta\in\wct(L)$, notons
$w(\delta)$ son {\it poids},
d\'efini par $w(\delta)=\delta'(1)$, d\'eriv\'ee\footnote{$\delta$ est automatiquement localement analytique, donc la d\'efinition a un sens.} de $\delta$ en $1$.
Si $\delta$ est {\it unitaire} (i.e.~si $\delta$ est \`a valeurs dans $\O_L^\dual$),
la th\'eorie locale du corps de classes associe \`a $\delta$ un caract\`ere continu de
 $\G_{\Q_p}$, que l'on note encore $\delta$. Le poids de Hodge-Tate g\'en\'eralis\'e de ce caract\`ere galoisien est alors $w(\delta)$.
On note juste $x\in\wct(L)$ le caract\`ere induit
par l'inclusion de $\Q_p$ dans $L$,
et $|x|$ le caract\`ere
envoyant $x\in\Q_p^\dual$ sur~$p^{-v_p(x)}$.
Notre convention est que
$x|x|$ correspond
au caract\`ere cyclotomique~$\chi$, son poids \'etant~$1$.

Soit $G={\bf GL}_2(\Q_p)$. Si $\delta\in \wct(L)$, on note
${\rm Rep}_L(\delta)$ la cat\'egorie dont les objets sont les
$L$-espaces de Banach $\Pi$, munis d'une action continue de $G$ telle que

$\bullet$ $\Pi$ a pour caract\`ere central $\delta$.

$\bullet$ L'action de $G$ est {\it unitaire},
 i.e. il existe une valuation $v_\Pi$ sur $\Pi$, qui d\'efinit
la topologie de $\Pi$ et telle que $v_\Pi(g\cdot v)=v_\Pi(v)$ pour tous $g\in G$ et $v\in \Pi$.

$\bullet$ $\Pi$ est {\it r\'esiduellement de longueur finie}, 
i.e. si $v_\Pi$ est comme ci-dessus, 
la r\'eduction mod $p$ de la boule unit\'e $\Pi_0$ de $\Pi$ pour $v_\Pi$ 
est un $\O_L[G]$-module de longueur finie.

Un morphisme entre deux objets $\Pi_1$ et $\Pi_2$ de ${\rm Rep}_L(\delta)$ est une application 
\linebreak
$L$-lin\'eaire, continue et
$G$-\'equivariante.

On note ${\rm Rep}_L(G)$ la cat\'egorie des repr\'esentations de $G$ unitaires, r\'esiduellement
de longueur finie, admettant un caract\`ere central; c'est la r\'eunion des ${\rm Rep}_L(\delta)$
pour $\delta\in \wct(L)$.

\begin{rema} 

 {\rm (i)} ${\rm Rep}_L(\delta)$ est vide quand $\delta$ n'est pas unitaire.

 {\rm (ii)} Il d\'ecoule des travaux de Barthel-Livn\'e \cite{BL} et Breuil \cite{Br1} que tout objet de ${\rm Rep}_L(\delta)$ est une repr\'esentation de Banach {\it admissible} (au sens de~\cite{ST1}) de $G$.

 {\rm (iii)}  Tout objet $\Pi$ de ${\rm Rep}_L(\delta)$ est topologiquement de longueur finie,
   car $\Pi_0/p\Pi_0$ l'est. En fait, si $p\geq 5$ (il est probable que cette hypoth\`ese est inutile),
   on peut d\'ecrire ${\rm Rep}_L(\delta)$ de mani\`ere \'equivalente comme la cat\'egorie des $L$-repr\'esentations de Banach unitaires et admissibles de $G$, \`a caract\`ere central et topologiquement de longueur finie (cela d\'ecoule d'un th\'eor\`eme profond de Paskunas \cite{Pa}).
    Nous avons besoin de cette hypoth\`ese de finitude pour $\Pi_0/p\Pi_0$ pour assurer que le $(\varphi,\Gamma)$-module
   attach\'e \`a $\Pi$ par la correspondance de Langlands locale $p$-adique est de dimension finie.
\end{rema}

\subsection{Compl\'etions unitaires et vecteurs localement analytiques}

  Si 
\linebreak
$\delta\in \wct(L)$ et si $\Pi\in {\rm Rep}_L(\delta)$, on note
  $\Pi^{\rm an}$ l'espace des vecteurs localement analytiques~\cite{STInv,Emlocan} de~ $\Pi$.
  C'est l'espace des vecteurs $v\in \Pi$ dont l'application orbite
  $$o_v: G\to \Pi, \quad g\to g\cdot v$$ est localement analytique.
C'est une sous-repr\'esentation de $\Pi$ et,
  d'apr\`es un r\'esultat g\'en\'eral
de Schneider et Teitelbaum \cite[th. 7.1]{STInv}, le sous-espace
  $\Pi^{\rm an}$ est dense\footnote{  Nous donnons une nouvelle preuve de cette densit\'e pour les objets de
  ${\rm Rep}_L(\delta)$, en utilisant la th\'eorie des $(\varphi,\Gamma)$-modules 
(cf.~cor.~\ref{dense1}).}  
dans~$\Pi$.

   L'espace $\Pi^{\rm an}$ a une topologie naturelle,  induite par l'injection $\Pi^{\rm an}\to {\cal C}^{\rm an}(G,\Pi)$,
   envoyant $v\in \Pi^{\rm an}$ sur $o_v$. Cette topologie est nettement plus forte que celle induite par l'inclusion 
   $\Pi^{\rm an}\subset \Pi$. 
 Le r\'esultat principal de cet article est alors le suivant (rappelons que $G={\rm GL}_2(\qp)$):

\begin{theo}\label{maincompletion}
Si $\Pi\in {\rm Rep}_L(G)$, alors 
$\Pi$ est le compl\'et\'e unitaire universel
de~$\Pi^{\rm an}$, i.e. pour toute $L$-repr\'esentation de Banach unitaire~$B$,
l'application naturelle 
${\rm Hom}_{L[G]}^{\rm cont}(\Pi, B)\to {\rm Hom}_{L[G]}^{\rm cont}(\Pi^{\rm an}, B)$,
induite par l'injection $\Pi^{\rm an}\to \Pi$,
est un isomorphisme.
\end{theo}

\begin{rema} 
 La notion de compl\'et\'e universel a \'et\'e d\'egag\'ee par
Emerton~\cite{Emcompl}. Le th\'eor\`eme ci-dessus r\'epond, dans
le cas de ${\bf GL}_2(\Q_p)$, \`a l'une de ses questions, \`a savoir
si le m\^eme \'enonc\'e est valable pour
${\bf GL}_n(\Q_p)$ ou, plus g\'en\'eralement, pour un groupe
r\'eductif d\'eploy\'e sur $\Q_p$ (l'application 
${\rm Hom}_{L[G]}^{\rm cont}(\Pi, B)\to {\rm Hom}_{L[G]}^{\rm cont}(\Pi^{\rm an}, B)$ est
 injective pour 
tout groupe de Lie $p$-adique $G$ et toute repr\'esentation de Banach admissible $\Pi$ de $G$, car 
$\Pi^{\rm an}$ est dense dans $\Pi$ dans ces cas~\cite{STInv}). 

(ii) Dans l'autre sens, la situation est nettement plus compliqu\'ee:
si $\Pi$ est une repr\'esentation localement analytique admissible de
${\bf GL}_2(\Q_p)$ admettant un compl\'et\'e unitaire universel $\widehat\Pi$,
la sous-repr\'esentation $\widehat\Pi^{\rm an}$ de $\widehat\Pi$
n'est pas forc\'ement \'egale \`a~$\Pi$. 
Le cas des composantes de Jordan-H\"older de la s\'erie principale analytique
est assez \'eclairant.
 Si $\delta_1,\delta_2:\qpet\to L^\dual $ sont des caract\`eres continus (et donc localement analytiques),
on note ${\rm Ind}^{\rm an}(\delta_1\otimes\delta_2)$
l'espace des fonctions $\phi: G\to L$, localement analytiques, telles que
$\phi\left(\left(\begin{smallmatrix} a & b \\0 & d\end{smallmatrix}\right)g\right))=\delta_1(a)
\delta_2(d)f(g)$ pour tous $a,d\in\qpet$, $b\in\qp$ et $g\in G$, que l'on munit
de l'action de $G$ d\'efinie par $(h\cdot\phi)(g)=\phi(gh)$; si $\delta_2=x^k\delta_1$
avec $k\in\N$, alors ${\rm Ind}^{\rm an}(\delta_1\otimes\delta_2)$ contient
une sous-repr\'esentation de dimension $k+1$ et le quotient est une {\it steinberg analytique}.
En utilisant les r\'esultats de~\cite{Br,Cserieunit,Cvectan,Emcompl,LXZ}, on voit qu'il
peut, en particulier, se passer les choses suivantes:

$\bullet$ $\widehat\Pi=0$ et donc $\widehat\Pi^{\rm an}=0$: c'est
le cas si le caract\`ere central n'est pas unitaire (ce qui \'equivaut \`a
$v_p(\delta_1(p))+v_p(\delta_2(p))\neq 0$) ou s'il est unitaire
mais $v_p(\delta_1(p))>0$.

$\bullet$ $\widehat\Pi$ n'est pas admissible (c'est le cas des steinberg analytiques avec $k\geq 1$).

$\bullet$ $\widehat\Pi$ est non nul et admissible, mais $\widehat\Pi^{\rm an}$ est strictement
plus grand que $\Pi$: c'est le cas si $\Pi={\rm Ind}^{\rm an}(\delta_1\otimes\delta_2)$,
si le caract\`ere central est unitaire, et si
$v_p(\delta_2(p))>0$ et $w(\delta_2)-w(\delta_1)\notin\N$.
\end{rema}

Mentionnons un corollaire imm\'ediat du th.~\ref{maincompletion}, qui ne semble pas facile \`a d\'emontrer directement:

\begin{coro}
 Le foncteur $\Pi\mapsto\Pi^{\rm an}$ de la cat\'egorie ${\rm Rep}_L(\delta)$ dans la cat\'egorie des $L$-repr\'esentations localement analytiques
 admissibles de $G$ est pleinement fid\`ele. 
 \end{coro}

  En utilisant ce corollaire et les r\'esultats de \cite{DBenjamin}, on d\'eduit que $\Pi^{\rm an}$ admet un caract\`ere infinit\'esimal
  pour tout objet absolument irr\'eductible $\Pi$ de ${\rm Rep}_L(G)$ (cela n'est pas une
cons\'equence
  formelle du r\'esultat principal de loc.cit., car $\Pi^{\rm an}$ peut fort bien 
ne pas \^etre irr\'eductible si $\Pi$ est absolument irr\'eductible).

\smallskip
  La suite de cette introduction explique les \'etapes de la preuve 
du th.~\ref{maincompletion}, dont la correspondance
  de Langlands locale $p$-adique  \cite{Cbigone} pour $G$ 
  est l'ingr\'edient cl\'e. 
  
 \subsection{Un raffinement du foncteur $\Pi\mapsto \Pi^{\rm an}$} \label{pihintro}
 
   Dans ce \S, on consid\`ere un groupe de Lie $p$-adique $G$ arbitraire, 
   un sous-groupe $H$ de $G$ qui est un pro-$p$-groupe uniforme, 
    et une $L$-repr\'esentation de Banach admissible
   $\Pi$ de $G$, pas forc\'ement unitaire. On choisit un syst\`eme minimal de g\'en\'erateurs topologiques
   $h_1,...,h_d$ de $H$ et on note 
   $$b^{\alpha}=(h_1-1)^{\alpha_1}\cdots (h_d-1)^{\alpha_d}\in \zp[H]$$
   pour $\alpha=(\alpha_1,...,\alpha_d)\in \mathbf{N}^d$. 
   Pour tout entier $h\geq 1$, on note 
   $r_h=\frac{1}{p^{h-1}(p-1)}$ et (en posant $|\alpha|=\alpha_1+\cdots+\alpha_d$)
   $$\Pi_H^{(h)}=\{v\in \Pi,\  \lim_{|\alpha|\to\infty} p^{-r_h|\alpha|} b^{\alpha}v=0\}.$$
    
    Alors $\Pi_H^{(h)}$ est naturellement un espace de Banach. Il ne d\'epend pas du choix 
des g\'en\'erateurs
    $h_1,...,h_d$ et il est stable par $H$. Par ailleurs, les $b^\alpha v$ sont les
coefficients de Mahler de $o_v$:
$$o_v(h_1^{x_1}\cdots h_d^{x_d})=\sum_{\alpha\in\N^d}\tbinom{x_1}{\alpha_1}\cdots
\tbinom{x_d}{\alpha_d}\cdot b^\alpha v,\quad{\text{pour tout $(x_1,\dots,x_d)\in\Z_p^d$.}}$$
Il r\'esulte donc du th\'eor\`eme d'Amice~\cite{Am64} que
    $\Pi^{\rm an}$ est la limite inductive des espaces $\Pi_H^{(h)}$, et cela pour tout sous-groupe $H$ de $G$ qui est un pro-$p$-groupe uniforme. 
Le r\'esultat suivant (cor~\ref{exac1} et prop.~\ref{dcalage})
peut, au langage pr\`es, se trouver dans~\cite{STInv}.
    
    \begin{theo}\label{dec}
Soient $H$ un pro-$p$ sous-groupe uniforme de $G$ et $h\geq 1$.

    {\rm (i)} 
 Le foncteur $\Pi\mapsto \Pi_H^{(h)}$ de la cat\'egorie des $L$-repr\'esentations de Banach admissibles de $G$ dans la cat\'egorie
    des $L$-espaces de Banach est exact.
    
   {\rm (ii)} On a $\Pi_H^{(h+1)}=\Pi_{H^p}^{(h)}$ 
    pour toute repr\'esentation de Banach $\Pi$ de $G$. 
    \end{theo}
   
      On dit que $\Pi$ est {\it coh\'erente} (ou $H$-coh\'erente si on veut pr\'eciser
le sous-groupe $H$ de r\'ef\'erence) s'il existe $h_0$ tel que, pour tout $h\geq h_0$, on ait
      $$\Pi_H^{(h+1)}=\sum_{H^p\subset gHg^{-1}} g\cdot \Pi_H^{(h)}.$$
Comme $g\cdot\Pi^{(h)}_H=\Pi^{(h)}_{gHg^{-1}}$ et $\Pi^{(h)}_{H_1}\subset \Pi^{(h)}_{H_2}$ si $H_2\subset H_1$,
     on d\'eduit du th\'eor\`eme pr\'ec\'edent que le terme de droite est toujours contenu dans celui de gauche. Le m\^eme th\'eor\`eme
      permet de montrer que la coh\'erence est une propri\'et\'e stable par extensions, ce qui joue un r\^ole important dans la preuve
      du th.~\ref{maincompletion}. Notre int\'er\^et pour la notion de coh\'erence vient du r\'esultat suivant:
      
      \begin{prop}
       Si $\Pi$ est $H$-coh\'erente, alors $\Pi^{\rm an}$ admet un compl\'et\'e unitaire universel. Plus pr\'ecis\'ement, si 
       $\Pi_0^{(h)}$ est la boule unit\'e de $\Pi_H^{(h)}$ et ${\cal L}_h=\sum_{g\in G} g\cdot \Pi_0^{(h)}$, alors pour tout 
       $h$ assez grand on a un isomorphisme de $L[G]$-modules de Banach      
         $$\widehat{\Pi^{\rm an}}\simeq L\otimes_{\O_L} (\linv {\cal L}_h/p^n{\cal L}_h).$$
             \end{prop}
  
       Au vu de la proposition pr\'ec\'edente, 
le th.~\ref{maincompletion} est une cons\'equence du r\'esultat suivant:
       
\begin{theo}
Si $G={\rm GL}_2(\qp)$ et si $\Pi\in {\rm Rep}_L(G)$, alors:
        
{\rm (i)} $\Pi$ est coh\'erente.
        
{\rm (ii)} Si $\Pi_0$ est un $\O_L$-r\'eseau de $\Pi$, ouvert, born\'e et $G$-stable, alors 
${\cal L}_h$ est commensurable avec $\Pi_0\cap \Pi^{\rm an}$ pour tout $h$ assez grand.
        
       \end{theo} 
  
        La preuve de ce th\'eor\`eme utilise de mani\`ere cruciale la th\'eorie des $(\varphi,\Gamma)$-modules. Plus pr\'ecis\'ement, si 
        $D$ est le $(\varphi,\Gamma)$-module attach\'e \`a $\Pi$ par la correspondance de Langlands locale $p$-adique, 
        on d\'ecrit l'espace $\Pi^{(h)}$ (et m\^eme la boule unit\'e $\Pi_0^{(h)}$) directement en termes de $D$. Cela demande 
        d'\'etendre et de raffiner bon nombre de r\'esultats des
chap.~II, IV et~V de \cite{Cbigone}, et les paragraphes suivants expliquent de quelle mani\`ere 
un plus en d\'etail.

 \subsection{Description de ${\rm Rep}_L(\delta)$ en termes de $(\varphi,\Gamma)$-modules}

 Soient ${\cal R}$
   l'anneau de Robba\footnote{ Il s'agit de l'anneau des s\'eries de Laurent $\sum_{n\in\mathbf{Z}} a_nT^n$ \`a coefficients dans $L$, qui convergent sur une couronne du type
   $0<v_p(T)\leq r$, o\`{u} $r$ d\'epend de la s\'erie.}, ${\cal E}^\dagger$ le sous-anneau de ${\cal R}$ des \'el\'ements born\'es
(c'est un corps) et ${\cal E}$ le compl\'et\'e de ${\cal E}^\dagger$ pour la valuation $p$-adique.
On munit ces anneaux d'actions continues de $\Gamma$ et d'un frobenius $\varphi$, commutant entre elles,
en posant $\varphi(T)=(1+T)^p-1$ et $\sigma_a(T)=(1+T)^a-1$ si $a\in\zpet$.

 Si $\Lambda\in \{{\cal E}, {\cal E}^{\dagger}, {\cal R}\}$, on note $\fget(\Lambda)$ la
 cat\'egorie des $(\varphi,\Gamma)$-modules \'etales sur $\Lambda$. Ce sont des $\Lambda$-modules libres de type fini
 $D$, munis d'actions de $\varphi$ et $\Gamma$, continues, semi-lin\'eaires, qui commutent et telles que $\varphi$ soit
 de pente nulle. Ces cat\'egories sont toutes \'equivalentes \`a la cat\'egorie des $L$-repr\'esentations
 de ${\rm Gal}(\overline{\qp}/\qp)$ (cf.~\cite{FoGrot,CCsurconv,KK1}); en particulier elles sont \'equivalentes entre elles, et on note
 $D^{\dagger}\in \fget({\cal E}^{\dagger})$, 
\linebreak
$D_{\rm rig}\in \fget({\cal R})$
 les $(\varphi,\Gamma)$-modules attach\'es \`a $D\in\fget({\cal E})$, de telle sorte que
\linebreak
 $D={\cal E}\otimes_{{\cal E}^{\dagger}} D^{\dagger}$ et $D_{\rm rig}={\cal R}\otimes_{{\cal E}^{\dagger}} D^{\dagger}$.

    Si $\delta$ est un caract\`ere unitaire et si $D\in\fget({\cal E})$,
    on peut construire \cite[chap.~II]{Cbigone} un faisceau $G$-\'equivariant $U\to D\boxtimes_{\delta} U$ sur
    $\p1(\qp)$ (muni de l'action usuelle d\'efinie par $\matrice{a}{b}{c}{d}\cdot x=\frac{ax+b}{cx+d}$),
 dont les sections sur $\zp$ sont $D$ (i.e. $D\boxtimes_\delta\Z_p=D$).
Par ailleurs, si $U$ est un ouvert compact de $\piqp$, l'extension par $0$ induit
une inclusion de $D\boxtimes_\delta U$ dans l'espace $D\boxtimes_\delta\piqp$ des sections globales.
 Les formules d\'ecrivant l'action de $G$ sont
    tr\`es compliqu\'ees (et inutiles dans la plupart des situations) en g\'en\'eral, mais on a par exemple, pour $z\in D=D\boxtimes_\delta\zp$, $a\in\zpet$ et $b\in\zp$,
    $$\left(\begin{smallmatrix} p & 0 \\0 & 1\end{smallmatrix}\right)z=\varphi(z), \quad \left(\begin{smallmatrix} a & 0 \\0 & 1\end{smallmatrix}\right)z=
    \sigma_a(z), \quad \left(\begin{smallmatrix} 1 & b \\0 & 1\end{smallmatrix}\right)z=(1+T)^b\cdot z.$$

   On dispose aussi \cite[chap.~IV]{Cbigone} d'un foncteur $\Pi\mapsto D(\Pi)$, contravariant, exact, de
   ${\rm Rep}_L(\delta)$ dans $\fget({\cal E})$. On note ${\cal C}_L(\delta)$ son image essentielle.
   Si $D\in\fget({\cal E})$, on note $\check{D}$ le dual de Cartier de $D$:
 si $D$ est attach\'e
   \`a une repr\'esentation galoisienne $V$, alors $\check{D}$ est attach\'e \`a $V^\dual \otimes\chi$.
   Le r\'esultat suivant fait le lien entre les constructions pr\'ec\'edentes et d\'ecrit ${\rm Rep}_L(\delta)$ (\`a des morceaux de dimension finie pr\`es)
   en termes de $(\varphi,\Gamma)$-modules, ce qui est fondamental pour la preuve du th.~\ref{maincompletion}.

\begin{theo}\label{Repdelta}
Soit $\delta:\qpet\to \O_L^\dual $ un caract\`ere unitaire.  Alors:

{\rm (i)} ${\cal C}_L(\delta)$ est stable par sous-quotients.

{\rm (ii)} Si $D\in {\cal C}_L(\delta)$, alors $\check{D}\in {\cal C}_L(\delta^{-1})$.

{\rm (iii)} Il existe un foncteur covariant $D\to \Pi_{\delta} (D)$ de ${\cal C}_L(\delta)$ dans ${\rm Rep}_L(\delta)$
tel que, pour tout $D\in {\cal C}_L(\delta)$, on ait une suite exacte de $G$-modules topologiques
 $$0\to \Pi_{\delta^{-1}}(\check{D})^\dual \to D\boxtimes_{\delta}\p1\to \Pi_{\delta}(D)\to 0.$$

{\rm (iv)} Les foncteurs $\Pi\mapsto D(\Pi)$ et $D\to \Pi_{\delta}(\check{D})$ induisent
des anti-\'equivalences quasi-inverses exactes entre ${\rm Rep}_L(\delta)/S$ et ${\cal C}_L(\delta)$,
o\`{u} $S$ est la sous-cat\'egorie de ${\rm Rep}_L(\delta)$ form\'ee des repr\'esentations de dimension finie.
\end{theo}
   Ce th\'eor\`eme est essentiellement d\'emontr\'e dans \cite{Cbigone}, mais il n'est pas facile
   de l'en extraire sous cette forme. Nous reprenons et \'etendons les arguments de
   loc.cit pour l'obtenir sous cette forme, plus adapt\'ee aux applications 
\'eventuelles.

    \subsection{Description de $\Pi^{(h)}$}
  Soit ${\cal E}^{(0,r_b]}$ le sous-anneau de $\edag$
des fonctions analytiques born\'ees sur la couronne
    $0<v_p(T)\leq r_b$. On note $D^{(0,r_b]}$ le plus grand sous-${\cal E}^{(0,r_b]}$-module de type fini $M$
    de $D$ tel que $\varphi(M)\subset {\cal E}^{(0,r_{b+1}]}\otimes_{{\cal E}^{(0,r_b]}} M$
    (son existence est un r\'esultat standard de la th\'eorie des $(\varphi,\Gamma)$-modules).
    Le th\'eor\`eme de surconvergence \cite{CCsurconv,BC} 
montre que $D^{(0,r_b]}$ est libre de rang $\dim_{{\cal E}} (D)$
    sur ${\cal E}^{(0,r_b]}$ et engendre $D$, si $b$ est assez grand. 
Comme $D\boxtimes_{\delta}\zp=D$ et comme $\piqp=\Z_p\cup \matrice{0}{1}{1}{0}\Z_p$, l'application
$z\mapsto\big({\rm Res}_{\Z_p}z, {\rm Res}_{\Z_p}\matrice{0}{1}{1}{0}z\big)$ est
une injection de
$D\boxtimes_{\delta}\p1$ dans $D\times D$, ce qui permet de d\'efinir le module
      $$D^{(0,r_b]}\boxtimes_{\delta}\p1=(D\boxtimes_{\delta}\p1)\cap (D^{(0,r_b]}\times D^{(0,r_b]}).$$
    Il est muni de la topologie induite par l'inclusion dans $D^{(0,r_b]}\times D^{(0,r_b]}$, le module
    $D^{(0,r_b]}$ \'etant muni de sa topologie naturelle.

    \begin{prop}
      Soit $D\in {\cal C}_L(\delta)$. Si $b$ est assez grand, le sous-$L$-espace
vectoriel $D^{(0,r_b]}\boxtimes_{\delta}\p1$
      de $D\boxtimes_{\delta}\p1$ est stable sous l'action de ${\rm GL}_2(\zp)$, et ${\rm GL}_2(\zp)$ agit contin\^ument pour la topologie naturelle
      de $D^{(0,r_b]}\boxtimes_{\delta}\p1$. De plus, $\Pi_{\delta^{-1}}(\check{D})^\dual $ est un sous-module ferm\'e de
      $D^{(0,r_b]}\boxtimes_{\delta}\p1$.
    \end{prop}

Soit $K_m=1+p^m{\rm M}_2(\zp)$ avec $m\geq 1$ (resp.~$m\geq 2$ si $p=2$).
 C'est un pro-$p$-groupe uniforme de dimension $4$,
     auquel les constructions du \S~\ref{pihintro} s'appliquent. Pour simplifier les notations, on note simplement 
     $\Pi^{(b)}=\Pi_{K_m}^{(b-m)}$ pour $b>m$ (le th.~\ref{dec} montre que le terme de droite ne d\'epend pas du choix de 
     $m<b$). Le r\'esultat technique principal de l'article est alors la description de
     $\Pi^{(b)}$, pour tout $b$ assez grand, si $\Pi\in{\rm Rep}_L (G)$.
Le th.~\ref{Repdelta} implique en particulier que tout objet de ${\rm Rep}_L (G)$ est
de la forme $\Pi_{\delta}(D)$ \`a des repr\'esentations de dimension finie pr\`es,
et pour une repr\'esentation de la forme $\Pi_{\delta}(D)$, on a le r\'esultat suivant.

     \begin{theo}\label{mainrayon}
      Soit $D\in {\cal C}_L(\delta)$. Alors pour tout $b$ assez grand, l'inclusion de $D^{(0,r_b]}\boxtimes_{\delta}\p1$ dans
      $D\boxtimes_{\delta}\p1$ induit une suite exacte de ${\rm GL}_2(\zp)$-modules topologiques
       $$0\to \Pi_{\delta^{-1}}(\check{D})^\dual \to D^{(0,r_b]}\boxtimes_{\delta}\p1\to \Pi_{\delta}(D)^{(b)}\to 0.$$
     \end{theo}

     Les m\'ethodes utilis\'ees pour l'\'etude de $\Pi^{(b)}$
  sont sensiblement diff\'erentes de celles de~\cite[chap.~V]{Cbigone}:
 les arguments de bidualit\'e de loc.cit. sont remplac\'es
    par une \'etude directe des rayons d'analyticit\'e des vecteurs de $\Pi_{\delta}(D)$, \`a travers l'\'etude de la croissance
    des coefficients de Mahler de $o_v$. Cette
    \'etude est grandement facilit\'ee par la prop.~\ref{relevanal}, qui est aussi utilis\'ee dans la preuve
    du cor.~\ref{uwuintro} ci-dessous. 

       Une cons\'equence imm\'ediate du th.~\ref{mainrayon} est la g\'en\'eralisation suivante du r\'esultat principal du chap.~V de \cite{Cbigone}.

       \begin{coro}
         Si $D\in {\cal C}_L(\delta)$, le sous-faisceau 
$U\mapsto D^\dagger\boxtimes_\delta U$
du faisceau 
\linebreak
$U\mapsto D\boxtimes_\delta U$ est stable par $G$,
qui agit contin\^ument pour la topologie naturelle de~$D^\dagger$,
         et on a suite exacte de $G$-modules topologiques
          $$0\to \Pi_{\delta^{-1}}(\check{D})\to D^{\dagger}\boxtimes_{\delta}\p1\to \Pi_{\delta}(D)^{\rm an}\to 0.$$
       \end{coro}
Mentionnons que
l'action de $G$ sur le faisceau $U\mapsto D^{\dagger}\boxtimes_{\delta}U$ s'\'etend
par continuit\'e en une action sur un faisceau $U\mapsto D_{\rm rig}\boxtimes_{\delta}U$,
et les sections globales
$D_{\rm rig}\boxtimes_{\delta}\p1$ de ce faisceau fournissent (chap.~\ref{drigp1})
une extension de $\Pi_{\delta}(D)^{\rm an}$ par $(\Pi_{\delta^{-1}}(\check D)^{\rm an})^\dual $
qui est tr\`es utile pour l'\'etude de $\Pi_{\delta}(D)^{\rm an}$.

\smallskip
Une autre cons\'equence est le r\'esultat suivant qui renforce le th\'eor\`eme
de Schneider et Teitelbaum sur la densit\'e des vecteurs      
localement analytiques.
\begin{coro}\label{dense1}
 Si $\Pi\in {\rm Rep}_L(G)$ il existe
  $m_0\geq 2$ tel que $\Pi^{(b)}$ soit dense dans~$\Pi^{\rm an}$ 
(et donc aussi dans $\Pi$) pour tout $b> m_0$. 
\end{coro}
En utilisant le fait que les orbites des \'el\'ements de $\Pi^{(b)}$
sont somme de leur s\'erie de Taylor sur $K_{b-1}$, et le cor.~\ref{dense1},
on en d\'eduit le r\'esultat suivant. 
  
  \begin{coro}\label{dense3}
   Soient $\Pi_1, \Pi_2\in {\rm Rep}_L(G)$ et soit 
   $f: \Pi_1^{\rm an}\to \Pi_2^{\rm an}$ une application continue, lin\'eaire et 
   $\mathfrak{gl}_2$-\'equivariante. Alors il existe un sous-groupe ouvert 
compact~$H$ de $G$ tel que $f$ soit $H$-\'equivariante. 
  \end{coro}

\begin{ques}
 Les cor.~\ref{dense1} et \ref{dense3} sont-ils valables
 pour les repr\'esentations de Banach admissibles, topologiquement
de longueur finie,
 d'un groupe de Lie $p$-adique quelconque?
\end{ques}

        Signalons aussi un sous-produit de la preuve, 
pour lequel nous ne connaissons pas de d\'emonstration plus simple. Une telle d\'emonstration
        simplifierait consid\'erablement l'\'etude de $\Pi_{\delta}(D)^{\rm an}$.

        \begin{coro}\label{uwuintro}
         Soient $\Pi\in {\rm Rep}_L(G)$ et $v\in \Pi$. Si les applications
         $x\mapsto \left(\begin{smallmatrix} 1 & x \\0 & 1\end{smallmatrix}\right)v$
   et $x\mapsto \left(\begin{smallmatrix} 1 & 0 \\x & 1\end{smallmatrix}\right)v$ sont localement
   analytiques (de $\qp$ dans $\Pi$), alors $v\in \Pi^{\rm an}$.
        \end{coro}

\makeatletter
\def\thesection{\Roman{section}}
\def\thesubsection{\thesection.\arabic{subsection}}
\def\thesubsubsection{\arabic{subsubsection}}
\makeatother

\section{Anneaux de fonctions analytiques}
\label{Fctan}

Ce chapitre peu \'eclairant introduit un certain nombre d'anneaux de s\'eries de Laurent et \'etablit certains r\'esultats 
 techniques dont on aura besoin dans le chap.~\ref{vectan}. Rappelons que $L$ est une extension finie de $\qp$, 
 dont on note $\O_L$ l'anneau des entiers.
 Pour $b\in\mathbf{N}^\dual $ on note $n_b=p^{b-1}(p-1)$ et $r_b=\frac{1}{n_b}$.

\subsection{Topologies sur les anneaux}\label{top}

   On munit l'anneau $$\oe =\{\sum_{n\in\mathbf{Z}} a_nT^n,\ a_n\in \O_L \ \text{et}\
\lim_{n\to-\infty} v_p(a_n)=\infty\}$$ de la topologie
faible, dont une base de voisinages de $0$ est constitu\'ee
des \linebreak
$p^{n}\oe +T^m\O_L[[T]]$, avec
$m,n\in\mathbf{N}$. On munit son corps des fractions 
\linebreak
${\cal E}=\oe [1/p]=\cup_{n\geq 0} p^{-n}\oe $ de la topologie
limite inductive.

  Si $a\geq b\geq 1$, on note ${\cal E}^{[r_a,r_b]}$ l'anneau des 
$f=\sum_{n\in\mathbf{Z}} a_nT^n$,
 analytiques sur la couronne $r_a\leq v_p(T)\leq r_b$, d\'efinies sur $L$,
   que l'on munit de la valuation $$v^{[r_a,r_b]}(f)=\inf_{r_a\leq v_p(x)\leq r_b} v_p(f(x))
=\inf_{n\in\Z}(v_p(a_n)+\min(nr_a,nr_b)).$$
   On note $\oe^{[r_a,r_b]}$ l'anneau de valuation de ${\cal E}^{[r_a,r_b]}$, et on 
pose
   ${\cal E}^{]0,r_b]}=\linv_{a\geq b} {\cal E}^{[r_a,r_b]}$ (c'est l'anneau
   des fonctions analytiques sur la couronne $0<v_p(T)\leq r_b$, d\'efinies sur~$L$). 
Soit $\oe^{\dagger,b}$ le
compl\'et\'e de $\O_L[[T]][\frac{p}{T^{n_b}}]$ pour la topologie $p$-adique.
La preuve du r\'esultat suivant est laiss\'ee
au lecteur.

 \begin{lemm}\label{oedagger}
{\rm (i)}
  $\oe^{\dagger,b}$ est l'anneau des s\'eries de Laurent $\sum_{k\in\mathbf{Z}} a_kT^k\in  \O_L[[T,T^{-1}]]$ telles que
  la suite\footnote{On note $[\ ]$ la partie enti\`ere.}
 $([v_p(a_k)]+kr_b)_{k\leq 0}$ est positive et tend vers $+\infty$ quand $k\to-\infty$.

{\rm (ii)} $\oe^{[r_a,r_b]}$ est 
l'anneau des s\'eries de Laurent $\sum_{k\in\mathbf{Z}} a_kT^k\in  L[[T,T^{-1}]]$ telles que
  la suite
 $(v_p(a_k)+\min(k r_b,kr_b))_{k\in\Z}$ est positive et tend vers $+\infty$ quand $k\to \pm\infty$.
   \end{lemm}

  On note ${\cal E}^{(0,r_b]}={\cal E}^{]0,r_b]}\cap {\cal E}$
(c'est le sous-anneau de ${\cal E}^{]0,r_b]}$ form\'e des fonctions analytiques born\'ees) et
$\oe^{(0,r_b]}$ le r\'eseau de ${\cal E}^{(0,r_b]}$ form\'e
des s\'eries \`a coefficients dans~$\O_L$. 
 On d\'eduit du lemme \ref{oedagger} que 
$\oe^{(0,r_b]}=\oe^{\dagger,b}[\frac{1}{T}]$ et que
 $\oe^{\dagger,b}$ est s\'epar\'e et complet pour la topologie $T$-adique. Cela munit $\oe^{(0,r_b]}$
d'une topologie naturelle et ${\cal E}^{(0,r_b]}$ de la topologie limite inductive, en \'ecrivant
${\cal E}^{(0,r_b]}=\cup_{k\geq 0} p^{-k}\oe^{(0,r_b]}$.

  Enfin, l'anneau de Robba ${\cal R}$ est la r\'eunion des ${\cal E}^{]0,r_b]}$, muni de la topologie limite inductive
 et le corps ${\cal E}^{\dagger}=\cup_{b\geq 1} {\cal E}^{(0,r_b]}$ est le sous-anneau de ${\cal R}$
   des \'el\'ements born\'es. Il est dense dans ${\cal R}$ et ${\cal E}$ s'identifie au compl\'et\'e de ${\cal E}^\dagger$ pour la valuation $p$-adique.
Si $\Lambda\in \{{\cal E},{\cal R}\}$, on pose $\Lambda^+=\Lambda\cap L[[T]]$.

\subsection{Quelques calculs...}  Les lemmes techniques suivants seront utilis\'es dans l'\'etude
des vecteurs localement analytiques des repr\'esentations unitaires admissibles
de ${\rm GL}_2(\qp)$.

\begin{lemm} \label{intersect}
 On a $p\oe \cap {\cal E}^{\dagger}\subset \cup_{n\geq 1} \oe^{\dagger, n}$.

\end{lemm}

\demo
Soit $f=\sum_{n\in\mathbf{Z}} a_n T^n\in p\oe \cap {\cal E}^{\dagger}$. Il existe 
$b$ tel que $f$ converge sur $0<v_p(T)\leq r_b$. On a donc 
$\lim_{k\to\infty} v_p(a_{-k})-kr_b=\infty$. En particulier, il existe 
$b_1$ tel que si $k\geq n_{b_1}$, alors $v_p(a_{-k})\geq 1+k r_b$. 
Puisque $v_p(a_{-k})\geq 1$ pour tout $k$, on en d\'eduit que 
$[v_p(a_{-k})]\geq k r_{b+b_1}$ pour tout $k\geq 0$, donc $f\in \oe^{\dagger, b+b_1}$ (lemme 
\ref{oedagger}).

   \begin{lemm}\label{valrarb}
    {\rm (i)} $\oe^{[r_a,r_b]}\cap \oe^{(0,r_b]}\subset \frac{1}{p}\oe^{\dagger,b}$.

    {\rm (ii)} Si $f\in \oe^{(0,r_b]}$ satisfait $v^{[r_a,r_b]}(f)\geq N$ pour un $N\in\mathbf{N}$, alors
    $f\in \frac{1}{p} T^{Nn_b} \oe^{\dagger,b}$.
   \end{lemm}

   \demo {\rm (i)} C'est une cons\'equence imm\'ediate du lemme \ref{oedagger}. 
   
   {\rm (ii)} Si $r_a\leq v_p(x)\leq r_b$, on a $$v_p(f(x))-Nn_bv_p(x)\geq v^{[r_a,r_b]}(f)-Nn_br_b\geq 0,$$
   donc $T^{-Nn_b}f\in \oe^{[r_a,r_b]}\cap \oe^{(0,r_b]}$ et on conclut en utilisant le {\rm (i)}.

  \begin{lemm}\label{infrarb}
    {\rm (i)} Si $(f_k)_k$ est une suite d'\'el\'ements de $\oe^{\dagger,b}$ qui converge vers $0$ pour la topologie
    $p$-adique, alors la s\'erie $\sum_{k\geq 0} \big(\tfrac{T^{n_a}}{p}\big)^k f_k$ converge dans
    $\oe^{[r_a,r_b]}$.

    {\rm (ii)} Si $f\in \oe^{[r_a,r_b]}$, alors il existe une suite $(f_k)_k$ comme dans {\rm (i)} et telle que
    $$pf=\sum_{k\geq 0} \big(\tfrac{T^{n_a}}{p}\big)^k f_k.$$

   \end{lemm}

   \demo
    {\rm (i)} Il suffit de constater que $$v^{[r_a,r_b]}\big( \big(\tfrac{T^{n_a}}{p}\big)^kf_k\big)\geq v^{[r_a,r_b]}(f_k)$$
    et que, par hypoth\`ese, la derni\`ere quantit\'e est positive et tend vers $\infty$ pour $k\to\infty$.

    {\rm (ii)} Posons $pf=\sum_{k\in\mathbf{Z}} b_kT^k$ et, pour $k\geq 0$, posons $g_k=\sum_{j=0}^{n_a-1} p^{k} b_{kn_a+j}T^j$.
    Alors $g_k\in \oe^+$ tend vers $0$ pour la topologie $p$-adique (car
     $v_p(b_k)+kr_a\geq 1$ pour tout $k$ et $\lim_{k\to+\infty} v_p(b_k)+kr_a=+\infty$)
     et on a  $$\sum_{k\geq 0} b_kT^k=\sum_{k\geq 0}
    \big(\tfrac{T^{n_a}}{p}\big)^k g_k.$$ Pour conclure, il suffit de v\'erifier que $\sum_{k\leq 0} b_kT^k\in
    \oe^{\dagger,b}$. Cela d\'ecoule du lemme \ref{oedagger}.

\subsection{Actions de $\varphi,\psi,\Gamma$}

   On munit les anneaux $ {\cal E}^+, {\cal R}^+, \oe , {\cal E},
{\cal E}^{\dagger}, {\cal R}$ d'actions continues de $\Gamma={\rm Gal}(\qp(\mu_{p^{\infty}})/\qp)$ et d'un Frobenius $\varphi$, commutant entre elles,
en posant $\varphi(T)=(1+T)^p-1$ et $\sigma_a(T)=(1+T)^a-1$ si $a\in\zpet$. L'op\'erateur $\varphi$ ne laisse pas stable les anneaux
$\oed{b}$, $\oe^{(0,r_b]}$, $\oe^{[r_a,r_b]}$ et ${\cal E}^{]0,r_b]}$;
il les envoie respectivement dans
$\oed{{b+1}}$, $\oe^{(0,r_{b+1}]}$, $\oe^{[r_{a+1},r_{b+1}]}$ et ${\cal E}^{]0,r_{b+1}]}$.
Ces anneaux sont, en revanche, stables sous l'action
de $\Gamma$.

Le corps ${\cal E}$ est une extension de degr\'e~$p$ de
$\varphi({\cal E})$, ce qui permet de d\'efinir un inverse
\`a gauche $\psi$ de $\varphi$ par la
formule
$$\psi(f)=p^{-1}\varphi^{-1}({\rm Tr}_{{\cal E}/\varphi({\cal E})}f).$$
Alors $\psi$ laisse stable $\O_{\cal E}$ et~$\edag$,
 s'\'etend par continuit\'e \`a ${\cal R}$,
et envoie les anneaux
$\oed{{b+1}}$, $\oe^{(0,r_{b+1}]}$, $\oe^{[r_{a+1},r_{b+1}]}$
et~${\cal E}^{]0,r_{b+1}]}$
dans $\oed{b}$, $\oe^{(0,r_b]}$, $\oe^{[r_a,r_b]}$ et ${\cal E}^{]0,r_b]}$
respectivement.
De plus, $\psi$ commute \`a $\Gamma$ et
 $\psi\big(\sum_{i=0}^{p-1}(1+T)^i\varphi(f_i)\big)=f_0$ (tout \'el\'ement
de ${\cal E}$ ou ${\cal R}$ peut s'\'ecrire sous cette forme, et une telle
\'ecriture est unique).

Le r\'esultat suivant est parfaitement classique. 

 \begin{lemm}\label{hyper}

  $\frac{\varphi^{n}(T)}{T^{p^n}}$ est une unit\'e de $\oe^{\dagger,b}$ si $b>n$.

\end{lemm}

\demo Voir le lemme II.5.2 de \cite{CCsurconv}.

\section{$(\varphi,\Gamma)$-modules}
  Ce chapitre est aussi pr\'eliminaire. On rappelle quelques r\'esultats standard de la th\'eorie des 
  $(\varphi,\Gamma)$-modules et on \'etablit deux r\'esultats techniques qui seront utilis\'es dans l'\'etude des vecteurs localement analytiques
  de la repr\'esentation $\Pi_{\delta}(D)$.

\label{surc0}

\subsection{$(\varphi,\Gamma)$-modules et faisceaux $P^+$-\'equivariants sur $\Z_p$}\label{faisceau}
  Soit $A$ un anneau topologique, commutatif,
 muni d'un endomorphisme continu $\varphi$ et d'une
action continue de $\Gamma$, qui commutent. Un
   $(\varphi,\Gamma)$-module sur $A$ est un $A$-module de type fini muni d'un endomorphisme
   semi-lin\'eaire $\varphi$ et d'une action semi-lin\'eaire de~$\Gamma$,
 commutant entre elles.

Un $(\varphi,\Gamma)$-module $D$ sur $\oe $ est dit \textit{\'etale} si $\varphi(D)$ engendre $D$ sur $\oe $.
 Un $(\varphi,\Gamma)$-module $D$ sur ${\cal E}$ est dit \textit{\'etale} s'il admet un $\oe $-r\'eseau stable par
  $\varphi$ et $\Gamma$ et qui est \'etale en tant que $(\varphi,\Gamma)$-module sur $\oe $.

 On note $\fget_{\rm tors}$ (resp. $\fget(\oe )$) la cat\'egorie
 des $(\varphi,\Gamma)$-modules \'etales sur $\oe $, qui sont de torsion (resp. libres) comme
 $\oe $-module. Enfin, on note $\fget({\cal E})$ la cat\'egorie des
 $(\varphi,\Gamma)$-modules \'etales sur ${\cal E}$. 

\medskip
Soit $D$ un $(\varphi,\Gamma)$-module \'etale.
Alors $D $ est muni d'une action de $P^+$ donn\'ee par
$\matrice{p^ka}{b}{0}{1}\cdot z=(1+T)^b\varphi^k\circ\sigma_a(z)$,
si $k\in\N$, $a\in\Z_p^\dual$ et $b\in\Z_p$, et d'un inverse
\`a gauche $\psi$ de $\varphi$ qui commute \`a l'action
de $\Gamma$ et qui est d\'efinie par
$\psi(\sum_{i=0}^{p-1}(1+T)^i\varphi(x_i))=x_0$.
On utilise ces donn\'ees pour associer \`a $D $
un faisceau $U\mapsto D \boxtimes U$ sur $\Z_p$ (o\`u $U$ d\'ecrit
les ouverts compacts de $\Z_p$),  \'equivariant sous l'action
de $P^+$, o\`u $P^+$ agit sur $\Z_p$ par la formule
$\matrice{a}{b}{0}{1}\cdot x={ax+b}$ habituelle.  De mani\`ere pr\'ecise:

$\bullet$ $D \boxtimes\Z_p=D $ et $D \boxtimes\emptyset=0$,

$\bullet$ $D \boxtimes (i+p^k\Z_p)=\matrice{p^k}{i}{0}{1}D \subset D $

$\bullet$ La restriction 
${\rm Res}_{i+p^k\Z_p}:D \boxtimes\Z_p\to D \boxtimes (i+p^k\Z_p)$
est d\'efinie par la formule
${\rm Res}_{i+p^k\Z_p}=\matrice{1}{i}{0}{1}\circ\varphi^k\circ\psi^k
\circ \matrice{1}{-i}{0}{1}$.

\begin{rema}
Soit ${\cal C}$ le faisceau sur $\Z_p$ des fonctions continues \`a valeurs
dans~$L$ et soit ${\cal D}_0$ le faisceau des mesures (i.e le dual de ${\cal C}$).
Le dictionnaire d'analyse fonctionnelle $p$-adique
fournit une suite exacte
$0\to {\cal D}_0\to D \to {\cal C}\boxtimes\chi^{-1}\to 0$
de faisceaux
$P^+$-\'equivariants sur $\Z_p$ si $D={\cal E}$ est le $(\varphi,\Gamma)$-module
trivial (la torsion par~$\chi^{-1}$ signifie que l'action
de $\matrice{a}{b}{0}{1}$ est multipli\'ee par $\chi(a)^{-1}$).
\end{rema}

\subsection{Surconvergence}\label{surc2}

  Soit $D\in\fget(\oe)$.  Si $b\in\N^\dual $, on note $D^{\dagger,b}$
le plus grand sous-$\oed{b}$-module $M$ de type fini
de $D$ tel que $\varphi(M)\subset \oed{b+1}\cdot M$ (le tout \`a
l'int\'erieur de $D$). 
On renvoie \`a \cite[prop. 4.2.6]{BC} pour une preuve du r\'esultat suivant:

 \begin{prop}\label{CCF}
  Si $D\in\fget(\oe )$,
 il existe $m(D)$ tel que $D^{\dagger,m(D)}$ soit libre de rang ${\rm rg}_{\oe }(D)$
  sur $\oe^{\dagger, m(D)}$, et $D^{\dagger,b}=\oe^{\dagger,b}\otimes_{\oe^{\dagger, m(D)}}D^{\dagger, m(D)}$
  pour tout $b\geq m(D)$.

 \end{prop}

   La prop.~\ref{CCF} permet de d\'efinir, pour $a\geq b\geq m(D)$, des modules
$D^{(0,r_b]}$, $D^{[r_a,r_b]}$, $D^{]0,r_b]}$, $D^{\dagger}$ et $D_{\rm rig}$,
en tensorisant $D^{\dagger, m(D)}$ par $\oe^{(0,r_b]}$, $\oe^{[r_a,r_b]}$,
${\cal E}^{]0,r_b]}$, ${\cal E}^{\dagger}$, ${\cal R}$ respectivement.
Ils ne d\'ependent pas du choix de $m(D)$ et ils sont libres de m\^eme rang que $D$ sur les anneaux correspondants.
Le choix d'une base permet de munir ces modules de topologies naturelles (induites par celles  
des anneaux de s\'eries de Laurent, voir le \S~\ref{top}), qui ne d\'ependent pas du choix de la base.

Tous les modules d\'efinis ci-dessus sont munis d'une action
de $\Gamma$, les modules $D^\dagger$ et $D_{\rm rig}$
sont aussi munis d'actions de $\varphi$ et $\psi$ commutant \`a celle
de $\Gamma$ et v\'erifiant $\psi\circ\varphi={\rm id}$
Le sous-faisceau $U\mapsto D^\dagger\boxtimes U$ de $U\mapsto D\boxtimes U$
est donc stable par $P^+$, et il s'\'etend en un faisceau $U\mapsto D_{\rm rig}\boxtimes U$.
 Par contre, $\varphi$ ne pr\'eserve pas les autres modules:
il envoie $D^{\dagger,b}$ dans
$D^{\dagger,b+1}$, $D^{[r_a,r_b]}$ dans $D^{[r_{a+1},r_{b+1}]}$,
et $D^{]0,r_b]}$ dans $D^{]0,r_{b+1}]}$. De mani\`ere analogue,
$\psi$ laisse stable $D^{\dagger}$ et $D_{\rm rig}$, mais il envoie (pour $a\geq b\geq m(D)$)
$D^{\dagger,b+1}$ dans
$D^{\dagger,b}$, $D^{[r_{a+1},r_{b+1}]}$ dans $D^{[r_a,r_{b}]}$,
et $D^{]0,r_{b+1}]}$ dans $D^{]0,r_{b}]}$.

 Nous aurons besoin
de l'estim\'ee plus pr\'ecise ci-dessous.

\begin{lemm} \label{psi}
  Soit $D\in\fget(\oe)$. Il existe $l(D)\geq 1$ tel que,
 pour tous $a\in\mathbf{Z}$, $k\in\mathbf{N}^\dual $ et $b\geq m(D)+k$,
  $$\psi^k(T^a D^{\dagger,b})\subset T^{[\frac{a}{p^{k}}]-l(D)}D^{\dagger, b-k}.$$
\end{lemm}

\demo Si $a\in\mathbf{Z}$ et si $c=[\frac{a}{p^k}]$, 
 le lemme \ref{hyper} montre que
    $$\psi^k(T^a D^{\dagger,b})\subset \psi^k ( T^{p^k c} D^{\dagger,b})=
    \psi^k(\varphi^k(T)^c D^{\dagger,b})= T^c\psi^{k}(D^{\dagger,b}).$$
 On peut donc se contenter de traiter le cas $a=0$.
Fixons une base $e_1,\dots, e_d$ de $D^{\dagger, m(D)}$ sur $\oed{m(D)}$;
c'est aussi une base de $D^{\dagger,b}$ sur $\oed{b}$ pour tout $b\geq m(D)$.

 Soit $l\geq 1$ tel que
$p$ divise $l$ et $\psi((1+T)^{j}e_i)\in T^{-l} D^{\dagger, m(D)}$ pour $(i,j)\in [1,d]\times [0,p-1]$.
Alors $\psi(D^{\dagger,b})\subset T^{-l}D^{\dagger, b-1}$ pour tout $b>m(D)$, car
$\psi(\oe^{\dagger,b})\subset \oe^{\dagger,b-1}$ et donc
   $$D^{\dagger,b}= \sum_{i=1}^{d} \oe^{\dagger,b}\cdot e_i=
   \sum_{i=1}^{d}\sum_{j=0}^{p-1} (1+T)^j\varphi(\oe^{\dagger, b-1})e_i.$$

   Posons $l(D)=2l$ et montrons par r\'ecurrence sur $k$ que
$\psi^k(D^{\dagger,b})\subset T^{-l(D)}D^{\dagger,b-k}$ pour $b\geq m(D)+k$.
   Pour $k=1$, on vient de le faire. Pour passer de $k$ \`a $k+1$, on utilise l'hypoth\`ese de r\'ecurrence et le lemme
   \ref{hyper}, ce qui donne pour $b>m(D)+k$
     $$\psi^{k+1}(D^{\dagger, b})\subset \psi(\varphi(T)^{-\frac{l(D)}{p}}D^{\dagger, b-k})=
 T^{-\frac{l(D)}{p}}\psi(D^{\dagger, b-k}).$$
   On conclut en utilisant l'inclusion $\psi(D^{\dagger, b-k})\subset T^{-\frac{l(D)}{2}} D^{\dagger, b-k-1}$ (second paragraphe)
   et l'in\'egalit\'e $l(D)\geq\frac{l(D)}{2}+\frac{l(D)}{p}$.

\subsection{Dualit\'e}
Le module $\Omega^1_{\oe}$ des $\O_L$-diff\'erentielles continues
de $\oe$ est naturellement un $(\varphi,\Gamma)$-module
\'etale libre de rang $1$, une base \'etant
$\tfrac{dT}{1+T}$ et les actions de $\varphi$ et $\Gamma$ \'etant\footnote{La formule
$\varphi\big(\tfrac{dT}{1+T}\big)=p\tfrac{dT}{1+T}$, qui semblerait naturelle,
ne fournit pas un $(\varphi,\Gamma)$-module \'etale.}
$$\sigma_a\big(\tfrac{dT}{1+T}\big)=a\tfrac{dT}{1+T},\ {\text{ si $a\in\Z_p^\dual$,}}
\quad{\rm et}\quad
\varphi\big(\tfrac{dT}{1+T}\big)=\tfrac{dT}{1+T}.$$
Si $D$ est un objet de $\fget(\oe)$, (resp.~$\fget({\cal E})$,
resp.~$\fget_{\rm tors}$), on note
$\check D$ le $(\varphi,\Gamma)$-module des morphismes $\oe$-lin\'eaires
de $D$ dans $\O_{\cal E}\tfrac{dT}{1+T}$ (resp.~${\cal E}\tfrac{dT}{1+T}$,
resp.~$({\cal E}/\O_{\cal E})\tfrac{dT}{1+T}$),
les actions de $\varphi$ et
$\Gamma$ \'etant d\'efinies par\footnote{La condition {\og $D$ \'etale\fg} est pr\'ecis\'ement ce qu'il faut
pour garantir l'existence et l'unicit\'e d'un tel $\varphi$
sur $\check D$, si $D$ est un $(\varphi,\Gamma)$-module sur $\oe$.}
$$\langle\sigma_a (x),\sigma_a (y)\rangle=\sigma_a (\langle x,y\rangle),\ {\text{ si $a\in\Z_p^\dual$,}}
\quad{\rm et}\quad
\langle\varphi(x),\varphi(y)\rangle=\varphi(\langle x,y\rangle),$$
l'accouplement
$\langle\ ,\ \rangle$ sur $\check{D}\times D$
 \'etant l'accouplement naturel.
Le foncteur $D\to \check D$ est involutif et exact. Par extension des scalaires et fonctorialit\'e, l'accouplement $\langle\ ,\ \rangle$ induit des accouplements (pour 
$a\geq m(D)$) 
$$\langle\ ,\ \rangle:\check D^{(0,r_a]}\times D^{(0,r_a]}
\to{\cal E}^{(0,r_a]}\tfrac{dT}{1+T}, \quad \langle\ ,\ \rangle:\check D^{\dagger,a}\times D^{\dagger,a}\to \oed{a}\tfrac{dT}{1+T}, $$
et $ \langle\ ,\ \rangle:\check D_{\rm rig}\times D_{\rm rig}\to {\cal R}\tfrac{dT}{1+T}$.

  L'application r\'esidu 
  $$\reso: \O_L[[T, T^{-1}]]dT\to \O_L, \quad \reso \big(\big(\sum_{n\in\mathbf{Z}} a_n T^n\big)dT \big)=a_{-1}$$
  induit une application $\reso: \oe \tfrac{dT}{1+T}\to \O_L$ 
 et donc des applications $\reso: {\cal E}\tfrac{dT}{1+T}\to L$ et
 $\reso: {\cal E}/\oe \tfrac{dT}{1+T}\to  L/\O_L$. 
    
    Si $\check{z}\in \check{D}$ et $z\in D$, on pose
$$\{\check{z},z\}=\reso\big(\langle \sigma_{-1}\cdot \check{z}, z \rangle\big).$$
On obtient ainsi un accouplement \`a valeurs dans $L/\O_L$ (resp. $\O_L$, resp. $L$)
si $D\in \fget_{\rm tors}$ (resp. $D\in \fget(\oe )$, resp.
$D\in \fget({\cal E})$). Cet accouplement est parfait, i.e.~l'application $\iota$ qui envoie
$x$ sur $\iota(x)=(y\mapsto \{x,y\})$ identifie $\check{D}$ et $D^\dual $ (le dual \'etant muni de la topologie de la convergence simple). 
On d\'efinit par la m\^eme formule un accouplement parfait
$\{\,\,,\,\}$ entre $\check{D}_{\rm rig}$ et $D_{\rm rig}$.

\medskip

  Le r\'esultat suivant sera utilis\'e
 dans l'\'etude des vecteurs localement analytiques des 
objets de ${\rm Rep}_L(G)$.

  \begin{lemm}\label{residupadic}
   Si $D\in\fget(\oe )$, $a_1,a_2\in\mathbf{Z}$ et $b>\max(m(D),m(\check{D}))$, alors $$\{T^{a_1}\check{D}^{\dagger,b}, T^{a_2}D^{\dagger,b}\}
   \subset \{x\in \O_L, v_p(x)\geq (a_1+a_2)r_b\}.$$
  \end{lemm}

\demo
   Soient $\check{z}\in \check{D}^{\dagger,b}$, $z\in D^{\dagger,b}$. Puisque 
   $\tfrac{\sigma_{-1}(T)}{T}$ est inversible dans $\oe^{\dagger, b}$, il existe 
   $f\in \oe^{\dagger,b}$ tel que 
   $$ \big(\tfrac{\sigma_{-1}(T)}{T}\big)^{a_1} \langle \sigma_{-1}(\check{z}), z\rangle=f \tfrac{dT}{1+T}.$$
   Puisque $\langle \, , \, \rangle$ est $\oe^{\dagger, b}$-lin\'eaire, on a 
     $$\{T^{a_1}\check{z},T^{a_2}z\}=\reso\big(T^{a_1+a_2}f \tfrac{dT}{1+T}\big).$$ 
     En \'ecrivant $f=\sum_{n\in\mathbf{Z}} b_nT^n$, un petit calcul montre
que $$\reso\big( T^{a_1+a_2}f\tfrac{dT}{1+T}\big)=\sum_{j\geq 0} (-1)^j b_{-1-(a_1+a_2)-j},$$
la convergence de la s\'erie \'etant assur\'ee par l'in\'egalit\'e $v_p(b_n)\geq -nr_b$,
si $n\leq 0$. Cette in\'egalit\'e permet aussi de montrer que 
$$ v_p \big(  \reso\big( T^{a_1+a_2}f\tfrac{dT}{1+T}\big) \big)\geq (a_1+a_2)r_b,$$
si $a_1+a_2\geq 0$; le cas $a_1+a_2<0$ \'etant trivial, cela permet de conclure.

\subsection{Les modules $D^{\rm nr}$, $D^{\sharp}$ et $D^{\natural}$}\label{sharp}
Les modules ci-dessous font l'objet d'une \'etude d\'etaill\'ee
   dans \cite[chap.~II]{Cmirab}.

 \begin{defi} {\rm (i)} Si $D\in\fget(\oe )\cup \fget_{\rm tors}$ on note
 $D^{\rm nr}=\cap_{n\geq 1} \varphi^n(D)$ et
   $$D^{++}=\{x\in D| \lim_{n\to\infty} \varphi^n(x)=0\}, \quad
 D^+=D^{++}\oplus D^{\rm nr}.$$

 {\rm (ii)} Si $D\in \fget_{\rm tors}$, on note
$D^{\natural}$ et $D^{\sharp}$ les orthogonaux respectifs de $\check{D}^+$ et $\check{D}^{++}$, pour l'accouplement
$\{\,,\,\}$. Si $D\in \fget (\oe )$, on pose $D^{?}=\varprojlim_{k} (D/p^kD)^{?}$, pour $?\in \{\natural,\sharp\}$.
  \end{defi}

  On \'etend ces d\'efinitions aux $(\varphi,\Gamma)$-modules
  sur ${\cal E}$, en choisissant des r\'eseaux stables par $\varphi$ et $\Gamma$ et en tensorisant par $L$
   (les objets obtenus ne d\'ependent pas
  des choix).

   \begin{prop}\label{dnr}
  Si $D\in\fget(\oe )\cup \fget_{\rm tors}$, alors:

  {\rm (i)} $D^{\rm nr}$ et $D^{\sharp}/D^{\natural}$ sont des $\O_L$-modules de type fini. Si $D$ est de torsion, alors
  $\check{D}^{\rm nr}$ est le dual de $D^{\sharp}/D^{\natural}$.

  {\rm (ii)} $D^{\natural}$ et $D^{\sharp}$ sont des sous $\O_L[[T]]$-modules compacts de $D$, qui engendrent $D$ et sur lesquels
  $\psi$ est surjectif.

  {\rm (iii)}  Si $D$ est de torsion
   ou si $D$ est irr\'eductible de rang $\geq 2$, alors $D^{\sharp}/D^{\natural}$ est
   un $\O_L$-module de longueur finie.

   \end{prop}

\demo
 Toutes les r\'ef\'erences sont \`a \cite{Cmirab}. Le {\rm (i)} suit de la prop. II.2.2 et de la prop. II.5.19.
 Le {\rm (ii)} d\'ecoule de la prop. II.6.3. Enfin, {\rm (iii)} est le cor. II.5.21.

\medskip
 On d\'eduit de la proposition ci-dessus que si $D\in\fget({\cal E})$, alors
$D^{\rm nr}$ et $D^{\sharp}/D^{\natural}$ sont des $L$-espaces vectoriels de dimension finie et
que $\check{D}^{\rm nr}$ est le $L$-dual de $D^{\sharp}/D^{\natural}$. De plus, si $D$ est irr\'eductible
de dimension $\geq 2$, alors $D^{\natural}=D^{\sharp}$ car $\check{D}^{\rm nr}=0$.
 Cela est faux si $D$ est de dimension $1$, car dans ce cas
$D^{\sharp}/D^{\natural}$ est un $L$-espace vectoriel de dimension $1$ 
puisque $\check{D}^{\rm nr}$
est de dimension~$1$.

\section{L'image du foncteur $\Pi\mapsto D(\Pi)$}
 Dans ce chapitre on d\'emontre le th\'eor\`eme~\ref{Repdelta}
  de l'introduction (ainsi que les versions enti\`ere et de torsion de ce th\'eor\`eme).
 Les r\'esultats obtenus n'ont pas d'hypoth\`ese sur
 $p$, car ils n'utilisent pas
 \cite[th.~II.3.3]{Cbigone}, \cite[th.~0.1.1]{KiAst} ou \cite{Pa}.
Beaucoup des arguments qui suivent sont tir\'es de~\cite{Cmirab} et
des chap.~II et~IV de~\cite{Cbigone} mais nous avons explicit\'e
certains r\'esultats implicites dans~\cite{Cbigone}
(comme ceux du \S~\ref{RPi} qui ne sont r\'edig\'es que dans le cas de torsion
dans~\cite{Cbigone}), rajout\'e des sorites sur les invariants
par ${\rm SL}_2(\Q_p)$, simplifi\'e la
d\'emonstration de r\'esultats clefs comme les th.~\ref{dualite}
et~\ref{recoverD},
et introduit la notion de paire $G$-compatible qui rend la pr\'esentation
des r\'esultats plus agr\'eable.

 \subsection{Repr\'esentations de $G$}

 \quad Si $A$ est un anneau commutatif et si $H$ est un groupe topologique, 
une $A$-repr\'esentation de $H$ est un $A[H]$-module \`a gauche.
   Une telle repr\'esentation $\Pi$ est dite \textit{lisse} si
  le stabilisateur de tout $v\in\Pi$ est ouvert dans $H$ et \textit{lisse admissible} si de plus
 $\Pi^K$ est un $A$-module de type fini pour tout sous-groupe ouvert compact $K$ de $H$. 

\smallskip
Nous aurons besoin des cat\'egories suivantes de repr\'esentations 
de $G={\rm GL}_2(\Q_p)$:

$\bullet$  ${\rm Rep}_{\rm tors}(G)$ est la cat\'egorie des $\O_L$-repr\'esentations lisses de $G$, de longueur finie et ayant un
    caract\`ere central\footnote{Qui n'est pas forc\'ement unique.}. 
Tout $\Pi\in {\rm Rep}_{\rm tors}(G)$ est de torsion comme $\O_L$-module, 
et admissible
d'apr\`es les travaux de Barthel-Livn\'e \cite{BL} et Breuil \cite{Br1}. 

$\bullet$ ${\rm Rep}_{\O_L}(G)$ est la cat\'egorie des $\O_L$-repr\'esentations $\Pi$ de $G$, ayant un caract\`ere central et telles que
$\Pi$ est un $\O_L$-module s\'epar\'e et complet pour la topologie $p$-adique
(i.e.~$\Pi=\linv\,\Pi/p^n\Pi$), sans $p$-torsion et tel que
$\Pi/p^n\Pi\in {\rm Rep}_{\rm tors}(G)$ pour tout $n$.

$\bullet$ ${\rm Rep}_L(G)$ est la cat\'egorie des $L$-repr\'esentations 
de Banach de $G$ qui admettent un $\O_L$-r\'eseau
ouvert, born\'e, stable par $G$ et appartenant \`a ${\rm Rep}_{\O_L}(G)$. 
${\rm Rep}_L(G)$ est donc la cat\'egorie des $L$-repr\'esentations de Banach de $G$, qui sont unitaires, admissibles au sens de~\cite{ST1},
   r\'esiduellement de longueur finie\footnote{
La condition {\og r\'esiduellement de longueur finie\fg} implique
{\og topologiquement de longueur finie\fg} de mani\`ere \'evidente.
Paskunas~\cite{Pa} a montr\'e (au moins si $p\geq 5$)
que ces deux conditions sont en fait \'equivalentes.}
 et \`a caract\`ere central. 

\smallskip
  Si $\Pi\in {\rm Rep}_{\rm tors}(G)$ (resp. ${\rm Rep}_{\O_L} (G)$, resp. ${\rm Rep}_L (G)$),
   on note $\Pi^\dual $ le dual de Pontryagin (resp. le $\O_L$ ou $L$-dual continu)
  de $\Pi$, que l'on munit de la topologie faible (i.e. celle de la convergence simple) et de l'action \'evidente de $G$. 

\smallskip
Si $\delta:\Q_p^\dual\to \O_L^\dual$ est un caract\`ere unitaire, on note
${\rm Rep}_{\rm tors}(\delta)$ (resp. ${\rm Rep}_{\O_L} (\delta)$, resp. ${\rm Rep}_L (\delta)$)
la sous-cat\'egorie de
${\rm Rep}_{\rm tors}(G)$ (resp. ${\rm Rep}_{\O_L} (G)$, resp. ${\rm Rep}_L (G)$)
des repr\'esentations sur lesquelles $\matrice{a}{0}{0}{a}$ agit
par multplication par $\delta(a)$.  Par d\'efinition, 
${\rm Rep}_{\rm tors}(G)$ (resp. ${\rm Rep}_{\O_L} (G)$, resp. ${\rm Rep}_L (G)$)
est la r\'eunion des
${\rm Rep}_{\rm tors}(\delta)$ (resp. ${\rm Rep}_{\O_L} (\delta)$, resp. ${\rm Rep}_L (\delta)$),
pour $\delta$ unitaire.

\smallskip
Si $\eta_1,\eta_2$ sont des caract\`eres continus de $\Q_p^\dual$, \`a valeurs
dans $k_L^\dual$ (resp.~$\O_L^\dual$),
on note ${\rm Ind}(\eta_1\otimes\eta_2)$ 
l'espace des fonctions $\phi: G\to k_L$ (resp.~$\phi: G\to L$), continues, telles que
$\phi\left(\left(\begin{smallmatrix} a & b \\0 & d\end{smallmatrix}\right)g\right))=\eta_1(a)
\eta_2(d)f(g)$ pour tous $a,d\in\qpet$, $b\in\qp$ et $g\in G$, que l'on munit
de l'action de $G$ d\'efinie par $(h\cdot\phi)(g)=\phi(gh)$.
Alors ${\rm Ind}(\eta_1\otimes\eta_2)$ est un objet de ${\rm Rep}_{\rm tors}(\eta_1\eta_2)$
(resp.~${\rm Rep}_L(\eta_1\eta_2)$).  Le r\'esultat suivant est parfaitement classique.

\begin{prop}\label{ordi1}
{\rm (i)}
Si $\eta_1\neq\eta_2$, la repr\'esentation
${\rm Ind}(\eta_1\otimes\eta_2)$ est irr\'eductible (resp.~topologiquement
irr\'eductible).

{\rm (ii)} Si $\eta_1=\eta_2$, la fonction $g\mapsto \eta_1\circ\det g$
engendre une sous-repr\'esentation de dimension~$1$ sur laquelle $G$ agit
\`a travers le caract\`ere $\eta_1\circ\det g$, et le quotient
est une repr\'esentation irr\'eductible (resp.~topologiquement
irr\'eductible) de $G$, de la forme ${\rm St}\otimes (\eta_1\circ\det g)$, o\`u
${\rm St}$ est {\emph {la steinberg}} (resp.~{\emph {la steinberg continue}}).
\end{prop}

Les composantes de Jordan-H\"older des ${\rm Ind}(\eta_1\otimes\eta_2)$
sont dites {\it ordinaires}; les objets absolument irr\'eductibles
de ${\rm Rep}_{\rm tors}(G)$ ou ${\rm Rep}_L(G)$ qui ne sont pas ordinaires
sont dits {\it supersinguliers}.
Il n'est pas tr\`es facile de construire des $L$-repr\'esentations
supersinguli\`eres par de purs proc\'ed\'es de th\'eorie des repr\'esentations,
 mais les th.~\ref{princip} et~\ref{LLp} en donnent
une classification
compl\`ete en termes de $(\varphi,\Gamma)$-modules.

\subsection{Le foncteur $\Pi\mapsto D(\Pi)$} \label{Montrealf}
On note $P=\left(\begin{smallmatrix} \qpet & \qp \\0 & 1\end{smallmatrix} \right)$
le sous-groupe {\it mirabolique} de $G$ et
$P^+$ le sous-semi-groupe $\left(\begin{smallmatrix} \zp-\{0\} & \zp \\0 & 1\end{smallmatrix} \right)$ de $P$.

 \quad Soit $\Pi$ un objet de ${\rm Rep}_{\rm tors}(G)$. Si
 $W\subset \Pi$ est un sous-$\O_L$-module de type fini, stable sous l'action de ${\rm GL}_2(\zp)$ et
 qui engendre $\Pi$ comme $G$-module (un tel $W$ existe car $\Pi$ est de longueur finie, cf. \cite[lemme III.1.6]{Cbigone}), on note:

  $\bullet$ $D_{W}^{\natural}(\Pi)$ le dual de Pontryagin de $P^+\cdot W$.

  $\bullet$ $D_{W}^{+}(\Pi)$ l'ensemble des $\mu\in \Pi^\dual $ nuls sur $g\cdot W$
  pour tout $g\in P-P^+$.

   $D_{W}^{+}(\Pi)$ est stable par $P^+$ car $P-P^+$ est stable par multiplication
par $g^{-1}$ si $g\in P^+$; il admet donc une structure naturelle\footnote{Les actions de
   $\varphi$ et $\Gamma$ sont celles de $\left(\begin{smallmatrix} p & 0\\0 & 1\end{smallmatrix} \right)$ et
   $\left(\begin{smallmatrix} \zpet & 0 \\0 & 1\end{smallmatrix} \right)$; la structure de
   $\O_L[[T]]$-module est induite par l'action de $\left(\begin{smallmatrix} 1 & \zp \\0 & 1\end{smallmatrix} \right)$
   et l'isomorphisme standard $\O_L[[T]]\simeq \O_L\left[\left[\left(\begin{smallmatrix} 1 & \zp \\0 & 1\end{smallmatrix} \right)\right]\right]$.} de
   $(\varphi,\Gamma)$-module sur $\O_L[[T]]$. On d\'efinit alors
   $$D(\Pi)=\oe \otimes_{\O_L[[T]]} D_W^{+}(\Pi)$$
   et on montre \cite[th. IV.2.13]{Cbigone} 
que $D(\Pi)\in \fget_{\rm tors}$
(la seule difficult\'e est de v\'erifier que $D(\Pi)$ est de longueur finie). De plus,
   $D(\Pi)$ ne d\'epend pas du choix de $W$ et $\Pi\mapsto D(\Pi)$ est un foncteur
   foncteur exact contravariant
    de ${\rm Rep}_{\rm tors}(G)$ dans $\fget_{\rm tors}$.

     Si $\Pi\in {\rm Rep}_{\O_L}(G)$, on pose
      $$D(\Pi)=\varprojlim_{n} D(\Pi/p^n\Pi).$$
    Enfin, si $\Pi\in {\rm Rep}_L(G)$, on choisit un $\O_L$-r\'eseau ouvert $\Pi_0$, 
 born\'e et stable par $G$, et on pose
     $D(\Pi)=L\otimes_{\O_L} D(\Pi_0)$ (cela ne d\'epend pas du choix de $\Pi_0$). On obtient ainsi des foncteurs
     exacts contravariants ${\rm Rep}_{\O_L}(G)\to \fget(\oe )$ et ${\rm Rep}_L(G)\to \fget({\cal E})$.

\begin{rema}
Il est clair que le foncteur $\Pi\mapsto D(\Pi)$ tue les objets de type fini
sur $\O_L$ (ou $L$) mais, comme on le verra (th.~\ref{princip}),
c'est la seule information que l'on perd en utilisant ce foncteur,
ce qui est assez remarquable car la construction
de $D(\Pi)$ n'utilise que peu d'information sur $\Pi$.
\end{rema}

 \subsection{Le r\'esultat principal}

    Soit $\delta:\qpet\to \O_L^\dual $ un caract\`ere unitaire et soit $D$ un $(\varphi,\Gamma)$-module \'etale\footnote{Cela signifie que $D$
    est un objet d'une des cat\'egories $\fget_{\rm tors}$, $\fget(\oe )$ ou 
    $\fget({\cal E})$.}.
Les constructions de \cite[chap. II]{Cbigone} fournissent un faisceau
$G$-\'equivariant sur $\p1=\p1(\qp)$, dont
l'espace des sections sur $U$ est not\'e $D\boxtimes_{\delta} U$. Par construction,
on a $D\boxtimes_{\delta}\zp=D$ (et la restriction du faisceau \`a $\Z_p$ muni de l'action de $P^+$ est
le faisceau du \S~\ref{faisceau}) et le caract\`ere central du $G$-module
$D\boxtimes_{\delta}\p1$ est~$\delta$.  De plus, si $U$ est un ouvert compact
de $\Q_p$, l'extension par $0$ permet de consid\'erer $D\boxtimes_\delta U$
comme un sous-module de $D\boxtimes_{\delta}\p1$.
Le module $D\boxtimes_\delta U$ est alors stable sous l'action du
stabilisateur de $U$; en particulier, $D=D\boxtimes_\delta\Z_p$ est stable
par $1+p{\bf M}_2(\Z_p)$ puisque ce groupe stabilise $\Z_p\subset\p1$.

Si $w=\matrice{0}{1}{1}{0}$,
   l'application $z\to ({\rm Res}_{\zp}(z), {\rm Res}_{\zp}(wz))$ induit une injection
   de $D\boxtimes_{\delta} \p1$ dans $D\times D$, ce qui permet de munir $D\boxtimes_{\delta}\p1$
   d'une structure de $G$-module topologique ($D$ \'etant muni de la topologie faible).
Plus pr\'ecis\'ement, si on note 
$w_{\delta}$ la restriction de l'action de $w$ \`a 
$D^{\psi=0}=D\boxtimes_{\delta} \zpet$, alors
$D\boxtimes_{\delta} \p1$ s'identifie au sous-ensemble de
$D\times D$ des $(z_1,z_2)$ v\'erifiant ${\rm Res}_{\zpet}z_2=
w_\delta({\rm Res}_{\zpet}z_1)$.
 \begin{rema}\label{squelette}
Comme $G$ est engendr\'e par $\matrice{p}{0}{0}{1}$, $\matrice{\Z_p^\dual}{0}{0}{1}$,
$w$ et $\matrice{1}{1+p\Z_p}{0}{1}$ et comme
$\p1=\Z_p\cup w\cdot p\Z_p$,
    l'action de $G$ sur $D\boxtimes_{\delta}\p1$ est 
compl\`etement d\'ecrite par les formules suivantes.

  $\bullet$ Si $z\in D=D\boxtimes_\delta\Z_p$, si $a\in\zpet$ et si $b\in\zp$, on a
    $$\left(\begin{smallmatrix} p & 0 \\0 & 1\end{smallmatrix}\right)z=\varphi(z), \quad \left(\begin{smallmatrix} a & 0 \\0 & 1\end{smallmatrix}\right)z=
    \sigma_a(z), \quad \left(\begin{smallmatrix} 1 & b \\0 & 1\end{smallmatrix}\right)z=(1+T)^b\cdot z.$$

    $\bullet$ Si $z=(z_1,z_2)\in D\boxtimes_{\delta}\p1$,
 on a $wz=(z_2,z_1)$, ${\rm Res}_{\zp}\left(w\left(\begin{smallmatrix} p & 0 \\0 & 1\end{smallmatrix}\right)z\right)=\delta(p)\psi(z_2)$ et, si $b\in p\zp$,
       ${\rm Res}_{p\zp}\left( w\left(\begin{smallmatrix} 1 & b \\0 & 1\end{smallmatrix}\right)z\right)=u_{b}\left({\rm Res}_{p\zp}(z_2)\right)$,
       o\`{u}\footnote{La formule de \cite[pag. 325]{Cbigone} comporte quelques fautes de frappe.}
       $$u_b=\delta(1+b)\left(\begin{smallmatrix} 1 & -1 \\0 & 1\end{smallmatrix}\right)\circ w_{\delta}\circ
       \left(\begin{smallmatrix} (1+b)^{-2} & b(1+b)^{-1} \\0 & 1\end{smallmatrix}\right)\circ w_{\delta}\circ
       \left(\begin{smallmatrix} 1 & (1+b)^{-1} \\0 & 1\end{smallmatrix}\right)$$ sur $D\boxtimes_\delta p\zp$.

  \end{rema}

         Soit ${\cal C}_{\rm tors}(\delta)\subset \fget_{\rm tors}$ 
l'image de ${\rm Rep}_{\rm tors}(\delta)$ par le foncteur $\Pi\mapsto D(\Pi)$.
     On d\'efinit de mani\`ere analogue les cat\'egories ${\cal C}_{\O_L}(\delta)$ et 
     ${\cal C}_L(\delta)$. 
Il r\'esulte du th.~\ref{princip} ci-dessous
que les objets de ces cat\'egories sont exactement
les $(\varphi,\Gamma)$-modules $D$ tels que $(D,\delta)$ soit
$G$-compatible (cf.~def.~\ref{Pid}).

\begin{theo}\label{princip}
Si $\delta:\qpet\to \O_L^\dual $ est un caract\`ere unitaire, alors:

{\rm (i)} ${\cal C}_{\rm tors}(\delta)$ est stable par sous-quotients.

{\rm (ii)} Si $D\in {\cal C}_{\rm tors}(\delta)$, alors $\check{D}\in {\cal C}_{\rm tors}(\delta^{-1})$.

{\rm (iii)} Il existe un foncteur covariant $D\to \Pi_{\delta} (D)$ de ${\cal C}_{\rm tors}(\delta)$ dans ${\rm Rep}_{\rm tors}(\delta)$
tel que pour tout $D\in {\cal C}_{\rm tors}(\delta)$ on ait une suite exacte de $G$-modules topologiques
 $$0\to \Pi_{\delta^{-1}}(\check{D})^\dual \to D\boxtimes_{\delta}\p1\to \Pi_{\delta}(D)\to 0.$$

{\rm (iv)} Les foncteurs $\Pi\mapsto D(\Pi)$ et $D\to \Pi_{\delta}(\check{D})$ induisent
des anti-\'equivalences quasi-inverses exactes entre ${\rm Rep}_{\rm tors}(\delta)/S$ et ${\cal C}_{\rm tors}(\delta)$,
o\`{u} $S$ est la sous-cat\'egorie de ${\rm Rep}_{\rm tors}(\delta)$ form\'ee des repr\'esentations de type fini comme
$\O_L$-module.

{\rm (v)} Les r\'esultats pr\'ec\'edents restent valables si on remplace ${\cal C}_{\rm tors}(\delta)$ par ${\cal C}_{\O_L}(\delta)$ (resp.
${\cal C}_{L}(\delta)$) et ${\rm Rep}_{\rm tors}(\delta)$ par ${\rm Rep}_{\O_L}(\delta)$ (resp. ${\rm Rep}_L(\delta)$)
et $\O_L$ par $\O_L$ (resp.~$L$) dans la d\'efinition de $S$. 
\end{theo}
\begin{rema}
La suite exacte $0\to \Pi_{\delta^{-1}}(\check{D})^\dual \to 
D\boxtimes_{\delta}\p1\to \Pi_{\delta}(D)\to 0$ ne d\'etermine pas uniquement
$\Pi_{\delta}(D)$ et $\Pi_{\delta^{-1}}(\check{D})$ mais presque (en fait,
si $D\in {\cal C}_{\rm tors}(\delta)$ (resp.~$D\in {\cal C}_L(\delta)$)
n'a pas de sous-quotient isomorphe
\`a $k_{\cal E}(\eta)$ (resp.~${\cal E}(\eta)$),
avec $\delta=\eta^2$ ou $\delta=\eta^2\chi^{-2}$, alors 
$\Pi_{\delta}(D)$ et $\Pi_{\delta^{-1}}(\check{D})$ sont uniquement d\'etermin\'es
par la suite exacte).
Nous donnerons une construction explicite de ces repr\'esentations (cf.~def.~\ref{Pid})
ce qui permet de restaurer l'unicit\'e dans tous les cas.
\end{rema}

  La preuve de ce th\'eor\`eme occupe la quasi-totalit\'e de ce chapitre. 
Les {\rm (i)}, {\rm (ii)} et~{\rm (iii)} s'obtiennent 
  en m\'elangeant le cor.~\ref{checkstab}, les prop.~\ref{sousquot} et~\ref{betap1}, 
ainsi que les th.~\ref{recoverPi} et~\ref{recoverD}.   
  Pour le {\rm (iv)}, il faut en plus utiliser la prop.~\ref{exact}.

\subsection{Paires $G$-compatibles}

 \quad Soit $D$ un $(\varphi,\Gamma)$-module \'etale. L'application $x\to \left({\rm Res}_{\zp} \left(\begin{smallmatrix} p^n & 0 \\0 & 1\end{smallmatrix}\right)x\right)_{n\geq 0}$ induit un isomorphisme
  $$D\boxtimes_{\delta}\qp\cong 
\{(x_n)_{n\in\mathbf{N}}, \ x_n\in D \ \text{et} \ \psi(x_{n+1})=x_n\},$$
 et on munit $D\boxtimes_{\delta}\qp$ de la topologie induite par la topologie produit sur $D^{\mathbf{N}}$. 
\begin{rema}\label{resqp}
Comme on passe de $\p1$ \`a $\qp$ en n'enlevant qu'un point, la restriction
\`a $\Q_p$ est presque injective~\cite[prop.~II.1.14]{Cbigone}:
$${\rm Ker}\big({\rm Res}_{\Q_p}:D\boxtimes_{\delta}\p1\to D\boxtimes_{\delta}\qp\big)=
\{(0,z_2),\ z_2\in D^{\rm nr}\}.$$
\end{rema}
Si $?\in \{\natural, \sharp\}$, on pose
 $$D^{?}\boxtimes_{\delta}\qp=(D\boxtimes_{\delta}\qp)\cap (D^{?})^{\mathbf{N}}.$$ 
Si $D\in \fget(\oe )\cup \fget_{\rm tors}$,
c'est un sous-module compact de $D\boxtimes_{\delta}\qp$.
\begin{prop}\label{stabilite1}
Soit $D\in\fget_{\rm tors}\cup \fget(\oe )$.

{\rm (i)} Si
 $M$ est un sous $\O_L$-module ferm\'e de $D\boxtimes_{\delta} \qp$, stable par $P$,
il existe un sous-objet $D_1$ de $D$ tel que 
$$D_1^{\natural}\boxtimes_{\delta}\qp\subset M\subset D_1^{\sharp}\boxtimes_{\delta}\qp.$$
En particulier, $M\subset D^{\sharp}\boxtimes_{\delta}\qp$
et $D^{\natural}\boxtimes_{\delta}\qp\subset M$,
si ${\rm Res}_{\zp} (M)$ engendre $D$ en tant que $(\varphi,\Gamma)$-module.

{\rm (ii)} $(D^{\sharp}\boxtimes_{\delta}\qp)/(D^{\natural}\boxtimes_{\delta}\qp)$
est isomorphe \`a $D^\natural/D^\sharp$ et est de type fini sur $\O_L$.

{\rm (iii)} Le foncteur $D\mapsto D^{\sharp}\boxtimes_{\delta}\qp$ est exact.
\end{prop}

\demo
 Le (i) correspond au th.~III.3.8 de \cite{Cmirab} (noter que le caract\`ere
 $\delta$ ne joue aucun r\^{o}le quand on consid\`ere la restriction \`a $P$). 
Le (ii) correspond aux prop.~III.3.1 et cor.~III.3.2 de~\cite{Cmirab},
et le (iii) au th.~III.3.5 de~\cite{Cmirab}.

\medskip
On d\'efinit des sous-$B$-modules 
[ferm\'es si $D\in\fget_{\rm tors}\cup \fget(\oe )$]
 de $D\boxtimes_{\delta}\p1$ par:
 $$D^{\sharp}\boxtimes_{\delta}\p1={\rm Res}_{\qp}^{-1}(D^{\sharp}\boxtimes_{\delta}\qp), \quad (D^{\natural}\boxtimes_{\delta}\p1)_{\rm ns}=
{\rm Res}_{\qp}^{-1}(D^{\natural}\boxtimes_{\delta}\qp).$$

 On pose $D^{\natural}\boxtimes_{\delta}\p1=(D^{\natural}\boxtimes_{\delta}\p1)_{\rm ns}$
  si $D\in \fget_{\rm tors}\cup \fget({\cal E})$, et\footnote{Le sous-module $(D^{\natural}\boxtimes_{\delta}\p1)_{\rm ns}$
  de $D\boxtimes_{\delta}\p1$ n'est pas forc\'ement satur\'e $p$-adiquement, voir la rem.~VII.4.28 de \cite{Cbigone}.} on d\'efinit
$D^{\natural}\boxtimes_{\delta}\p1$ comme
  le satur\'e du $\O_L$-module $(D^{\natural}\boxtimes_{\delta}\p1)_{\rm ns}$ si
  $D\in \fget(\oe )$.

  \begin{rema}\label{fonct1}
$(D^{\sharp}\boxtimes_{\delta}\p1)/(D^{\natural}\boxtimes_{\delta}\p1)$
et $(D^{\natural}\boxtimes_{\delta}\p1)/(D^{\natural}\boxtimes_{\delta}\p1)_{\rm ns}$
sont
de type fini (sur $\O_L$ ou $L$ suivant les cas) et sont nuls si
$D^\natural=D^\sharp$.  En effet,
$D^\sharp\boxtimes_\delta\qp$ est satur\'e, et donc
$D^{\sharp}\boxtimes_{\delta}\p1$ aussi, ce qui fait que
$D^{\natural}\boxtimes_{\delta}\p1\subset D^{\sharp}\boxtimes_{\delta}\p1$
et que $(D^{\sharp}\boxtimes_{\delta}\p1)/(D^{\natural}\boxtimes_{\delta}\p1)$
est un quotient de 
$(D^{\sharp}\boxtimes_{\delta}\p1)/(D^{\natural}\boxtimes_{\delta}\p1)_{\rm ns}$.
Or, par d\'efinition, ${\rm Res}_{\Q_p}$ induit une injection
de $(D^{\sharp}\boxtimes_{\delta}\p1)/(D^{\natural}\boxtimes_{\delta}\p1)_{\rm ns}$
dans $(D^{\sharp}\boxtimes_{\delta}\qp)/(D^{\natural}\boxtimes_{\delta}\qp)\cong D^\sharp/D^\natural$
(cf.~(ii) de la prop.~\ref{stabilite1}).
\end{rema}

\begin{defi}\label{Pid}
 On dit que $(D,\delta)$ est \textit{$G$-compatible} si $D^{\natural}\boxtimes_{\delta}\p1$
  est stable par~$G$. Dans ce cas, on pose $$\Pi_{\delta}(D)=(D\boxtimes_{\delta}\p1)/(D^{\natural}\boxtimes_{\delta}\p1).$$
 \end{defi}
\begin{rema}\label{fonct2}
  Si $f: D_1\to D_2$ est un morphisme de $(\varphi,\Gamma)$-modules, 
  $f$ induit un morphisme \'equivariant du faisceau attach\'e \`a
  $(D_1,\delta)$ dans le faisceau attach\'e \`a $(D_2, \delta)$
(cela r\'esulte de la construction du faisceau $D\to D\boxtimes_{\delta} U$).
 En particulier, $f$ induit des morphismes
  de $G$ (resp. $B$)-modules topologiques $f: D_1\boxtimes_{\delta}\p1\to D_2\boxtimes_{\delta}\p1$
  (resp. $f: D_1\boxtimes_{\delta}\qp\to D_2\boxtimes_{\delta}\qp$). Si $?\in\{\sharp,\natural\}$, alors $f$ envoie
  $D_1^{?}$ dans $D_2^{?}$ et donc $f$ envoie $D_1^{?}\boxtimes_{\delta}\qp$ dans $D_2^{?}\boxtimes_{\delta}\qp$
  et $D_1^{?}\boxtimes_{\delta}\p1$ dans $D_2^{?}\boxtimes_{\delta}\p1$.
Il en r\'esulte que si $(D_1,\delta)$ et $(D_2,\delta)$ sont $G$-compatibles,
alors $f$ induit un morphisme $G$-\'equivariant de $\Pi_\delta(D_1)$
dans $\Pi_\delta(D_2)$
  \end{rema}

\begin{prop}
Si $(D,\delta)$ est $G$-compatible,
alors $\Pi_\delta(D)$ est un objet de $\rm{Rep}_{\rm{tors}}(G)$,
$\rm {Rep}_{\O_L}(G)$ ou $\rm {Rep}_L(G)$, suivant
que $D\in \fget_{\rm tors}$, $D\in \fget(\oe )$,
 ou $D\in \fget({\cal E})$.
\end{prop}
\demo
Cf.~\cite[lemme II.2.10]{Cbigone}: la seule difficult\'e est de prouver
que les objets obtenus sont (r\'esiduellement) de longueur finie (voir
la prop.~\ref{longueur} pour une justification de cette finitude).

\begin{rema}\label{dimfinPi}
 Soit $(D,\delta)$ une paire $G$-compatible.

{\rm (i)} Si $D\in\fget_{\rm tors}\cup \fget(\oe )$,
alors $D^{?}\boxtimes_{\delta}\p1$ 
est compact, si $?\in\{\natural,\sharp\}$.
En effet,
 $z\mapsto({\rm Res}_{\Z_p}z,{\rm Res}_{\Z_p}w\cdot z)$
permet d'identifier $D^\sharp\boxtimes_\delta\p1$ \`a un sous-module ferm\'e de $D^\sharp\times D^\sharp$,
ce qui prouve qu'il est compact.  Le m\^eme argument montre que
$(D^{\natural}\boxtimes_{\delta}\p1)_{\rm ns}$ est compact, et la rem.~\ref{fonct1}
permet d'en d\'eduire le r\'esultat pour $D^{\natural}\boxtimes_{\delta}\p1$.

{\rm (ii)} Si $D\in\fget_{\rm tors}\cup \fget(\oe )$
est non nul, alors
$\Pi_{\delta}(D)$ n'est pas de type fini comme $\O_L$-module.
En effet $D^{\natural}\boxtimes_{\delta}\p1$ est compact
et donc son intersection $M$
avec $D=D\boxtimes_{\delta}\Z_p$ aussi, 
ainsi que l'image $\overline M$ de $M$ dans $k_L\otimes D$. Il en r\'esulte que
$(k_L\otimes D)/\overline M$ est de dimension infinie sur $k_L$ et donc que l'image de $D$
dans $\Pi_{\delta}(D)$ n'est pas de type fini sur $\O_L$.

{\rm (iii)} Si $D\in\fget({\cal E})$ est non nul, alors
$\Pi_{\delta}(D)$ est de dimension infinie sur $L$ (cela r\'esulte du (ii)
en tensorisant par $L$).
\end{rema}

\begin{prop}\label{twist}
 Soit $D$ un $(\varphi,\Gamma)$-module \'etale et soient $\delta,\eta:\qpet\to \O_L^\dual $ des caract\`eres
 unitaires. Si $(D,\delta)$ est une paire $G$-compatible, il en est de m\^eme de $(D(\eta),\delta \eta^2)$
 et on a un isomorphisme de $G$-modules de Banach $$\Pi(D(\eta),\delta\eta^2)\simeq \Pi(D,\delta)\otimes (\eta\circ\det).$$
\end{prop}

\demo
 C'est une cons\'equence de l'isomorphisme (cf.~\cite[prop. II.1.11]{Cbigone})
$D(\eta)\boxtimes_{\delta\eta^2}\p1\cong
(D\boxtimes_\delta\p1)\otimes (\eta\circ\det)$.

\begin{prop}\label{dim1} 
Si $\Lambda\in\{k_{\cal E},{\cal E})$ et si $D$ est de rang $1$ sur $\Lambda$, alors
$(D,\delta)$ est $G$-compatible pour tout $\delta$.
Plus pr\'ecis\'ement, si $\delta_1,\delta_2$ sont deux caract\`eres
unitaires:

{\rm (i)} On a un isomorphisme de $G$-modules topologiques
$$\Lambda(\delta_1)^{\natural}\boxtimes_{\delta_1\delta_2\chi^{-1}}\p1
\cong
\big({\rm Ind}(\delta_1\otimes \delta_2\chi^{-1})\big)^\dual \otimes (\delta_1\delta_2\chi^{-1}\circ\det).$$

{\rm (ii)} L'application 
$z\mapsto \phi_z$, avec $\phi_z(g)={\rm res}_0\big({\rm Res}_{\zp}(wgz)\tfrac{dT}{1+T}\big)$,
induit un isomorphisme de $G$-modules topologiques
$$\Pi_{\delta_1\delta_2\chi^{-1}}(\Lambda(\delta_1))\cong
{\rm Ind}(\delta_2\otimes \delta_1\chi^{-1}).$$
\end{prop}

\demo
Il s'agit d'une traduction de l'analyse fonctionnelle sur $\zp$: 
voir la rem.~II.1.1 de \cite{Cbigone} ou la prop.~4.9 de~\cite{Cvectan}.

\begin{theo}\label{LLp}
{\rm (i)} Si\footnote{Un des ingr\'edients de la d\'emonstration de
la $G$-compatibilit\'e de $(D,\delta_D)$ est la zariski-densit\'e de l'ensemble
des $(\varphi,\Gamma)$-modules cristallins dans l'espace des
$(\varphi,\Gamma)$-modules.  Si $p=2$, cette densit\'e n'est pas connue
dans le cas ou le semi-simplifi\'e de la r\'eduction modulo~$p$ de $D$ est
la somme directe de deux caract\`eres \'egaux. Les r\'esultats d'unicit\'e
reposent sur les travaux de Paskunas qui supposent $p\geq 5$.  Il y a donc quelques
restrictions provisoires \`a la validit\'e de l'\'enonc\'e.}
$D$ est de rang~$2$ et si $\delta_D$ est le caract\`ere $\chi^{-1}\,\det D$,
alors $(D,\delta_D)$ est $G$-compatible et si $D$ est ind\'ecomposable,
$\delta_D$ est l'unique caract\`ere~$\delta$
de $\Q_p^\dual$ tel que $(D,\delta)$ soit $G$-compatible.

{\rm (ii)} Si $D$ est absolument irr\'eductible de rang~$\geq 3$, alors
$(D,\delta)$ n'est $G$-compatible pour aucun choix de $\delta$.
\end{theo}
\demo
La $G$-compatibilit\'e de $(D,\delta_D)$ est le r\'esultat principal
du chap.~IV de~\cite{Cbigone}.
Le reste de l'\'enonc\'e est une cons\'equence des travaux de Paskunas~\cite{Pa}.
\begin{rema}
La $G$-compatibilit\'e de $(D,\delta_D)$
est valable, plus g\'en\'eralement, pour une d\'eformation d'un
$(\varphi,\Gamma)$-module de rang~$2$, i.e.~pour un $(\varphi,\Gamma)$-module
de rang~$2$ sur $A\otimes_L{\cal E}$ o\`u $A$ est une $L$-alg\`ebre artinienne.
Ce genre de consid\'eration joue d'ailleurs un grand r\^ole dans la d\'emonstration
de la $G$-compatibilit\'e de $(D,\delta_D)$ pour $D\in\fget({\cal E})$.
\end{rema}

\begin{prop}\label{stabilite2}
{\rm (i)} Soit $D\in\fget(\oe )$.
Alors $(D,\delta)$ est $G$-compatible si et seulement si
  $(D/p^kD, \delta)$ est  $G$-compatible pour tout $k\geq 0$.

{\rm (ii)}
 Une paire $(D,\delta)$ est $G$-compatible si et seulement si 
$D^{\sharp}\boxtimes_{\delta}\p1$ est stable par~$G$. 

{\rm (iii)} Si $D\in \fget_{\rm tors}\cup \fget(\oe )$,
et si $(D,\delta)$ est $G$-compatible,
le module $D^{\sharp}\boxtimes_{\delta}\p1$ est le plus grand sous-module compact de
$D\boxtimes_{\delta}\p1$ stable par $G$.
\end{prop}

\demo
Ces \'enonc\'es sont contenus dans la prop.~II.2.6 de~\cite{Cbigone}
(et sa preuve).

\subsection{Le sous-module $\tilde{D}^+$ de $D^{\natural}\boxtimes_{\delta}\p1$}\label{tildeplus}

   On note $\tilde{\mathbf{E}}_{\qp}$ le compl\'et\'e de la cl\^{o}ture radicielle de $\mathbf{F}_p((T))$.
   Il est muni d'actions continues de $\varphi$ et $\Gamma$ (on a $\varphi(T)=T^p$ et $\sigma_a(T)=(1+T)^a-1$).
   Soit $\tilde{\mathbf{A}}_{\qp}=W(\tilde{\mathbf{E}}_{\qp})$ l'anneau des vecteurs de Witt \`a coefficients dans
   $\tilde{\mathbf{E}}_{\qp}$. Si $x\in \tilde{\mathbf{E}}_{\qp}$, on note $[x]$ le repr\'esentant de Teichmuller
   de $x$ dans $\tilde{\mathbf{A}}_{\qp}$. L'anneau $\tilde{\mathbf{A}}_{\qp}$ est naturellement muni d'actions de $\varphi$ et $\Gamma$, que l'on \'etend
   par $\O_L$-lin\'earit\'e \`a l'anneau $\tilde{\O}_{{\cal E}}:=\O_L\otimes_{\zp} \tilde{\mathbf{A}}_{\qp}$.
    Si $b\in \qp$ et $n\geq 1$ est tel que $p^n b\in \zp$, on pose
     $$[(1+T)^b]=\varphi^{-n} ((1+T)^{p^n b})=\varphi^{-n} 
\big(\sum_{k=0}^{\infty} \tbinom{p^n b}{k} T^k\big)\in \tilde{\O}_{{\cal E}}.$$
     Cela ne d\'epend pas du choix de $n$ et on a $[(1+T)^{b+c}]=[(1+T)^b]\cdot [(1+T)^c]$ si $b,c\in \qp$.

   Si $D$ est
   un $(\varphi,\Gamma)$-module \'etale, on pose $\tilde{D}=\tilde{\O}_{{\cal E}}\otimes_{\oe } D$, que l'on munit
   d'une action du mirabolique $P$
  en posant, si $k\in\mathbf{Z}$, $a\in\zpet$, $b\in\qp$,
   $$\left(\begin{smallmatrix} p^ka & b \\0 & 1\end{smallmatrix}\right) \tilde{z}=[(1+T)^{b}]\varphi^{k}(\sigma_a(\tilde{z})).$$
 Cette action laisse stable le sous-module $\tilde{D}^+$ de $\tilde{D}$, form\'e des $x\in \tilde{D}$ tels que la suite
 $(\varphi^n(x))_{n}$ soit born\'ee dans $\tilde D$.
\begin{prop}\label{td1}
Si $D\in\fget_{\rm tors}\cup\fget({\cal E})$, alors
$\tilde D/\tilde D^+$ est un $\O_L[B]$-module (resp.~un
$L[B]$-module topologique) de longueur \'egale \`a celle
de $D$.  En particulier, si $D$ est irr\'eductible, il en est
de m\^eme de $\tilde D/\tilde D^+$ comme $B$-module (topologique).
\end{prop}
\demo
C'est une reformulation de~\cite[prop.~IV.5.6]{Cmirab} et de sa preuve.
\begin{rema}\label{td2}
La d\'emonstration de~\cite[prop.~IV.5.6]{Cmirab} repose sur le fait
que $\tilde D/\tilde D^+$ est le dual de $\check D^\natural\boxtimes\Q_p$
(cf.~\cite[prop.~IV.5.4]{Cmirab}), ce qui permet d'utiliser
le (i) de la prop.~\ref{stabilite1} (i.e.~\cite[th.~III.3.8]{Cmirab})
pour d\'eterminer la longueur de $\tilde D/\tilde D^+$.
\end{rema}

    Soit $I_n=[0,1[\cap p^{-n}\zp$ et soit $I$ la r\'eunion croissante des $I_n$. C'est
     un syst\`eme de repr\'esentants de $\qp/\zp$ et
    on montre \cite[lemme IV.1.2]{Cmirab} que tout \'el\'ement $z$ de~$\tilde{D}$ 
s'\'ecrit, de mani\`ere unique, sous la forme $z=\sum_{i\in I} [(1+T)^i]z_i$, avec $z_i\in D$ et
     $\lim_{i\to \infty} z_i=0$.
 D'apr\`es \cite[lemme II.1.16]{Cbigone},
 la suite de terme g\'en\'eral $$
 \sum_{i\in I_n} \left(\begin{smallmatrix} 1 & i \\0 & 1\end{smallmatrix}\right) z_i
 \in D\boxtimes_{\delta} p^{-n}\zp\subset D\boxtimes_{\delta}\p1$$ converge dans $D\boxtimes_{\delta}\p1$
 et on note $i(z)$ sa limite. On montre \cite[lemme IV.2.2]{Cmirab} que $i: \tilde{D}\to D\boxtimes_{\delta}\p1$ est
 une injection $P$-\'equivariante, qui envoie $\tilde{D}^+$ dans $(D^{\natural}\boxtimes_{\delta}\p1)_{\rm ns}$
car ${\rm Res}_{\Z_p}\tilde D^+\subset D^\natural$.
(Tout ceci ne suppose pas que $(D,\delta)$ soit $G$-compatible.)

\subsection{Dualit\'e}

\quad  Si $D$ est un $(\varphi,\Gamma)$-module \'etale, on
\'etend  \cite[th.~II.1.13]{Cbigone} l'accouplement $\{\,\,,\,\}$ sur $ \check{D}\times D$ en un accouplement $G$-\'equivariant et parfait
  $\{\,\,,\,\}_{\p1}$ sur $ (\check{D}\boxtimes_{\delta^{-1}}\p1)\times (D\boxtimes_{\delta}\p1),$ en posant
    $$\{(\check{z}_1, \check{z}_2), (z_1,z_2)\}_{\p1}=\{\check{z}_1, z_1\}+\{\psi(\check{z}_2), \psi(z_2)\}.$$
Dans cet accouplement, $\check{D}\boxtimes_{\delta^{-1}}U$ et $D\boxtimes_{\delta} V$ sont orthogonaux
si $U$ et $V$ sont des ouverts compacts de $\p1$ tels que $U\cap V=\emptyset$
(on se ram\`ene \`a~\cite[prop.~III.2.3]{Cmirab} en utilisant la $G$-\'equivariance).

  \begin{lemm}\label{orthoincl}
   Soit $D$ un $(\varphi,\Gamma)$-module \'etale et $\delta:\qpet\to \O_L^\dual $ un caract\`ere unitaire.
   Alors l'orthogonal de $\tilde{D}^+$ dans $\check{D}\boxtimes_{\delta^{-1}}\p1$ est
inclus dans $\check{D}^{\natural}\boxtimes_{\delta^{-1}}\p1$.
  \end{lemm}

\demo
 Le cas $D\in\fget({\cal E})$ se d\'eduit par tensorisation par $L$;
 on suppose donc que
$D\in\fget_{\rm tors}\cup\fget(\oe )$.
 Soit $N$ l'orthogonal de $\tilde{D}^+$ dans $\check{D}\boxtimes_{\delta^{-1}}\p1$. Il est stable par
$P$, car $\tilde{D}^+$ l'est.

Supposons que $D$ est de torsion. Si $x=(x_1,x_2)\in N$, alors pour tout $y\in D^+\subset
\tilde{D}^+$ on a $\{x_1,y\}=\{x,y\}_{\p1}=0$, donc $x_1$ est orthogonal \`a $D^+$ et $x_1\in \check{D}^{\natural}$. En appliquant
ceci \`a $\left(\begin{smallmatrix} p^n & 0\\0 & 1\end{smallmatrix}\right)x$ pour tout $n\geq 0$, on obtient $x\in \check{D}^{\natural}\boxtimes_{\delta^{-1}}\p1$ et
donc $N\subset \check{D}^{\natural}\boxtimes_{\delta^{-1}}\p1$.

 Supposons maintenant que $D\in\fget(\oe )$ et soient $D_k=D/p^k D$ et $c$ tel que
 $p^c$ tue le conoyau de $\tilde{D}^+\to \tilde{D}_k^+$ pour tout $k$ 
(l'existence de $c$ d\'ecoule de \cite[lemme~IV.5.1]{Cmirab}). 
Donc, si $x\in N$, alors $p^cx\pmod {p^k}$ est orthogonal
 \`a $\tilde{D}_k^+$ et, d'apr\`es le cas de torsion, on a $p^c x\pmod {p^k}\in
 \check{D}_k^{\natural}\boxtimes_{\delta^{-1}}\p1$. En passant \`a la limite projective
 on obtient 
\linebreak
$p^cx\in (\check{D}^{\natural}\boxtimes_{\delta^{-1}}\p1)_{\rm ns}$ et donc
 $x\in\check{D}^{\natural}\boxtimes_{\delta^{-1}}\p1$, ce qui permet de conclure.

\begin{theo}\label{dualite}  Si $(D,\delta)$ est une paire $G$-compatible, alors $\check{D}^{\natural}\boxtimes_{\delta^{-1}}\p1$ est l'orthogonal
de $\tilde{D}^+$ et de $D^{\natural}\boxtimes_{\delta}\p1$ dans $\check{D}\boxtimes_{\delta^{-1}}\p1$.

\end{theo}

\demo
 En utilisant le lemme~\ref{orthoincl}
 et l'inclusion $\tilde{D}^+\subset D^{\natural}\boxtimes_{\delta}\p1$, il suffit de
 montrer que $D^{\natural}\boxtimes_{\delta}\p1$ est orthogonal \`a $\check{D}^{\natural}\boxtimes_{\delta^{-1}}\p1$. Quitte \`a remplacer
  $D$ par $D/p^kD$ et \`a passer \`a la limite, on peut supposer
que $D$ est de torsion.

  Soit $M$ l'orthogonal de $D^{\natural}\boxtimes_{\delta}\p1$. Alors $M$ est ferm\'e
(c'est un orthogonal)
dans $\check{D}^{\natural}\boxtimes_{\delta^{-1}}\p1$ (lemme~\ref{orthoincl}), 
et donc $M$ est compact (rem.~\ref{dimfinPi}). 
Il s'ensuit que ${\rm Res}_{\qp}(M)$ est
  un sous $P$-module compact de $\check{D}^{\natural}\boxtimes_{\delta}\qp$, 
et comme ${\rm Res}_{\zp}(M)$ contient $\check{D}^{++}$ (car $M$ contient
  $\check{D}^{++}\subset \check{D}^{\natural}\boxtimes_{\delta^{-1}}\p1$), qui engendre $\check{D}$,
  on en d\'eduit (cf.~(i) de la prop.~\ref{stabilite1}) 
que ${\rm Res}_{\qp}(M)=\check{D}^{\natural}\boxtimes_{\delta}\qp$ et donc 
(rem.~\ref{resqp}) $\check{D}^{\natural}\boxtimes_{\delta^{-1}}\p1\subset
  M+(0, \check{D}^{\rm nr})$. Ainsi, $\check{D}^{\natural}\boxtimes_{\delta^{-1}}\p1$ 
est lui-m\^eme compact\footnote{Cela n'a rien de trivial
  \`a cet instant, car nous ne savons pas encore que $(\check{D},\delta^{-1})$ est $G$-compatible. C'est d'ailleurs ce que nous cherchons
  \`a montrer...} car
$(0, \check{D}^{\rm nr})$ est de longueur finie sur $\O_L$ (prop.~\ref{dnr}).
 Le m\^eme argument montre alors que ${\rm Res}_{\qp}(\check{M})=D^{\natural}\boxtimes_{\delta}\qp$,
  si $\check{M}$ est l'orthogonal de $\check{D}^{\natural}\boxtimes_{\delta^{-1}}\p1$.

  On a donc
   $\check{D}^{\natural}\boxtimes_{\delta}\p1\subset
  M+(0, \check{D}^{\rm nr})$ et $D^{\natural}\boxtimes_{\delta}\p1\subset \check{M}+(0, D^{\rm nr})$, et il reste \`a voir que $M+(0, \check{D}^{\rm nr})$ est orthogonal \`a $\check{M}+(0, D^{\rm nr})$. Or, on a vu que $\check{M}+(0,D^{\rm nr})\subset D^{\natural}\boxtimes_{\delta}\p1$ et, par d\'efinition, $M$ est orthogonal \`a $D^{\natural}\boxtimes_{\delta}\p1$, donc
  $M$ est orthogonal \`a $\check{M}+(0,D^{\rm nr})$. En faisant la m\^eme chose avec $\check{M}$ et en utilisant le fait que
  $D^{\rm nr}$ est orthogonal \`a $\check{D}^{\rm nr}$
car $\check{D}^{\rm nr}\subset \check D^+$ et $D^{\rm nr}\subset D^\natural$,
 cela permet de conclure.

\begin{coro}\label{checkstab}
 Soit $(D,\delta)$ une paire $G$-compatible. Alors:
 
 {\rm (i)} $(\check{D}, \delta^{-1})$ est $G$-compatible.
 
 {\rm (ii)} 
On a un isomorphisme
 de $G$-modules topologiques 
 $\Pi_{\delta}(D)^\dual \simeq \check{D}^{\natural}\boxtimes_{\delta^{-1}}\p1$
et $\tilde{\check{D}}^+$ est dense\footnote{Rappelons que tous les duaux sont munis de la topologie faible.} dans $\Pi_{\delta}(D)^\dual $.
 
 {\rm (iii)} On a une suite exacte de $G$-modules topologiques $$0\to \Pi_{\delta^{-1}}(\check{D})^\dual \to D\boxtimes_{\delta}\p1\to
 \Pi_{\delta}(D)\to 0.$$

\end{coro}

\demo
 Cela d\'ecoule du th\'eor\`eme pr\'ec\'edent et du fait que $\{\,,\,\}_{\p1}$ est 
$G$-\'equivariant (pour le (i)) et parfait (pour le reste).

\subsection{Un mod\`ele de $\Pi_\delta(D)$}

L'application ${\rm Res}_{\qp}: D\boxtimes_{\delta}\p1\to D\boxtimes_{\delta}\qp$
 est $B$-\'equivariante et son noyau
est inclus dans $(D^{\natural}\boxtimes_{\delta}\p1)_{\rm ns}$ par d\'efinition de ce module.
Si $D\in \fget({\cal E})$ et si $?\in \{\natural, \sharp\}$, on pose
  $(D^{?}\boxtimes_{\delta}\qp)_{\rm b}=(D_0^{?}\boxtimes_{\delta}\qp)\otimes_{\O_L} L$, pour n'importe quel r\'eseau $D_0$ de $D$, stable par $\varphi$ et
$\Gamma$.

\begin{prop}\label{resqp1} Soit $(D,\delta)$ une paire $G$-compatible.
  Posons $X=D^{\natural}\boxtimes_{\delta}\qp$ si $D\in\fget_{\rm tors}\cup \fget(\oe )$ et $X=(D^{\natural}\boxtimes_{\delta}\qp)_{\rm b}$
  si $D\in\fget({\cal E})$. Alors
   ${\rm Res}_{\qp}$ induit une suite exacte  $$0\to (0,D^{\rm nr})\to (D^{\natural}\boxtimes_{\delta}\p1)_{\rm ns}\to
  X\to 0.$$

\end{prop}

\demo
  L'exactitude \`a gauche est imm\'ediate, celle au milieu
r\'esulte de la rem.~\ref{resqp}.
 Par d\'efinition,  ${\rm Res}_{\qp}((D^{\natural}\boxtimes_{\delta}\p1)_{\rm ns})\subset X$.
 Pour d\'emontrer la surjectivit\'e, on peut supposer que $D\in\fget_{\rm tors}\cup \fget(\oe )$.
 Alors $(D^{\natural}\boxtimes_{\delta}\p1)_{\rm ns}$ est compact (rem.~\ref{dimfinPi})
 donc son image par ${\rm Res}_{\qp}$ est un sous-$P(\qp)$-module compact de $D^{\natural}\boxtimes_{\delta}\qp$,
 qui contient $\tilde{D}^+$; le (i) de la prop.~\ref{stabilite1}
permet donc de montrer que cette image contient $D^{\natural}\boxtimes_{\delta}\qp$.

\begin{coro}\label{bisom}
  Soit $(D,\delta)$ une paire $G$-compatible telle que $\check{D}^{\rm nr}=0$.

  {\rm (i)} Si $D\in\fget_{\rm tors}$, on a un isomorphisme
$\Pi_{\delta}(D)^\dual \simeq \check{D}^{\natural}\boxtimes_{\delta^{-1}}\qp$ de $B$-modules compacts.

  {\rm (ii)} Si $D\in\fget({\cal E})$, on a un isomorphisme 
  $\Pi_{\delta}(D)^\dual \simeq (\check{D}^{\natural}\boxtimes_{\delta^{-1}}\qp)_{\rm b}$ de $B$-modules topologiques.
\end{coro}

\demo
 Cela d\'ecoule du cor.~\ref{checkstab} et de la prop.~\ref{resqp1}.

\begin{coro}\label{lienKir}
 Si $(D,\delta)$ est $G$-compatible, avec $D\in \fget_{\rm tors}\cup
 \fget({\cal E})$, alors l'inclusion de $\tilde{D}$ dans $D\boxtimes_{\delta}\p1$ induit une suite exacte de $B$-modules topologiques
  $$0\to \tilde{D}/\tilde{D}^+\to \Pi_{\delta}(D)\to D^{\sharp}/D^{\natural}\to 0.$$
\end{coro}

\demo
  Commen\c{c}ons par le cas de torsion. Alors $D^{\sharp}/D^{\natural}$ est le dual (de Pontryagin) de $\check{D}^{\rm nr}$ et
  $\tilde{D}/\tilde{D}^+$ est le dual de $\check{D}^{\natural}\boxtimes_{\delta^{-1}}\qp$
(rem.~\ref{td2}).
  En utilisant le th.~\ref{dualite}, on voit que la suite exacte demand\'ee est obtenue en dualisant la suite
  exacte $$0\to (0,\check{D}^{\rm nr})\to \check{D}^{\natural}\boxtimes_{\delta^{-1}}\p1\to \check{D}^{\natural}\boxtimes_{\delta}\qp\to 0$$
  de la prop.~\ref{resqp1}, ce qui permet de conclure.

    Supposons que $D\in \fget(\oe )$ et posons $D_k=D/p^k D$. Alors
    $D^{\sharp}/D^{\natural}$ est la limite projective des $D_k^{\sharp}/D_k^{\natural}$ et
    $\tilde{D}/\tilde{D}^+$ est la limite projective des $\tilde{D}_k/\tilde{D}_k^+$ (\cite[lemme IV.5.3]{Cmirab}
    pour ce dernier). L'application naturelle $\tilde{D}_{k+1}\to \tilde{D}_k$ est surjective, car $D_{k+1}$ se surjecte sur $D_k$. Il en est donc de
    m\^eme de l'application $\tilde{D}_{k+1}/\tilde{D}_{k+1}^+\to \tilde{D}_k/\tilde{D}_{k}^+$.
    Ainsi, en passant \`a la limite dans $$0\to \tilde{D}_k/\tilde{D}_k^+
    \to \Pi_{\delta}(D_k)\to D_k^{\sharp}/D_k^{\natural}\to 0$$ on obtient bien une
     suite exacte
     $$0\to \tilde{D}/\tilde{D}^+\to \varprojlim \Pi_{\delta}(D_k)\to D^{\sharp}/D^{\natural}\to 0.$$
On a $\Pi_\delta(L\otimes D)=L\otimes\Pi_\delta(D)$ puisque
 $\Pi_{\delta}(D)$ est le quotient de $\varprojlim \Pi_{\delta}(D_k)$ par son $\O_L$-module de torsion;
on en d\'eduit le r\'esultat pour un objet de
$\fget({\cal E})$.

\begin{rema}
Dans le cas $D\in \fget(\oe )$, il r\'esulte de la preuve
ci-dessus qu'il faut modifier la suite exacte de la proposition
en rempla\c cant $\Pi_{\delta}(D)$ par $\varprojlim \Pi_{\delta}(D_k)$
ou bien en quotientant $D^\sharp/D^\natural$ par son $\O_L$-module de torsion.
\end{rema}

\begin{prop}\label{longueur}
 Si $(D,\delta)$ est $G$-compatible, avec $D\in 
\fget_{\rm tors}\cup \fget({\cal E})$, alors
 $\Pi_{\delta}(D)$ est un $B$-module (topologiquement) de longueur finie,
et donc aussi un 
\linebreak 
$G$-module (topologiquement) de longueur finie.
\end{prop}
\demo
C'est une cons\'equence du cor.~\ref{lienKir}, de la finitude
de la longueur
de $\tilde D/\tilde D^+$ (prop.~\ref{td1}) et de celle
de $D^\sharp/D^\natural$ (prop.~\ref{dnr}).

\subsection{Presque exactitude de $D\to \Pi_{\delta}(D)$}

  \quad Si $(D_1,\delta)$ et $(D_2,\delta)$ sont des paires $G$-compatibles (noter que $\delta$ est le m\^eme dans les deux paires)
  et si $f: D_1\to D_2$ est un morphisme de $(\varphi,\Gamma)$-modules, la rem.~\ref{fonct2} montre que
  $f$ induit un morphisme $G$-\'equivariant continu $f: \Pi_{\delta}(D_1)\to \Pi_{\delta}(D_2)$. 
  
  \begin{prop}\label{sousquot}
  Soit $0\to D_1\to D\to D_2\to 0$ une suite exacte dans une des cat\'egories
 $\fget_{\rm tors}$, $\fget(\oe )$,
 $\fget({\cal E})$. Si $(D,\delta)$ est une paire $G$-compatible, alors $(D_1,\delta)$ et $(D_2,\delta)$ sont des paires $G$-compatibles.

  \end{prop}
  
  \demo
Pour montrer que $(D_1,\delta)$ est $G$-compatible,
il suffit (prop. \ref{stabilite2}) 
de montrer que\footnote{Si $D\in\fget({\cal E})$, remplacer
 $D^{\sharp}\boxtimes_{\delta}\qp$ par $(D^{\sharp}\boxtimes_{\delta}\qp)_{\rm b}$.}
 $D_1^{\sharp}\boxtimes_{\delta}\p1$ est stable par
 $G$, ce qui r\'esulte de ce que $D_1^{\sharp}\boxtimes_{\delta}\p1=
(D_1\boxtimes_{\delta}\p1)\cap (D^{\sharp}\boxtimes_{\delta}\p1)$
par exactitude du foncteur $D\to D^{\sharp}\boxtimes_{\delta}\qp$ (prop.~\ref{stabilite1}).

   Pour montrer que $(D_2,\delta)$ est $G$-compatible, on dualise la suite exacte
   $0\to D_1\to D\to D_2\to 0$ et on obtient une suite exacte $0\to \check{D}_2\to \check{D}\to \check{D}_1\to 0$. On conclut alors
   en utilisant ce que l'on vient de d\'emontrer et le cor.~\ref{checkstab}.

\begin{rema}
 La r\'eciproque la prop.~\ref{sousquot} est presque toujours fausse, la 
$G$-compatibilit\'e \'etant une contrainte tr\`es forte.
\end{rema}

  \begin{prop}\label{exact}
 Soit $0\to D_1\to D\to D_2\to 0$ une suite exacte dans 
 $\fget_{\rm tors}$ (resp. $\fget(\oe )$, resp.
 $\fget({\cal E})$). Si $(D,\delta)$ est une paire $G$-compatible, les groupes de cohomologie
 du complexe $0\to \Pi_{\delta}(D_1)\to\Pi_{\delta}(D)\to\Pi_{\delta}(D_2)\to 0$ sont des $\O_L$-modules de longueur finie
 (resp. de type fini sur $\O_L$, resp. de dimension finie sur $L$). 
\end{prop}
\demo  
Commen\c cons par traiter le cas de $(\varphi,\Gamma)$-modules sur $\oe$.
La suite $0\to D_1\boxtimes_{\delta}\p1\to D\boxtimes_{\delta}\p1\to
D_2\boxtimes_{\delta}\p1\to 0$ \'etant trivialement exacte,
il suffit de prouver que
la cohomologie du complexe $0\to D_1^{\natural}\boxtimes_{\delta}\p1\to
D^{\natural}\boxtimes_{\delta}\p1\to D_2^{\natural}\boxtimes_{\delta}\p1\to 0$ 
a les propri\'et\'es de finitude requises. 
L'exactitude du foncteur $D\to D^{\sharp}\boxtimes_{\delta}\qp$ et le fait que 
$\Delta^{\sharp}\boxtimes_{\delta}\qp/ \Delta^{\natural}\boxtimes_{\delta}\qp$ soit
un sous-quotient de $\Delta^{\sharp}/\Delta^{\natural}$, et donc un
$\O_L$-module de type fini si $\Delta\in\{D_1,D,D_2\}$,
 entrainent la finitude des groupes de cohomologie du complexe
$0\to D_1^{\natural}\boxtimes_{\delta}\qp\to D^{\natural}\boxtimes_{\delta}\qp\to D_2^{\natural}\boxtimes_{\delta}\qp\to 0$. 
Les suites exactes (pour $\Delta\in\{D_1,D,D_2\}$)
$$0\to (0, \Delta^{\rm nr})\to (\Delta^{\natural}\boxtimes_{\delta}\p1)_{\rm ns}\to \Delta^{\natural}\boxtimes_{\delta}\qp\to 0$$
fournies par la prop.~\ref{resqp1} et la finitude des $\O_L$-modules $\Delta^{\rm nr}$ montrent que les groupes de cohomologie
du complexe $0\to  (D_1^{\natural}\boxtimes_{\delta}\p1)_{\rm ns}\to  (D^{\natural}\boxtimes_{\delta}\p1)_{\rm ns}\to  (D_2^{\natural}\boxtimes_{\delta}\p1)_{\rm ns}\to 0$
sont de type fini sur $\O_L$. 

$\bullet$ Si $D_1,D,D_2$ sont des objets de $\fget_{\rm tors}$, cela permet de conclure. 

$\bullet$ Si 
$D_1,D,D_2$ sont des objets de $\fget(\oe )$, 
on conclut en utilisant la finitude de
$\Delta^{\natural}\boxtimes_{\delta}\p1/(\Delta^{\natural}\boxtimes_{\delta}\p1)_{\rm ns}$ 
(rem.~\ref{fonct1}).

$\bullet$ Si $D_1,D,D_2$ sont des objets de $\fget({{\cal E}})$,
il suffit de remplacer dans la d\'emonstration
$\Delta^{?}\boxtimes_{\delta}\p1$ par $(\Delta^{?}\boxtimes_{\delta}\p1)_{\rm b}$,
 si $\Delta=D_1,D,D_2$ et $?\in\{\natural,\sharp\}$.

\begin{rema}\label{exactnesspi}  
Si $D_1,D,D_2$ sont des objets de $\fget_{\rm tors}$ ou de
$\fget({{\cal E}})$,
il r\'esulte de la preuve que
la suite $0\to \Pi_{\delta}(D_1)\to \Pi_{\delta}(D)\to \Pi_{\delta}(D_2)\to 0$
est exacte si $D_j^{\rm nr}=0$ et
$\check{D}_j^{\rm nr}=0$ pour $j=1,2$.  En effet, dans ce cas
$\Delta^{\rm nr}$ et $\Delta^\natural/\Delta^\natural\cong (\check\Delta^{\rm nr})^\dual $
sont nuls si $\Delta\in\{D_1,D,D_2\}$, et donc les groupes
de cohomologie du complexe sont nuls.
\end{rema}

\begin{prop}\label{irred}
 Soit $(D,\delta)$ une paire $G$-compatible, avec 
$D\in\fget_{\rm tors}\cup\fget({\cal E})$.

{\rm (i)} Si $\Pi_\delta(D)$ est irr\'eductible, alors $D$ est irr\'eductible.

{\rm (ii)} Si $D$ est de dimension $\geq 2$, les assertions suivantes sont
 \'equivalentes:

\quad{\rm a)} $D$ est irr\'eductible.

\quad{\rm b)} $\Pi_{\delta}(D)$ est topologiquement irr\'eductible comme $G$-module.

\quad{\rm c)} $\Pi_{\delta}(D)$ est topologiquement irr\'eductible comme $B$-module, $B$ \'etant le Borel sup\'erieur.

\end{prop}

\demo 
(i)
 Si  $0\to D_1\to D\to D_2\to 0$ est une suite exacte dans $\fget({\cal E})$,
  la prop.~\ref{exact} et l'irr\'eductibilit\'e de $\Pi_{\delta}(D)$ montrent qu'une des repr\'esentations $\Pi_{\delta}(D_1)$ et 
  $\Pi_{\delta}(D_2)$ est de dimension finie sur $L$. La rem.~\ref{dimfinPi} 
permet d'en d\'eduire que $D_1=0$ ou $D_2=0$, et donc 
que $D$ est irr\'eductible. 

(ii)  
Supposons que $D$ est irr\'eductible et montrons le c). 
Comme $\dim_{{\cal E}} D\geq 2$, on a
$D^\sharp=D^\natural$ et donc
$\Pi_{\delta}(D)\simeq \tilde D/\tilde D^+$ 
en tant que $B$-modules topologique (cor.~\ref{lienKir}).
On conclut en utilisant la prop.~\ref{td1}.
L'implication b)$\Rightarrow$a) ayant \'et\'e prouv\'ee dans le (i), cela
permet de conclure puisque l'implication c)$\Rightarrow$b) est triviale.

\begin{rema}
Les conclusions de la proposition sont en d\'efaut en rang~$1$.

(i) La repr\'esentation $\Pi_{\delta}(D)$ n'est jamais 
topologiquement irr\'eductible comme 
\linebreak
$B$-module: il y a un quotient de
dimension~$1$ puisque la repr\'esentation obtenue est une induite d'un
caract\`ere de $B$ (prop.~\ref{dim1}).

(ii) Si $D={\cal E}(\eta)$, et si $\delta=\eta^2$,
alors $\Pi_{\delta}(D)$ n'est pas topologiquement irr\'eductible comme $G$-module
(il y a un sous-objet de dimension~$1$).
Par contre, si $\delta\neq\eta^2$, alors 
$\Pi_{\delta}(D)$ est topologiquement irr\'eductible.
\end{rema}

\begin{rema}
Si $(D,\delta)$ est $G$-compatible avec $D\in\fget_{\rm tors}\cup
\fget({\cal E})$, et si $0\to D_1\to D\to D_2\to 0$ est exacte,
les groupes de cohomologie de $0\to \Pi_\delta(D_1)\to
\Pi_\delta(D)\to\Pi_\delta(D_2)\to0$ sont des $G$-modules
dont les composantes de Jordan-H\"older sont parmi
celles de $D\boxtimes_\delta\p1$.
Comme, par ailleurs (prop.~\ref{exact}), elles sont de type fini (sur $\O_L$ ou $L$),
la suite $0\to \Pi_\delta(D_1)\to
\Pi_\delta(D)\to\Pi_\delta(D_2)\to0$
est exacte si $D\boxtimes_\delta\p1$ n'a pas de composante
de Jordan-H\"older de dimension finie (sur $k_L$ ou sur $L$).
En regroupant les r\'esultats des prop.~\ref{ordi1}, \ref{dim1}, \ref{irred}
et du cor.~\ref{checkstab}, on voit que c'est le cas si et seulement
si $D$ n'a pas de composante de Jordan-H\"older de la forme
$k_{\cal E}(\eta)$ ou $L(\eta)$, avec $\delta=\eta^2$ ou $\delta=\eta^2\chi^{-2}$.
\end{rema}

\subsection{Invariants sous ${\rm SL}_2(\qp)$}

  \quad Dans ce \S~on \'etudie
les ${\rm SL}_2(\qp)$-invariants d'une repr\'esentation de ${\rm Rep}_{\rm tors}(G)$ ou
${\rm Rep}_L(G)$.

 \begin{lemm}\label{finitors}
   Si $\Pi\in {\rm Rep}_{\rm tors}(G)$, alors 
$\Pi^{{\rm SL}_2(\qp)}$ est un $\O_L$-module de longueur
   finie.
  \end{lemm}

\demo Soit $K_m=1+p^m{\rm M}_2(\Z_p)$.
 Comme $\Pi$ est admissible,
il suffit de montrer que $\Pi^{{\rm SL}_2(\qp)}\subset \Pi^{K_m}$ pour $m$ assez grand.
 Soit $\delta$ un caract\`ere central de $\Pi$ et soit $n\geq 1$ tel que $\delta$ soit trivial
 sur $1+p^{n}\zp$. Si $x\in 1+p^{n+1}\zp$, il existe $y\in 1+p^n\zp$ tel que $x=y^2$. Si
 $v\in \Pi^{{\rm SL}_2(\qp)}$, alors $$\left(\begin{smallmatrix} x & 0 \\0 & 1\end{smallmatrix} \right)v=
 \left(\begin{smallmatrix} y & 0 \\0 & y\end{smallmatrix} \right)\left(\begin{smallmatrix} y & 0 \\0 & y^{-1}\end{smallmatrix} \right)v=
 \delta(y)v=v.$$ Comme $K_{n+1}\subset 
\matrice{1+p^{n+1}\zp}{0}{0}{1}\cdot {\rm SL}_2(\qp)$,
 on voit que l'on peut prendre $m=n+1$.

\begin{coro}\label{dimfinsl}
 Si $\Pi\in {\rm Rep}_{\O_L} (G)$ (resp.
 $\Pi\in {\rm Rep}_L (G)$), alors 
$\Pi^{{\rm SL}_2(\qp)}$ est un $\O_L$-module libre de type fini
 (resp. un $L$-espace vectoriel de dimension finie).
\end{coro}
\demo
 On peut supposer que $\Pi\in {\rm Rep}_{\O_L} (G)$. Dans ce cas, le r\'esultat 
est une cons\'equence du lemme pr\'ec\'edent
 et du fait que $\Pi^{{\rm SL}_2(\qp)}$ est un $\O_L$-module satur\'e
de $\Pi$.

\begin{rema}\label{finitors2}
La d\'emonstration ci-dessus utilise juste l'admissibilit\'e.
Or on a suppos\'e que les repr\'esentations sont de
longueur finie et toute composante de Jordan-H\"older
de $\Pi^{{\rm SL}_2(\qp)}$ est de dimension au plus~$2$
sur $k_L$ ou $L$ car le sous-groupe engendr\'e
par ${{\rm SL}_2(\qp)}$ et le centre est d'indice~$2$ dans $G$.
Cela prouve, non seulement que $\Pi^{{\rm SL}_2(\qp)}$ a les propri\'et\'es
de finitude du lemme~\ref{finitors} et du cor.~\ref{dimfinsl},
mais aussi que $(\Pi^\dual)^{{\rm SL}_2(\qp)}$ a les m\^emes propri\'et\'es.
\end{rema}

  \begin{prop}\label{isom}
Soit $\Pi\in{\rm Rep}_{\O_L}(G)$.  Si $\Pi^{{\rm SL}_2(\qp)}=0$, 
le $\O_L$-module
$((L/\O_L)\otimes\Pi)^{{\rm SL}_2(\qp)}$ est de type fini,
et donc inclus dans la $p^n$-torsion de $(L/\O_L)\otimes\Pi$ si $n$ est assez grand.
\end{prop}
\demo
Notons $H$ le groupe ${\rm SL}_2(\qp)$
et $\Pi_n$ la $p^n$-torsion de $(L/\O_L)\otimes\Pi$.
Alors $\Pi_n\cong\Pi/p^n\Pi$ est un objet de ${\rm Rep}_{\rm tors}(G)$
et donc $\Pi_n^H$ est de type fini sur $\O_L$ (lemme \ref{finitors}).
Il s'agit de prouver qu'il existe $n_0\in\N$ tel que
$((L/\O_L)\otimes\Pi)^H=\Pi_{n_0}^H$.
Dans le cas contraire,
il existe une partie infinie $I$ de $\N$ et, pour tout $n\in I$, un vecteur
$v_n\in \Pi_{n}^H$ n'appartenant pas \`a $\Pi_{n-1}$.
On peut donc trouver
$x_n\in \Pi-p\Pi$ tel que $v_n= p^{-n}x_n\pmod \Pi$ et $gx_n-x_n\in p^{n}\Pi$ pour tout
$g\in H$. Pour $n\in I\cap ]j,\infty]$ on a $x_n\pmod {p^j}\in \Pi_j^H$, 
qui est un ensemble fini
(lemme \ref{finitors}). Par extraction diagonale, on obtient ainsi l'existence d'une sous-suite
$(y_n)_n$ de $(x_n)_{n\in I}$ qui converge $p$-adiquement vers un $\alpha\in\Pi$.
En passant \`a la limite dans l'\'egalit\'e $gx_n-x_n\in p^{n}\Pi$, on obtient
$\alpha\in \Pi^H=0$. Mais cela contredit le fait que $y_n\notin p\Pi$ pour tout $n$, 
ce qui permet de conclure.

\begin{lemm}\label{sl2tors}
 Soit $M$ un $\O_L$-module tu\'e par une puissance de $p$ et muni
d'une action $\O_L$-lin\'eaire de ${\rm SL}_2(\qp)$. Alors
$M/M^{{\rm SL}_2(\qp)}$ n'a pas de ${\rm SL}_2(\qp)$-invariants non triviaux.
\end{lemm}

\demo
 Il faut montrer que si $x\in M$ et $(g-1)(h-1)x=0$ pour tous $g,h\in {\rm SL}_2(\qp)$, alors
 $(g-1)x=0$ pour tout $g\in {\rm SL}_2(\qp)$. 
Pour $h=g^n$, on obtient $g^{n+1}(x)-g^{n}(x)=g(x)-x$, et donc
$g^{n}(x)=n(g(x)-x)+x$ 
 pour tous $g\in {\rm SL}_2(\qp)$ et $n\geq 0$. Par hypoth\`ese
 il existe $n$ qui tue $M$. On a alors $g^{n}(x)=x$ pour tout $g\in {\rm SL}_2(\qp)$. On conclut
 en utilisant le fait que $g\to g^n$ est bijective sur les sous-groupes unipotents de~$G$,
  sous-groupes qui engendrent ${\rm SL}_2(\qp)$.

\begin{coro}\label{sl2tors1}
 Si $\Pi\in {\rm Rep}_L(G)$, alors $\Pi/\Pi^{{\rm SL}_2(\qp)}$ n'a pas
 de ${\rm SL}_2(\qp)$-invariants non triviaux.
\end{coro}

\demo
C'est une cons\'equence formelle du lemme \ref{sl2tors}.

\begin{lemm}\label{unitdimfin}
 {\rm (i)} Si $\Pi\in {\rm Rep}_{\rm tors}(G)\cup {\rm Rep}_{\O_L}(G)$ est
 un $\O_L$-module de type fini, alors $\Pi=\Pi^{{\rm SL}_2(\qp)}$.
  
 {\rm (ii)} Si $\Pi\in {\rm Rep}_L(G)$ est de dimension finie sur $L$, alors 
 $\Pi=\Pi^{{\rm SL}_2(\qp)}$.

\end{lemm}

\demo
  Si $\Pi\in {\rm Rep}_{\rm tors}(G)$, cela d\'ecoule de 
   \cite[lemme III.1.5]{Cbigone}. Les autres cas s'en d\'eduisent.

\begin{prop}\label{plusgrande}
 Si $(D,\delta)$ est $G$-compatible, avec $D\in\fget_{\rm tors}\cup \fget(\oe )$ (resp. 
 $D\in\fget({\cal E})$), alors
  $$\Pi_{\delta}(D)^{{\rm SL}_2(\qp)}=(D^{\sharp}\boxtimes_{\delta}\p1)/(D^{\natural}\boxtimes_{\delta}\p1)$$ et c'est la plus grande sous-repr\'esentation
  de type fini sur $\O_L$ (resp. de dimension finie sur $L$) de $\Pi_{\delta}(D)$.
\end{prop}

\demo On peut supposer que $D\in\fget_{\rm tors}\cup \fget(\oe )$.
 D'apr\`es la prop.~\ref{stabilite2},
$(D^{\sharp}\boxtimes_{\delta}\p1)/(D^{\natural}\boxtimes_{\delta}\p1)$ est un sous-$\O_L[G]$-module de $\Pi_{\delta}(D)$,
et il est de type fini sur $\O_L$ d'apr\`es la rem.~\ref{fonct1}.
On d\'eduit du lemme \ref{unitdimfin} que $X\subset \Pi_{\delta}(D)^{{\rm SL}_2(\qp)}$.

   Pour montrer l'inclusion inverse, soit
 $Y$ le sous-$\O_L$-module des $z\in D\boxtimes_{\delta}\p1$ dont l'image dans $\Pi_{\delta}(D)$ est invariante par ${\rm SL}_2(\qp)$.
  Si $z\in Y$, alors $(\left(\begin{smallmatrix} 1 & 1 \\0 & 1\end{smallmatrix}\right)-1)x\in D^{\natural}\boxtimes_{\delta}\p1$, pour $x\in\{z,wz\}$, et en appliquant 
  ${\rm Res}_{\zp}$ on obtient 
${\rm Res}_{\zp}(z), {\rm Res}_{\zp}(wz)\in \frac{1}{T}D^{\natural}$
car ${\rm Res}_{\zp}\big(\big(\matrice{1}{1}{0}{1}-1\big)x\big)=T\,{\rm Res}_{\zp}x$. 
   Donc $Y$ est compact. Comme $Y$ est stable par $G$, on obtient $Y\subset D^{\sharp}\boxtimes_{\delta}\p1$
d'apr\`es le (ii) de la prop.~\ref{stabilite2}, ce qui permet de conclure.

\subsection{Reconstruction de $\Pi$}\label{RPi}
Le but de ce \S~est d'expliquer comment reconstruire $\Pi$ \`a partir de
$D(\Pi)$ (\'eventuellement \`a des morceaux pr\`es de type fini sur $\O_L$ ou~$L$).

 On reprend les notations du  \S~\ref{Montrealf}. 
Soit $\Pi\in {\rm Rep}_{\rm tors}(\delta)$ et soit $W$ comme dans loc.cit.
La restriction \`a $P^+\cdot W$ induit une injection
   de $D_W^+(\Pi)$ dans $D_{W}^{\natural}(\Pi)$, dont l'image est d'indice fini
   dans $D_{W}^{\natural}(\Pi)$ (cf.~\cite[lemme~IV.1.4]{Cbigone}). On a donc un isomorphisme de
   $\oe $-modules $D(\Pi)\simeq \oe \otimes_{\O_L[[T]]} D_W^{\natural}(\Pi)$, ce qui permet
   de d\'efinir une application $\beta_{\zp}: \Pi^\dual \to D(\Pi)$, compos\'ee des
    $$\Pi^\dual \to D_W^{\natural}(\Pi)\to \oe \otimes_{\O_L[[T]]} D_{W}^{\natural}(\Pi)\simeq D(\Pi),$$
    la premi\`ere fl\`eche \'etant la restriction \`a $P^+\cdot W$, et la deuxi\`eme $x\to 1\otimes x$. On d\'efinit $$\beta_{\p1}: \Pi^\dual \to D(\Pi)\oplus D(\Pi), \quad
    \beta_{\p1}(x)=(\beta_{\zp}(x), \beta_{\zp}(w\cdot x)).$$
    Comme leurs noms l'indiquent, 
$\beta_{\zp}$ et $\beta_{\p1}$ ne d\'ependent pas
    du choix de $W$ et sont fonctorielles par fonctorialit\'e
de $\Pi\mapsto D(\Pi)$.

   Soient maintenant $\Pi\in {\rm Rep}_{\O_L} (\delta)$ et $D=D(\Pi)$.
On note $\Pi_n$ le sous-module de $p^n$-torsion de $(L/\O_L)\otimes\Pi$
et on pose
   $D_n=D(\Pi_n)\simeq D/p^n$. Alors $\Pi^\dual=\linv\,\Pi_n^\dual$
et $D=\linv\,D_n$.
      Les applications $\beta_{\p1}: \Pi_n^\dual \to D_n\oplus D_n$ sont compatibles, d'o\`u
   une application continue $\beta_{\p1}: \Pi^\dual \to D\oplus D$. Le cas 
   $\Pi\in {\rm Rep}_{L} (\delta)$ s'en d\'eduit en prenant un r\'eseau appartenant \`a 
   ${\rm Rep}_{\O_L}(\delta)$.

   \begin{prop}\label{betap1}

    Soit $\Pi\in {\rm Rep}_{\rm tors}(\delta)\cup {\rm Rep}_{\O_L}(\delta)\cup {\rm Rep}_L(\delta)$ et soit $D=D(\Pi)$.

    {\rm (i)} $(D,\delta^{-1})$ est une paire $G$-compatible.

  {\rm (ii)} $\beta_{\p1}$ est un morphisme $G$-\'equivariant $\Pi^\dual \to D^{\sharp}\boxtimes_{\delta^{-1}}\p1$, de noyau $(\Pi^\dual )^{{\rm SL}_2(\qp)}$.
  
    {\rm (iii)} $\beta_{\p1}$ envoie
 l'orthogonal de $\Pi^{{\rm SL}_2(\qp)}$ dans $D^{\natural}\boxtimes_{\delta^{-1}}\p1$.

   \end{prop}
   
   \demo
   Le cas $\Pi\in {\rm Rep}_{\rm tors}(\delta)$ est le contenu du th.~IV.4.7 de \cite{Cbigone}. Le cas 
   $\Pi\in {\rm Rep}_L(\delta)$ se d\'eduit du cas $\Pi\in {\rm Rep}_{\O_L}(\delta)$. Supposons donc que 
   $\Pi\in {\rm Rep}_{\O_L}(\delta)$ et consid\'erons les objets $D_n, \Pi_n$ introduits ci-dessus
de sorte que
   $D\simeq \linv\,  D_n$ et $\Pi^\dual \simeq \linv\,  \Pi_n^\dual $. Puisque chacune des paires $(D_n,\delta^{-1})$ est $G$-compatible, le module 
   $(D^{\natural}\boxtimes_{\delta^{-1}}\p1)_{\rm ns}\simeq \linv (D_n^{\natural}\boxtimes_{\delta^{-1}}\p1)$ est stable par 
   $G$, donc $(D,\delta^{-1})$ est $G$-compatible. Un argument identique d\'emontre le {\rm (ii)} \`a partir du cas de torsion. 
   
    Passons \`a la preuve du {\rm (iii)}. Soit $H={{\rm SL}_2(\qp)}$.
Si $\tilde\Pi=\Pi/\Pi^H$, on a 
     $D(\tilde\Pi)\simeq D(\Pi)$ (puisque $\Pi^H$ est un $\O_L$-module de type fini d'apr\`es le lemme \ref{unitdimfin}, 
donc tu\'e par le foncteur $\Pi\mapsto D(\Pi)$) et $\tilde\Pi^H=0$ (lemme \ref{sl2tors1}). 
De plus, $\tilde\Pi^\dual $ est l'orthogonal de $\Pi^H$, et on est ramn\'e \`a prouver
que $\beta_{\p1}(\tilde\Pi^\dual )\subset D(\tilde\Pi)^{\natural}\boxtimes_{\delta^{-1}}\p1$.
     Autrement dit, on peut supposer que $\Pi^H=0$.

 Notons $Z_n$ l'orthogonal de $\Pi_n^{H}$ dans $\Pi_n^\dual $, de telle sorte que la suite exacte
 $$0\to \Pi_n^H\to \Pi_n\to \Pi_n/\Pi_n^H\to 0$$ nous donne une exacte de $\O_L$-modules profinis $0\to Z_n\to \Pi_n^\dual \to (\Pi_n^H)^\dual \to 0$.
 En passant \`a la limite projective, on obtient
 $$0\to \lim_{\longleftarrow} Z_n\to \Pi^\dual \to (\lim\limits_{\longrightarrow} \Pi_n^H)^\dual \to 0.$$
 La prop.~\ref{isom} montre qu'il existe $N$ tel que $p^N$ tue
 $(\lim\limits_{\longrightarrow} \Pi_n^H)^\dual $; on en d\'eduit que
donc $p^N\Pi^\dual \subset \lim\limits_{\longleftarrow} Z_n$. Comme
   $\beta_{\p1}(Z_n)\subset D_n^{\natural}\boxtimes_{\delta^{-1}}\p1$ d'apr\`es le cas de torsion, 
 on obtient $$\beta_{\p1}(p^N\Pi^\dual )\subset
 \lim\limits_{\longleftarrow} (D_n^{\natural}\boxtimes_{\delta^{-1}}\p1)=
(D^{\natural}\boxtimes_{\delta^{-1}}\p1)_{\rm ns}$$
 et donc $\beta_{\p1}(\Pi^\dual )\subset D^{\natural}\boxtimes_{\delta^{-1}}\p1$, ce qui permet de conclure.

   \begin{theo}\label{recoverPi}
Si $\Pi\in {\rm Rep}_{\rm tors}(\delta)\cup {\rm Rep}_{\O_L}(\delta)$ 
(resp. ${\rm Rep}_L(\delta)$), la transpos\'ee $\beta_{\p1}^\dual$ de $\beta_{\p1}$
induit un morphisme $G$-\'equivariant
   $$\beta_{\p1}^\dual: \Pi_{\delta}(\check{D}(\Pi))\to \Pi/\Pi^{{\rm SL}_2(\qp)},$$
    dont les noyau et conoyau sont de type fini sur $\O_L$ (resp. de dimension finie sur $L$). 
    Plus pr\'ecis\'ement, ${\rm Coker}(\beta_{\p1}^\dual)$
 est un quotient de $((\Pi^\dual )^{{\rm SL}_2(\qp)})^\dual $.
  \end{theo}

  \demo 
   On peut supposer que $\Pi$ est une $\O_L$-repr\'esentation. Soit 
  $\tilde\Pi =\Pi/\Pi^{{\rm SL}_2(\qp)}$. La prop.~\ref{betap1} montre que 
  $\beta_{\p1}(\tilde\Pi ^\dual )\subset D(\Pi)^{\natural}\boxtimes_{\delta^{-1}}\p1=\Pi_{\delta}(\check{D}(\Pi))^\dual $
  (la derni\`ere \'egalit\'e suit du cor.~\ref{checkstab}).  
  Puisque le noyau de $\beta_{\p1}$ est un sous-$\O_L$-module de 
  $(\Pi^\dual )^{{\rm SL}_2(\qp)}$, le conoyau de $\beta_{\p1}^\dual$
 est un quotient de $((\Pi^\dual )^{{\rm SL}_2(\qp)})^\dual $,
 qui est un $\O_L$-module de type fini (rem.~\ref{finitors2}).

  Pour conclure, il nous reste \`a prouver que ${\rm Coker}(\beta_{\p1})$ est un $\O_L$-module de type fini.
    Comme $\tilde\Pi ^\dual $ est compact, $M={\rm Res}_{\qp} (\beta_{\p1}(\tilde\Pi ^\dual ))$ est un sous-$P$-module compact
de $D(\Pi)^{\natural}\boxtimes_{\delta^{-1}}\qp$. De plus, ${\rm Res}_{\zp}(M)$ engendre $D(\Pi)$
(car ${\rm Res}_{\zp}(\beta_{\p1}(\Pi_n^\dual ))=\beta_{\zp}(\Pi_n^\dual )$ engendre $D(\Pi_n)$ par construction
m\^eme) donc $M=D(\Pi)^{\natural}\boxtimes_{\delta^{-1}}\qp$
(cf.~(i) de la prop.~\ref{stabilite1}).
On en d\'eduit (rem.~\ref{resqp}) que $D(\Pi)^{\natural}\boxtimes_{\delta^{-1}}\p1\subset \beta_{\p1}(\tilde\Pi ^\dual )+(0, D(\Pi)^{\rm nr})$,
et donc ${\rm Coker}(\beta_{\p1})$ est un quotient
de $(0, D(\Pi)^{\rm nr})$
ce qui permet de conclure puisque $D(\Pi)^{\rm nr}$ est un $\O_L$-module de type fini.

  \begin{coro}\label{longf}
   Tout objet de ${\rm Rep}_{\rm tors}(G)$ ou
de ${\rm Rep}_L(G)$ est de longueur finie comme $B$-module (topologique).
  \end{coro}
\demo
 C'est une cons\'equence imm\'ediate du th.~\ref{recoverPi}
 et de la prop.~\ref{longueur}.

\begin{coro}\label{absirr}
 Soit $\Pi\in {\rm Rep}_L(G)$ supersinguli\`ere, de caract\`ere central $\delta$. Alors $D(\Pi)$ est absolument irr\'eductible de dimension
    $\geq 2$ et on a des isomorphismes topologiques de $G$-modules
     $$\Pi^\dual \simeq D(\Pi)^{\natural}\boxtimes_{\delta^{-1}}\p1, \quad \Pi\simeq \Pi_{\delta}(\check{D}(\Pi)).$$
\end{coro}

\demo 
Comme $\Pi$ est irr\'eductible de dimension infinie, on a $\Pi^{{\rm SL}_2(\qp)}=0$ et
$(\Pi^\dual )^{{\rm SL}_2(\qp)}=0$. Le th.~\ref{recoverPi} fournit une suite exacte $0\to K\to \Pi_{\delta}(\check{D}(\Pi))\to \Pi\to 0$, avec $\dim_L(K)<\infty$.

L'irr\'eductibilit\'e de $D(\Pi)$ est une cons\'equence du (i)
de la prop.~\ref{irred}.

  Si $\dim_{{\cal E}}(D(\Pi))=1$, il d\'ecoule du
th. \ref{recoverPi} et de la prop.
  \ref{dim1} que $\Pi$ est ordinaire, ce qui est contraire \`a l'hypoth\`ese.
Donc $\dim_{{\cal E}} D(\Pi)\geq 2$ et, puisque $D(\Pi)$ est irr\'eductible, 
on a $\check{D}(\Pi)^{\sharp}=
  \check{D}(\Pi)^{\natural}$. La prop.~\ref{plusgrande} permet de conclure que $K=0$, 
et donc $\Pi\simeq \Pi_{\delta}(\check{D}(\Pi))$.
  On conclut en utilisant le cor. \ref{checkstab}.

\begin{rema}
 On d\'eduit du cor. \ref{absirr} que si $\Pi_1$, $\Pi_2$ sont supersinguli\`eres, 
{\it de m\^eme caract\`ere central} et si $D(\Pi_1)\simeq D(\Pi_2)$, alors
 $\Pi_1\simeq \Pi_2$. Cette propri\'et\'e d'injectivit\'e du foncteur $\Pi\mapsto D(\Pi)$ est d\'emontr\'ee
 par voie tr\`es d\'etourn\'ee dans \cite{Pa} (voir la preuve du th.~10.4 de loc.cit.), mais l'approche de loc.cit.
 fournit plus d'informations. On y prouve que si $p\geq 5$ et si $\Pi\in {\rm Rep}_L(\delta)$ est supersinguli\`ere, alors
 $D(\Pi)$ est de dimension $2$ sur ${\cal E}$ et $\delta=\chi^{-1}\det \check{D}(\Pi)$.  
En particulier, l'image par
  le foncteur $\Pi\mapsto D(\Pi)$ suffit \`a retrouver le caract\`ere central, ce qui est assez surprenant
  car le $(\varphi,\Gamma)$-module attach\'e \`a un $\Pi\in {\rm Rep}_L(G)$ n'utilise que la restriction
  au mirabolique. 
  \end{rema}

\subsection{Reconstruction de $D$}

  Le but de ce paragraphe est de d\'emontrer que l'on peut r\'ecup\'erer $D$ \`a partir de $\Pi_{\delta}(D)$
  quand $(D,\delta)$ est une paire $G$-compatible.  
On note $K={\rm GL}_2(\zp)$ le sous-groupe compact maximal de $G$ et
$Z$ son centre.
  
  \begin{theo}\label{recoverD}
    Soit $(D,\delta)$ une paire $G$-compatible, avec $D$ dans une des cat\'egories $\fget_{\rm tors}$, 
    $\fget(\oe )$ ou $\fget({\cal E})$. Alors on a un isomorphisme canonique de
    $(\varphi,\Gamma)$-modules
    $D(\Pi_{\delta}(D))\simeq \check{D}$.
  \end{theo}
  
  \demo Le cas $D\in \fget({\cal E})$ se d\'eduit 
  du cas $D\in \fget(\oe )$ en tensorisant par $L$.
  Supposons que $D\in \fget(\oe )$ et posons 
  $D_n=D/p^n\in \fget_{\rm tors}$. On d\'eduit de la prop.~\ref{exact}
  que l'application naturelle
  $D\to D_n$ induit un isomorphisme $D(\Pi_{\delta}(D)/p^n \Pi_{\delta}(D))=D(\Pi_{\delta}(D_n))$ pour tout $n$, ce qui montre
  qu'il suffit de traiter le cas $D\in\fget_{\rm tors}$, ce que l'on supposera dans la suite. 

\smallskip
  Soit $W$ l'image de $\tilde W=\sum_{g\in K} g\cdot D^{\sharp}$ dans $\Pi_{\delta}(D)$. 
    
    \begin{lemm}
     $W$ est un sous-$KZ$-module de $\Pi_{\delta}(D)$, de longueur finie comme 
     $\O_L$-module et $W$ engendre $\Pi_{\delta}(D)$ comme $\O_L[G]$-module. 
    \end{lemm}
   
   \demo Il est clair que $W$ est stable par $KZ$.
  Soient
  $z_1,...,z_d\in D^{\sharp}$ tels que $D^{\sharp}=\cup_{i=1}^{d} (D^++z_i)$. 
Puisque $D^+\subset D^{\natural}\boxtimes_{\delta}\p1$ et $D^{\natural}\boxtimes_{\delta}\p1$
est stable par $G$, on a $g D^{+}\subset D^{\natural}\boxtimes_{\delta}\p1$ pour tout 
$g\in G$. Si $v_i$ est l'image de $z_i$ dans $\Pi_{\delta}(D)$, on conclut que 
$W$ est le sous-$\O_L[K]$-module de $\Pi_{\delta}(D)$ engendr\'e par 
$v_1,...,v_d$. Comme $\Pi_{\delta}(D)$ est lisse, le $\O_L$-module $W$ est de longueur finie.

   Montrons enfin que $W$ engendre $\Pi_{\delta}(D)$ comme $\O_L[G]$-module. Il suffit de v\'erifier
   que $D\boxtimes_{\delta}\p1=\sum_{g\in G} g\cdot D^{\sharp}$ et comme $D\boxtimes_{\delta}\p1=D+w\cdot D$, 
   il suffit de prouver l'inclusion $D\subset \sum_{g\in G} g\cdot D^{\sharp}$. 
   Or, si $z\in D$, alors $z_{n,i}:=\psi^{n}((1+T)^{-i} z)\in D^{\sharp}$ pour tout
    $n$ assez grand, uniform\'ement en $i\in \zp$, et 
    $$z=\sum_{i=0}^{p^n-1} (1+T)^i \varphi^{n} (z_{n,i})=\sum_{i=0}^{p^n-1} \matrice{1}{i}{0}{1}\cdot \matrice{p^n}{0}{0}{1}z_{n,i},$$
    ce qui permet de conclure.

\smallskip
     Rappelons que $D_W^+$ est l'orthogonal de $\sum_{g\in P-P^+} g\cdot W$ dans 
     $\Pi_{\delta}(D)^\dual =\check{D}^{\natural}\boxtimes_{\delta^{-1}} \p1$. 
L'isomorphisme $\oe \otimes_{\O_L[[T]]} D_W^+\cong \check D$
que l'on cherche \`a \'etablir est une cons\'equence de la platitude de $\oe $ sur $\O_L[[T]]$
et des deux lemmes suivants.
        
      \begin{lemm}
       On a $D_{W}^+\subset \check{D}^{++}$. 
      \end{lemm}
   \demo Si $\check{z}\in D_W^+$, alors $\{\check{z}, gk D^{\sharp}\}_{\p1}=0$ pour tout
   $g\in P-P^+$ et $k\in K$, donc  
    ${\rm Res}_{\zp}(k^{-1}g^{-1}\check{z})$ est orthogonal \`a 
   $D^{\sharp}$. On en d\'eduit que ${\rm Res}_{\zp}(gz)\in \check{D}^{++}$
   pour tout $g\in M:= \{zx h^{-1},\ z\in Z,\  k\in K,\  h\in P-P^+\}$. Si $n\geq 1$ et 
   $0\leq i<p^n$ est un multiple de $p$, alors\footnote{ Pour $i=0$ cela d\'ecoule de l'identit\'e $  \matrice{p^{-n}}{0}{0}{1}\cdot w=\matrice{p^{-n}}{0}{0}{p^{-n}} w \matrice{p^{-n}}{0}{0}{1}^{-1}$, et
   si $i\ne 0$, de l'identit\'e 
   (avec $N=n-2v_p(i)$)
   $$ \matrice{p^{-n}}{0}{0}{1}\cdot \matrice{1}{-i}{0}{1}\cdot w= \matrice{1/i}{0}{0}{1/i}\cdot \matrice{-i^2 p^{N-n}}{0}{p^Ni}{1} \matrice{ p^N}{1/i}{0}{1}^{-1}.$$}
    $ \matrice{p^{-n}}{0}{0}{1}\cdot \matrice{1}{-i}{0}{1}\cdot w\in M$.
     On en d\'eduit que, pour tout $0\leq i<p^n$ multiple de $p$, l'on a $$\psi^{n}((1+T)^{-i} {\rm Res}_{\zp}(w\check{z}))\in \check{D}^{++}.$$
   En faisant $n\to\infty$ dans l'\'egalit\'e 
   $$ {\rm Res}_{p\zp}(w\check{z})=
\sum_{i<p^n, p|i} (1+T)^i \varphi^n(\psi^n((1+T)^{-i} {\rm Res}_{\zp}(w\check{z}) ),$$
   on obtient ${\rm Res}_{p\zp}(w\check{z})=0$, i.e. $\check{z}\in \check{D}$. De plus, 
   $\varphi(\check{z})={\rm Res}_{\zp} ( \matrice{p^{-1}}{0}{0}{1}^{-1}\check{z})\in \check{D}^{++}$, donc 
   $\check{z}\in \check{D}^{++}$, ce qui permet de conclure.

\begin{lemm}
Si $n$ est assez grand, $\varphi^n(T)\check{D}^{++}\subset D_W^+$.
\end{lemm}
\demo
Il s'agit de prouver que $\varphi^n(T)\check{D}^{++}$ est orthogonal \`a
$g\cdot \tilde W$, si $g=\matrice{a}{b}{0}{1}$ avec $v_p(a)<0$ ou $v_p(b)<0$.
De mani\`ere \'equivalente, il s'agit de v\'erifier
que $\check{D}^{++}$ est orthogonal \`a $\big(\matrice{1}{-p^n}{0}{1}-1\big)\cdot g\cdot \tilde W$.
L'argument est diff\'erent suivant que $v_p(a)\geq 0$ ou $v_p(a)<0$.

$\bullet$ Si $v_p(a)\geq 0$, on a $v_p(b)<0$.
Or il existe un treillis\footnote{Un treillis est un sous-$\O_L[[T]]$-module
compact de $D$ dont l'image modulo $p^n$ est ouverte dans $D/p^nD$ pour tout $n\in\N$
(il suffit qu'elle le soit dans $k_L\otimes D$).}
  $M$ de $D$ tel que $\matrice{a}{0}{0}{1}\cdot \tilde W$
soit inclus dans $D+wM+(D^\natural\boxtimes_\delta\p1)$ pour tout $a\in\Z_p^\dual$.
Choisissons $n$ assez grand pour que $\big(\matrice{1}{-p^n}{0}{1}-1\big)\cdot wM\subset
D^\natural\boxtimes_\delta\p1$.
En \'ecrivant $\big(\matrice{1}{-p^n}{0}{1}-1\big)\cdot g$ sous la forme
$\matrice{1}{b}{0}{1}\cdot \big(\matrice{1}{-p^n}{0}{1}-1\big)\cdot \matrice{a}{0}{0}{1}$,
on voit que 
$$\big(\matrice{1}{-p^n}{0}{1}-1\big)\cdot g\cdot \tilde W\subset (D\boxtimes_\delta (b+\Z_p))+
(D^\natural\boxtimes_\delta\p1),$$ qui est orthogonal \`a
$\check{D}^{++}$ car $D\boxtimes_\delta (b+\Z_p)$ l'est puisque $\Z_p\cap(b+\Z_p)=\emptyset$
et $D^\natural\boxtimes_\delta\p1$ l'est puisque $\check{D}^{++}$ est inclus
dans $D^\natural\boxtimes_{\delta^{-1}}\p1$.

$\bullet$ Si $v_p(a)<0$, on \'ecrit $\big(\matrice{1}{-p^n}{0}{1}-1\big)\cdot g$
sous la forme $g\cdot \big(\matrice{1}{-a^{-1}p^n}{0}{1}-1\big)$, et on choisit
$n$ assez grand pour que $\big(\matrice{1}{c}{0}{1}-1\big)\cdot \tilde W\subset
(D^\natural\boxtimes_\delta\p1)$ pour tout $c\in p^n\Z_p$.
Les arguments ci-dessus montrent qu'alors
$\check{D}^{++}$ est orthogonal \`a $\big(\matrice{1}{-p^n}{0}{1}-1\big)\cdot g\cdot \tilde W$.

Ceci permet de conclure.

\begin{rema}
Le module $W$ utilis\'e dans la preuve du th.~\ref{recoverD}
est raisonnablement naturel, et garde un sens si $D$ n'est pas de torsion.
Cela semble un bon candidat si on veut \'ecrire $\Pi_\delta(D)$
comme quotient d'une induite de $KZ$ \`a $G$.
\end{rema}

\section{Repr\'esentations localement analytiques}
Dans ce chapitre, on revisite les travaux de Schneider et Teitelbaum~\cite{STInv}
sur les repr\'esentations localement analytiques en \'etudiant de plus pr\`es
la filtration naturelle par rayon d'analyticit\'e (on fera attention
que cette notion est l\'eg\`erement diff\'erente de celle \`a
laquelle on penserait naturellement (cf.~rem.~\ref{exac2})).

\subsection{Groupes uniformes}\label{uniformes}

   \quad On renvoie \`a \cite{Sautoy} pour les preuves des r\'esultats \'enonc\'es ci-dessous.
   Posons $\kappa=1$ si $p>2$ et $\kappa=2$ sinon. Dans ce chapitre, $H$ est un pro-$p$-groupe uniforme, i.e.
un pro-$p$-groupe topologiquement de type fini, sans $p$-torsion et tel que $[H,H]\subset H^{p\kappa}$.

Si $i\geq 0$, soit $H_i=\{g^{p^i}| g\in H\}=H^{p^i}$. Alors $H_i$ est un sous-groupe
ouvert distingu\'e de $H$, et  $(H_i)_{i\geq 0}$ est
 un syst\`eme fondamental de voisinages ouverts de $1$. 
  En posant $\omega(1)=\infty$ et $\omega(g)=i$ si $g\in H_{i-\kappa} \backslash H_{i-\kappa+1}$, on obtient une $p$-valuation (au sens de Lazard) 
 satisfaisant l'hypoth\`ese HYP de \cite{STInv}, ce qui nous permet d'utiliser directement les r\'esultats de loc.cit. Si $h_1,h_2,...,h_d$ est un syst\`eme minimal de g\'en\'erateurs topologiques
 de $H$, alors $\omega(h_i)=\kappa$ pour tout $i$, et l'application $\zp^d\to H$ d\'efinie par $(x_1,x_2,...,x_d)\mapsto h_1^{x_1}h_2^{x_2}\cdots h_d^{x_d}$
 est un hom\'eomorphisme. De plus, on a
    $$\omega(h_1^{x_1}h_2^{x_2}\cdots h_d^{x_d})=\kappa+\min_{1\leq i\leq d} v_p(x_i)$$
    pour tous
   $x_1,...,x_d\in \zp$.
   
   On utilise les notations standard pour les $d$-uplets: $|\alpha|=\sum_{i=1}^{d}\alpha_i$, $\binom{\alpha}{\beta}=\prod_{i} \binom{\alpha_i}{\beta_i}$, $h^{\alpha}=\prod_{i} h_i^{\alpha_i}$, etc. On \'ecrit $\alpha\leq\beta$ si $\alpha_i\leq \beta_i$ pour tout $i$. Si $\alpha=(\alpha_1,...,\alpha_d)\in\mathbf{N}^d$, on note
   $$b^{\alpha}=(h_1-1)^{\alpha_1} (h_2-1)^{\alpha_2}\cdots (h_d-1)^{\alpha_d}\in \zp[H].$$
   
   \subsection{Coefficients de Mahler}\label{coefMahler}

         L'espace ${\cal C}(H)$ des fonctions continues sur $H$, \`a valeurs dans $L$
     est une $L$-repr\'esentation de Banach de $H$, si on le munit de la norme sup
     et de l'action de $H$ d\'efinie par $(g\cdot \phi)(x)=\phi(xg)$. 
   Si $\phi\in {\cal C}(H)$, on note 
    $$a_{\alpha}(\phi)=(b^{\alpha} \phi)(1)=\sum_{\beta\leq \alpha} (-1)^{\alpha-\beta}\binom{\alpha}{\beta}\phi(h^{\beta})$$
    ses coefficients de Mahler, relativement au choix des coordonn\'ees $h_1,...,h_d$ sur $H$. Soit $\phi_{\alpha}\in{\cal C}(H)$ l'application d\'efinie par 
    $\phi_{\alpha}(h^x)=\binom{x}{\alpha}$ pour $x\in \zp^d$. Un th\'eor\`eme classique de Mahler 
    montre que $(\phi_{\alpha})_{\alpha\in\mathbf{N}^d}$ est une base orthonormale de 
    ${\cal C}(H)$, et pour tout $\phi\in {\cal C}(H)$ 
 $$\phi=\sum_{\alpha\in\mathbf{N}^d} a_{\alpha}(\phi)\cdot \phi_{\alpha}.$$

 \begin{defi}
  Pour tout $h\in \mathbf{N}^\dual $ on note 
    $${\rm LA}^{(h)}(H)=\{\phi\in {\cal C}(H)| \lim_{|\alpha|\to\infty} (v_p(a_{\alpha}(\phi))-r_h|\alpha|)=\infty\}.$$
   C'est un espace de Banach pour la valuation $v^{(h)}$ d\'efinie par 
   $$v^{(h)}(\phi)=\inf_{\alpha} (v_p(a_{\alpha}(\phi))-r_h|\alpha|).$$

 \end{defi}
 
   D'apr\`es le th\'eor\`eme d'Amice~\cite{Am64}, l'espace ${\cal C}^{\rm an}(H)$ des fonctions localement analytiques sur $H$ est la limite inductive 
   des espaces ${\rm LA}^{(h)}(H)$.

  \subsection{Compl\'etions de l'alg\`ebre des mesures de $H$}\label{dh}
 
 Le but de ce paragraphe est de rappeler un certain nombre de constructions et r\'esultats de \cite{STInv}, et d'\'etablir
    quelques estim\'ees techniques dont on aura besoin plus loin.

    Soit $\Lambda(H)$ le dual faible de ${\cal C}(H)$; on a aussi $\Lambda(H)=L\otimes_{\O_L}\big(\linv\,\O_L[H/H^{p^n}]\big)$. 
C'est une alg\`ebre (pour le produit de convolution) topologique localement compacte, et le \S~\ref{coefMahler} montre que 
    tout \'el\'ement de $\Lambda(H)$ s'\'ecrit de mani\`ere unique $\lambda=\sum_{\alpha\in\mathbf{N}^d}
   c_{\alpha} b^{\alpha}$, avec $(c_{\alpha})_{\alpha}$ une suite born\'ee dans $L$. La valeur de 
   $\lambda$ en $\phi\in{\cal C}(H)$ est donn\'ee par 
   $$\langle \lambda, \phi\rangle=\sum_{\alpha\in\mathbf{N}^d} c_{\alpha} a_{\alpha}(\phi).$$
   
      On note ${\cal D}_{1/p}(H)$ le compl\'et\'e de $\Lambda(H)$ pour la valuation d'alg\`ebre
(i.e.~$v_{1/p}(\lambda\mu)\geq v_{1/p}(\lambda)+v_{1/p}(\mu)$):
      $$v_{1/p} (\sum_{\alpha} c_{\alpha} b^{\alpha})=\inf_{\alpha} (v_p(c_{\alpha})+\kappa |\alpha|)$$
      et, si $h\in \mathbf{N}^\dual $, on d\'efinit\footnote{Rappelons que 
   $r_h=\frac{1}{p^{h-1}(p-1)}$.}
    $${\cal D}_{h}(H)=\{\sum_{\alpha} c_{\alpha} b^{\alpha}\in{\cal D}_{1/p}(H),\  \inf_{\alpha} (v_p(c_{\alpha})+r_h|\alpha|)>-\infty\},$$
    que l'on munit d'une valuation d'alg\`ebre $v^{(h)}$ en posant 
   $$v^{(h)}(\sum_{\alpha} c_{\alpha} b^{\alpha})=\inf_{\alpha} (v_p(c_{\alpha})+r_h|\alpha|).$$
   Enfin, on note ${\cal D}(H)$ le dual topologique (fort) de ${\cal C}^{\rm an}(H)$. C'est aussi la limite projective des
   ${\cal D}_{h}(H)$. L'accouplement $$\langle \lambda, \phi\rangle=\sum_{\alpha\in\mathbf{N}^d} c_{\alpha} a_{\alpha}(\phi)$$
   identifie ${\cal D}_{h}(H)$ au dual topologique de ${\rm LA}^{(h)}(H)$, ce qui permet de d\'efinir une topologie faible sur 
   ${\cal D}_{h}(H)$. Dans la suite on munit\footnote{L'alg\`ebre ${\cal D}_{h}(H)$ peut aussi \^etre munie d'une topologie d'alg\`ebre de Banach
  en utilisant la valuation $v^{(h)}$. C'est d'ailleurs ce qui est fait dans \cite{STInv}. Cette topologie 
est  trop forte pour les applications que nous avons en vue.}
 ${\cal D}_{h}(H)$ de cette topologie faible. La diff\'erence avec la topologie forte
   (qui, elle, est induite par la valuation $v^{(h)}$) est que $p^{-r_h|\alpha|}b^{\alpha}$ tend vers $0$ pour la topologie faible dans 
     ${\cal D}_{h}(H)$, mais pas pour la topologie forte. 
     
    \begin{prop} \label{Schneider} La multiplication dans $\Lambda(H)$ s'\'etend par continuit\'e \`a ${\cal D}_h(H)$ 
   et l'inclusion naturelle $\Lambda(H)\to {\cal D}_{h}(H)$
   est plate. De plus, l'inclusion $\Lambda(H)\to {\cal D}(H)$ est fid\`element plate.
    
    \end{prop}
    
    \demo
    Tout ceci est d\'emontr\'e dans le chap.~$4$ de 
    \cite{STInv}. Pour faciliter la comparaison, notons que l'on utilise comme $p$-valuation 
    celle introduite dans le $\S$ \ref{uniformes}, et que si l'on pose 
    $s_h=p^{-r_h \kappa^{-1}}\in ]\frac{1}{p},1[\cap p^{\mathbf{Q}}$, alors $p^{-v^{(h)}}$ correspond \`a $||\cdot ||_{s_h}$ 
    de loc.cit, donc ${\cal D}_{h}(H)$ correspond \`a $D_{<s_h}(H,L)$ de loc.cit.

 \begin{lemm}\label{ST}
 
 {\rm (i)} Soit $h_1',...,h_d'$ un autre syst\`eme minimal de g\'en\'erateurs topologiques de $H$ et soit
 $b'^{\alpha}=(h_1'-1)^{\alpha_1}\cdots (h_d'-1)^{\alpha_d}$. Il existe 
des
 $c_{\alpha,\beta}\in \O_L$ tels que $b'^{\alpha}=\sum_{\beta} c_{\alpha,\beta} b^{\beta}$ et
 $v_p(c_{\alpha,\beta})\geq  \max(0,\kappa (|\alpha|-|\beta|))$.
 
 {\rm (ii)} Pour tout $x\in H$ il existe des $c_{\alpha,\beta,x}\in \O_L$ tels que 
   $xb^{\alpha}=\sum_{\beta} c_{\alpha,\beta,x} b^{\beta}$ et
 $v_p(c_{\alpha,\beta,x})\geq \max(0, \kappa(|\alpha|-|\beta|))$.
 
 {\rm (iii)} Il existe des $c_{\alpha,\beta}\in \O_L$ tels que 
 $$(h_1^p-1)^{\alpha_1}\cdots (h_d^p-1)^{\alpha_d}=\sum_{\beta} c_{\alpha,\beta} b^{\beta}$$
 et $v_p(c_{\alpha,\beta})\geq \max(0, r_1(p|\alpha|-|\beta|))$.
 
 \end{lemm}
 
 \demo {\rm (i)} L'existence des $c_{\alpha,\beta}$ et le fait qu'ils appartiennent \`a $\O_L$ suivent du fait que 
 $b'^{\alpha}\in \O_L[H]\subset \O_L[[H]]$. L'in\'egalit\'e $v_p(c_{\alpha,\beta})\geq  \max(0,\kappa (|\alpha|-|\beta|))$
 d\'ecoule du fait que $v_{1/p}(g-1)\geq \kappa$ pour tout
 $g\in H$. La preuve du {\rm (ii)} est identique et laiss\'ee au lecteur.
 
 {\rm (iii)} Puisque $v^{(1)}$ est une valuation d'alg\`ebre sur $\Lambda(H)$, on a  
 $$\inf_{\beta} (v_p(c_{\alpha,\beta})+r_1|\beta|)\geq \sum_{i=1}^{d} \alpha_i v^{(1)}(h_i^p-1).$$
  Or
 $$v^{(1)}(h_i^p-1)=v^{(1)}\big(\sum_{k=1}^{p} \tbinom{p}{k} (h_i-1)^k\big)=\inf_{k\leq p} v_p\big(\tbinom{p}{k}\big)+kr_1=pr_{1},$$
ce qui permet de conclure.

\subsection{Le foncteur $\Pi\mapsto \Pi^{(h)}$}

 Soit $\Pi$ une 
   $L$-repr\'esentation de Banach de $H$ et soit $v_{\Pi}$ une valuation sur 
   $\Pi$ qui d\'efinit sa topologie. Si $h\geq 1$ on d\'efinit
     $$\Pi^{(h)}=\{v\in \Pi| \lim_{|\alpha|\to\infty} v_{\Pi}(b^{\alpha}v)-r_h|\alpha|=\infty\},$$
     et l'on munit de la valuation $v^{(h)}$ d\'efinie par
     $$v^{(h)}(v)= \inf_{\alpha\in\mathbf{N}^d} (v_{\Pi}(b^{\alpha}v)-r_h|\alpha|).$$
   Alors $\Pi^{(h)}$ est un $L$-espace de Banach, qui ne d\'epend pas du choix de 
   $v_{\Pi}$ (la valuation $v^{(h)}$ en d\'epend de mani\`ere \'evidente, mais 
   changer $v_{\Pi}$ remplace $v^{(h)}$ par une valuation \'equivalente). 

Les $b^\alpha v$ sont les coefficients de Mahler de la fonction $o_v:H\to \Pi$
d\'efinie par $o_v(g)=g\cdot v$, et le th\'eor\`eme d'Amice montre que le sous-espace
$\Pi^{\rm an}$ des vecteurs localement analytiques est la limite inductive
des $\Pi^{(h)}$.
  
       \begin{prop}\label{indepcoord}
   L'espace $\Pi^{(h)}$ et la valuation $v^{(h)}$ ne d\'ependent pas du choix du syst\`eme minimal
   de g\'en\'erateurs topologiques $h_1,...,h_d$ de $H$, que l'on utilise pour d\'efinir $b^{\alpha}$.
  
  \end{prop}
  
  \demo Soit $h_1',...,h_d'$ un autre syst\`eme minimal de g\'en\'erateurs topologiques
  et soit $c_{\alpha,\beta}$ comme dans le lemme \ref{ST}. Soient 
  $v\in\Pi^{(h)}$, $M\in\mathbf{R}$ et $N$ tels que $v_{\Pi}(b^{\alpha}v)-r_h|\alpha|\geq M$ pour tout
  $|\beta|\geq N$. Si $|\alpha|\geq N$, alors (en utilisant le lemme~\ref{ST})
  $$\inf_{|\beta|<N} (v_{\Pi}(c_{\alpha,\beta} b^{\beta}v)-r_h|\alpha|)\geq 
  \inf_{|\beta|<N} v_{\Pi}(b^{\beta}v)+\kappa(|\alpha|-N)-r_h|\alpha|,$$
  quantit\'e qui d\'epasse $M$ si $|\alpha|$ est assez grand, 
  et $$\inf_{|\beta|\geq N} (v_{\Pi}(c_{\alpha,\beta} b^{\beta}v)-r_h|\alpha|)\geq 
  M+r_h(|\beta|-|\alpha|)+\max(0, \kappa(|\alpha|-|\beta|))\geq M.$$
  Comme $b'^{\alpha}v=\sum_{\beta} c_{\alpha,\beta} b^{\beta} v$, 
on d\'eduit des in\'egalit\'es pr\'ec\'edentes 
que l'on a $\lim_{|\alpha|\to\infty} v_{\Pi}(b'^{\alpha}v)-r_h|\alpha|=\infty$ et 
 (en prenant $N=0$ et $M=v^{(h)}(v)$) 
$$\inf_{\alpha} (v_{\Pi}(b^{\alpha}v)-r_h|\alpha|)\leq \inf_{\alpha} (v_{\Pi}(b'^{\alpha}v)-r_h|\alpha|).$$
  Le r\'esultat s'en d\'eduit par sym\'etrie.

     \begin{prop} \label{casfacile}
   On a ${\cal C}(H)^{(h)}={\rm LA}^{(h)}(H)$.
       
  \end{prop}
  
  \demo
  Il est imm\'ediat de v\'erifier que ${\cal C}(H)^{(h)}\subset {\rm LA}^{(h)}(H)$. R\'eciproquement, supposons que 
  $\phi\in {\rm LA}^{(h)}(H)$. On veut montrer que 
  $$\lim_{|\alpha|\to\infty} \inf_{x\in H} v_p((b^{\alpha}\phi)(x))-r_h|\alpha|=\infty.$$
  Si $c_{\alpha,\beta,x}$ est comme dans le lemme \ref{ST}, alors $$(b^{\alpha}\phi)(x)=(xb^{\alpha}\phi)(1)=\sum_{\beta} c_{\alpha,\beta,x} a_{\beta}(\phi).$$
   On conclut en utilisant le lemme 
  \ref{ST}, comme dans la preuve de la 
  prop.~\ref{indepcoord}.

     \begin{theo}\label{carduale}
 {\rm (i)}     Si $v\in\Pi$, alors $v\in \Pi^{(h)}$ si et seulement si 
      la fonction $g\mapsto l(g\cdot v)$ appartient \`a ${\rm LA}^{(h)}(H)$ pour tout 
      $l\in \Pi^\dual $.

{\rm (ii)}  Si $l_1,\dots,l_r$ engendrent $\Pi^\dual$ comme $\Lambda(H)$-module,
alors $v\mapsto \iota(v)$, o\`u 
\linebreak
$\iota(v)\in {\rm LA}^{(h)}(H)^r$ est la fonction
$g\mapsto(l_1(g\cdot v),\dots,l_r(g\cdot v))$, est un plongement ferm\'e
de $\Pi^{(h)}$ dans ${\rm LA}^{(h)}(H)^r$.
     \end{theo}
     
     \demo
      Par d\'efinition, $v\in \Pi^{(h)}$ si et seulement si la suite 
     $x_{\alpha}=p^{-r_h|\alpha|}  b^{\alpha}v$ tend vers $0$ dans 
      $\Pi\otimes_{L} L(p^{r_h})$. Le lemme (surprenant) \ref{weak} ci-dessous montre que cela arrive si et seulement si 
   $p^{-r_h|\alpha|} l(b^{\alpha}v)$ tend vers $0$
      pour tout $l\in \Pi^\dual $. On conclut la preuve du (i)
en remarquant que $l(b^{\alpha}v)=a_{\alpha}(o_{l,v})$, o\`u
      $o_{l,v}: H\to L$ est la fonction $g\mapsto l(gv)$. 

Pour prouver le (ii), notons que $\Pi$ est admissible (car $\Pi^\dual $ est de type fini comme $\Lambda(H)$-module). Puisque 
$l_1,...,l_r$ engendrent $\Pi^\dual $ comme 
$\Lambda(H)$-module, l'application transpos\'ee 
$$\iota: \Pi\to {\cal C}(H)^r, \quad v\mapsto (g\mapsto (l_1(g\cdot v),...,l_r(g\cdot v)))$$
est un plongement ferm\'e \cite{ST1}. Il d\'ecoule de la prop.~\ref{casfacile} que 
$\iota$ envoie $\Pi^{(h)}$ dans ${\rm LA}^{(h)}(H)^r$ et il nous reste \`a montrer que 
l'application $\iota$ ainsi obtenue est un plongement ferm\'e. Supposons que les
$v_n\in \Pi^{(h)}$ sont tels que $\iota(v_n)\to f$ dans ${\rm LA}^{(h)}(H)^r$.
Alors $\iota(v_n)\to f$ dans $ {\cal C}(H)^r$, et donc 
il existe $v\in \Pi$ tel que $f=\iota(v)$. Ainsi $\iota(v)=f\in {\rm LA}^{(h)}(H)^r$. Autrement dit, 
les $g\mapsto l_i(gv)$ appartiennent \`a ${\rm LA}^{(h)}(H)$ pour tout $i$. Comme 
$l_1,...,l_r$ engendrent $\Pi^\dual $ en tant que $\Lambda(H)$-module, et comme 
${\rm LA}^{(h)}(H)$ est un
      $\Lambda(H)$-module, il s'ensuit que $g\mapsto l(gv)$ appartient \`a ${\rm LA}^{(h)}(H)$
 pour tout $l\in \Pi^\dual $. Le (i) permet  de conclure
que $v\in \Pi^{(h)}$.

\begin{lemm} \label{weak}
Dans un espace localement convexe sur un corps sph\'eriquement complet
      une suite converge vers $0$ si et seulement si elle converge faiblement vers $0$.
\end{lemm}
     \demo
 Voir par exemple \cite[th. 5.5.2]{Garcia}.

  \begin{coro}\label{Pihstable}
   $\Pi^{(h)}$ est stable sous l'action de $H$.  \end{coro}
  
\subsection{De $H$ \`a $H^p$}
Le but de ce \S~est d'\'etudier la variation de la filtration
par rayon d'analyticit\'e quand on remplace $H$ par un sous-groupe
(prop.~\ref{dcalage}).
\begin{prop}\label{extzero}
  Soit $\phi\in {\rm LA}^{(h)}(H^p)$ et soit $\tilde{\phi}\in {\cal C}(H)$ l'extension par $0$ de $\phi$.
  Alors $\tilde{\phi}\in {\rm LA}^{(h+1)}(H)$. 
\end{prop}

\demo
 Il s'agit de montrer que $\lim_{|\alpha|\to\infty} v_p(a_{\alpha}(\tilde{\phi}))-r_{h+1}|\alpha|=\infty$. 
  Si $\beta\in\mathbf{N}^d$, on a $h^{\beta}\in H^p$ si et seulement si 
    $p$ divise $\beta$ (i.e. $p$ divise chaque $\beta_i$). On en d\'eduit que\footnote{Attention au fait que 
    les coefficients de Mahler de $\phi$ sont calcul\'es par rapport \`a $(h_1^p,...,h_d^p)$. On a donc 
    $\phi(h^{p\beta})=\sum_{\gamma} a_{\gamma}(\phi)\binom{\beta}{\gamma}$.}
\begin{align*}
    a_{\alpha}(\tilde{\phi})& =\sum_{p\beta\leq \alpha} (-1)^{\alpha-p\beta}\tbinom{\alpha}{p\beta} \phi(h^{p\beta})=
    \sum_{p\beta\leq \alpha} (-1)^{\alpha-p\beta}\tbinom{\alpha}{p\beta} \cdot \sum_{\gamma} a_{\gamma}(\phi)\tbinom{\beta}{\gamma}\\
    &=\sum_{\gamma} a_{\gamma}(\phi)\cdot c_{\alpha,\gamma},\quad \text{avec}\ c_{\alpha,\gamma}=\sum_{p\beta\leq \alpha} (-1)^{\alpha-p\beta} \tbinom{\alpha}{p\beta}\tbinom{\beta}{\gamma}.
\end{align*}
    On conclut comme dans la preuve de la prop.~\ref{indepcoord}, en utilisant le lemme ci-dessous:
    
    \begin{lemm}
     Si $c_{\alpha,\gamma}=\sum_{p\beta\leq \alpha} (-1)^{\alpha-p\beta} \binom{\alpha}{p\beta}\binom{\beta}{\gamma}$, alors 
     $$v_p(c_{\alpha,\gamma})>\frac{|\alpha|}{p}-d-|\gamma|.$$
    \end{lemm}
    
\demo
Soit 
    $f_k=\psi(T^k)$ pour $k\in\mathbf{N}$. Alors $f_k\in \zp[T]$ et 
    $$ f_k(T^p-1)=\frac{1}{p}\sum_{\zeta^p=1} (\zeta T-1)^k=\sum_{pj\leq k} (-1)^{k-pj}\tbinom{k}{pj}T^{pj}.$$
    On en d\'eduit que pour tout $\alpha\in\mathbf{N}^d$ on a 
    $$\sum_{p\beta\leq \alpha} (-1)^{\alpha-p\beta} \tbinom{\alpha}{p\beta} T^{\beta}=\prod_{i=1}^{d} f_{\alpha_i}(T_i-1).$$
En d\'erivant cette relation $\gamma$-fois et en \'evaluant en $1$, on obtient: 
    $$c_{\alpha,\gamma}=\frac{1}{\gamma!} \frac{d^{\gamma}}{dT^{\gamma}}(\psi(T_1^{\alpha_1})\cdots\psi(T_d^{\alpha_d}))(0)=b_{\alpha_1, \gamma_1}\cdots 
    b_{\alpha_d,\gamma_d},$$
    o\`u $\psi(T^k)=\sum_{i\leq k/p} b_{k,i} T^i$. Le r\'esultat d\'ecoule alors du lemme I.8 de \cite{Cserieunit}, qui fournit l'in\'egalit\'e
    $v_p(b_{k,i})>\frac{k}{p}-1-i$.

  \begin{prop} \label{dcalage}
   Si $\Pi_1$ est la restriction de $\Pi$ \`a $H^p$, alors $\Pi^{(h+1)}= \Pi_1^{(h)}$.
     \end{prop}
  
  \demo
   On note $h'_i=h_i^p$ et $b'^{\alpha}=(h'_1-1)^{\alpha_1}\cdots (h_d'-1)^{\alpha_d}$.
  Alors $h_1',...,h_d'$ forment un syst\`eme minimal de g\'en\'erateurs topologiques
   de $H^p$. Commençons par montrer l'inclusion  $\Pi^{(h+1)}\subset \Pi_1^{(h)}$. 
   La prop.~\ref{indepcoord} montre qu'il suffit de prouver que 
   $\lim_{|\alpha|\to\infty} v_{\Pi}(b'^{\alpha}v)-r_h|\alpha|=\infty$ pour tout $v\in \Pi^{(h+1)}$. 
   La preuve est identique \`a celle de la prop.~\ref{indepcoord}, en utilisant le lemme \ref{ST}.

      Soit maintenant $v\in \Pi_1^{(h)}$ et soient $(g_i)_{i\in I}$ tels que $H=\coprod_{i\in I} H^p g_i$. Soit enfin 
      $l\in \Pi^\dual $ et notons $\phi: H\to L$  l'application $g\to l(gv)$. 
      Puisque 
      $H^p$ est distingu\'e dans $H$ et $v\in\Pi_1^{(h)}$, on montre comme dans la preuve du cor.~\ref{Pihstable} 
      que $(g_i\phi)|_{H^p}\in {\rm LA}^{(h)}(H^p)$. Si $\phi_i$ est l'extension par z\'ero de 
      $(g_i\phi)|_{H^p}$, la prop.~\ref{extzero} montre que $\phi_i\in {\rm LA}^{(h+1)}(H)$. On en d\'eduit que 
      $\phi=\sum_{i\in I} g_i^{-1}\cdot \phi_i\in {\rm LA}^{(h+1)}(H)$, puisque ${\rm LA}^{(h+1)}(H)$ est stable par 
      $H$. Le th.~\ref{carduale} permet de conclure que $v\in\Pi^{(h+1)}$.

  \begin{rema} \label{decompLA}
    On d\'eduit de la preuve de la prop.~\ref{dcalage} que, si 
on d\'ecompose $H$ sous la forme
    $H=\coprod_{i\in I} H^p g_i$, 
l'application $\phi\mapsto ((g_i\cdot \phi)|_{H^p})_{i\in I}$ induit un isomorphisme d'espaces de Banach
   $${\rm LA}^{(h+1)}(H)\simeq \bigoplus_{i\in I} {\rm LA}^{(h)}(H^p, L).$$
  \end{rema}

    \subsection{Exactitude du foncteur $\Pi\mapsto \Pi^{(h)}$}

  A partir de maintenant, on suppose que $\Pi$ une repr\'esentation de Banach {\it admissible} de
  $H$. Cela signifie que $\Pi^\dual $ est un $\Lambda(H)$-module de type fini (et donc de pr\'esentation finie, puisque
  $\Lambda(H)$ est noeth\'erien). Le choix d'une surjection $\Lambda(H)^r\to \Pi^\dual $ induit un plongement ferm\'e
  $\iota: \Pi\to{\cal C}(H, L)^r$, ainsi qu'une surjection $\sigma:{\cal D}_{h}(H)^r\to {\cal D}_{h}(H)\otimes_{\Lambda(H)} \Pi^\dual $.
On munit ${\cal D}_{h}(H)\otimes_{\Lambda(H)}\Pi^\dual $ de la topologie quotient, induite par $\sigma$ (l'espace ${\cal D}_{h}(H)$ ayant la topologie faible d\'efinie dans le \S~\ref{dh}).
  La topologie faible de $\Pi^\dual $ est alors la topologie quotient induite par la surjection $\Lambda(H)^r\to \Pi^\dual $.  
  
    Si $v\in\Pi^{(h)}$ et si $\lambda=\sum_{\alpha} c_{\alpha} b^{\alpha}\in{\cal D}_{h}(H)$, la s\'erie
    $\sum_{\alpha} c_{\alpha} b^{\alpha}$ converge dans $\Pi^{(h)}$. Cela permet de munir $\Pi^{(h)}$ (et donc $(\Pi^{(h)})^\dual $) d'une structure de 
    ${\cal D}_{h}(H)$-module. L'application naturelle
   $\Pi^\dual \to (\Pi^{(h)})^\dual $ induit une application continue \footnote{On munit dans la suite
  $\Pi^\dual $ et $(\Pi^{(h)})^\dual $ de la topologie faible de dual de Banach.}  
   ${\cal D}_{h}(H)$-lin\'eaire ${\cal D}_{h}(H)\otimes_{\Lambda(H)}\Pi^\dual \to
    (\Pi^{(h)})^\dual $. Par passage \`a la limite on obtient aussi une application 
    ${\cal D}(H)$-lin\'eaire ${\cal D}(H)\otimes_{\Lambda(H)} \Pi^\dual \to (\Pi^{\rm an})^\dual $.

    \begin{prop}\label{dualdepih}
     {\rm (i)} L'application ${\cal D}_{h}(H)\otimes_{\Lambda(H)}\Pi^\dual \to
    (\Pi^{(h)})^\dual $ est un isomorphisme de $L$-espaces vectoriels topologiques.
    
    {\rm (ii)} L'application ${\cal D}(H)\otimes_{\Lambda(H)} \Pi^\dual \to (\Pi^{\rm an})^\dual $ est un isomorphisme
    de $L$-espaces vectoriels topologiques.
    
    \end{prop}
    
    \demo
     La preuve est fortement inspir\'ee de la preuve du th.~7.1 de \cite{STInv} (qui est pr\'ecis\'ement la partie {\rm (ii)} de la proposition). 
    Le {\rm (ii)} se d\'eduit du {\rm (i)} par passage \`a la limite, et pour le {\rm (i)} 
    il suffit de montrer la bijectivit\'e de l'application en question.

     Commençons par la surjectivit\'e. Le plongement ferm\'e 
     $\Pi^{(h)}\to ({\rm LA}^{(h)}(H))^r$ du th.~\ref{carduale} se dualise
     en une surjection $ {\cal D}_{h}(H)^r\to (\Pi^{(h)})^\dual $. Par construction, cette surjection 
     se factorise par la surjection ${\cal D}_{h}(H)^r\to {\cal D}_{h}(H)\otimes_{\Lambda(H)}\Pi^\dual $, ce qui permet de conclure.
     
    Pour d\'emontrer l'injectivit\'e, compte tenu
     du th\'eor\`eme de Hahn-Banach et de la dualit\'e de Schikhof \cite{ST1}, il suffit de prouver la surjectivit\'e de l'application transpos\'ee $\Pi^{(h)}\to ({\cal D}_{h}(H)\otimes_{\Lambda(H)}\Pi^\dual )^\dual $.
     L'application naturelle
  $\Pi^\dual \to{\cal D}_{h}(H)\otimes_{\Lambda_L} \Pi^\dual $ est continue ($\Pi^\dual $ est muni de la topologie faible), d'image dense, car 
  $\Lambda(H)\to {\cal D}_{h}(H)$ a ces propri\'et\'es. Cela montre que si $F\in ({\cal D}_{h}(H)\otimes_{\Lambda(H)}\Pi^\dual )^\dual $, alors il existe 
     un unique $v\in \Pi$ tel que $F(1\otimes l)=l(v)$ pour tout $l\in \Pi^\dual $. La continuit\'e de $F$ combin\'ee au fait que
          $p^{-|\alpha|r_h} b^{\alpha}$ tend vers $0$ dans ${\cal D}_{h}(H)$ pour la topologie faible 
          montrent que $p^{-r_h|\alpha|} l(b^{\alpha}v)$ tend vers $0$ pour tout $l\in \Pi^\dual $. On d\'eduit du
        th.~\ref{carduale} que $v\in \Pi^{(h)}$ et on conclut en utilisant la densit\'e de l'image de  $\Pi^\dual \to {\cal D}_{h}(H)\otimes_{\Lambda(H)} \Pi^\dual $.
        
        Pour d\'emontrer le {\rm (ii)},  posons $M={\cal D}(H)\otimes_{\Lambda(H)} \Pi^\dual $. Puisque 
        $\Pi^\dual $ est un $\Lambda(H)$-module de pr\'esentation finie, 
        $M$ est un ${\cal D}(H)$-module coadmissible~\cite{STInv}, et donc
        $M$ est isomorphe \`a la limite inverse des ${\cal D}_h(H)\otimes_{\Lambda(H)} \Pi^\dual $.
        Le {\rm (ii)} se d\'eduit donc de {\rm (i)} et de l'isomorphisme $\Pi^{\rm an}\simeq \varinjlim_{h} \Pi^{(h)}$ 
fourni par le th\'eor\`eme d'Amice.

    \begin{coro}\label{exac1}
     Le foncteur $\Pi\mapsto \Pi^{(h)}$ est exact de la cat\'egorie des $L$-repr\'esentations de Banach admissibles de $H$ dans la cat\'egorie
     des $L$-espaces de Banach.
    
    \end{coro}
    
    \demo
     C'est une cons\'equence directe du th\'eor\`eme de l'image ouverte, de la proposition pr\'ec\'edente et de la platitude de 
     ${\cal D}_{h}(H)$ sur $\Lambda(H)$ (prop. \ref{Schneider}).

\begin{rema} \label{exac2}
Il serait plus naturel d'\'etudier le foncteur $\Pi\mapsto \Pi_h$, o\`u 
$\Pi_h$ est l'espace des vecteurs $v\in \Pi$ tels que la fonction
$x \mapsto h^x\cdot v$ soit analytique sur $x_0+p^h\zp^d$ pour tout 
$x_0\in \zp^d$. La principale raison pour ne pas prendre ce point de vue est que ce foncteur
n'est pas exact, ce qui est d\'esagr\'eable pour les applications. En effet, si ${\rm LA}_h(H)={\cal C}(H)_h$, le th\'eor\`eme d'Amice~\cite{Am64} montre que 
$\phi\in {\rm LA}_h(H)$
si et seulement si
$$\lim_{|\alpha|\to\infty} v_p(a_{\alpha}(\phi))-\sum_{i=1}^dv_p\left(\left[\frac{\alpha_i}{p^h}\right]!\right)=\infty,$$
donc le dual topologique de ${\rm LA}_h(H)$ est l'alg\`ebre \`a puissances
divis\'ees partielles\footnote{Ses \'el\'ements sont de la forme
$\sum_{\alpha\in\N^d}c_{\alpha}\frac{1}{[\frac{\alpha}{p^h}]!}b^{\alpha}$
o\`u $(c_{\alpha})_{\alpha}$ est une suite born\'ee de $L$.}. Or cette alg\`ebre n'est pas plate
sur $\Lambda(H)$ si $d\geq 2$. On remarquera toutefois que ${\rm LA}_h(H)^\dual $  est en fait tr\`es
proche de ${\cal D}_{<h+1}(H)$ puisque $v_p([\frac{n}{p^h}]!)\sim n r_{h+1}$, et donc
$\Pi^{(h+1)}$ est tr\`es semblable \`a $\Pi_h$.

\end{rema}

\subsection{Coh\'erence}

On suppose \`a partir de maintenant que $G={\rm GL}_n(\qp)$
 et on note $K={\rm GL}_n(\zp)$ et $Z$ le centre de $G$.
 Soit
 $H=1+p^{\kappa} {\rm M}_n(\zp)$; c'est un pro-$p$-groupe uniforme de dimension $n^2$
 (rappelons que $\kappa=1$ si $p>2$ et $\kappa=2$ si $p=2$). 
Si $g\in G$, on note $d(g,1)$ le plus petit entier $m$
tel que $H^{p^m}\subset gHg^{-1}$. Puisque $KZ$ normalise 
$H$ et $H^{p^m}$, $d(g,1)$ ne d\'epend que de $KZgKZ$. 

\begin{prop}
Soit $\Pi$ une $L$-repr\'esentation de Banach de $G$. Alors $\Pi^{(h)}$ est stable par 
$K$ et $g\cdot \Pi^{(h)}\subset \Pi^{(h+n)}$ si $d(g,1)\leq n$.
\end{prop}

\demo
 La premi\`ere assertion se d\'emontre comme le cor.~\ref{Pihstable}, en utilisant le fait que 
 $H$ est distingu\'e dans $K$. La seconde d\'ecoule de la prop.~\ref{dcalage}.

\medskip

  D'apr\`es la d\'ecomposition de Cartan, on a $G=\sqcup_{t\in T^+} KtK$, o\`u $T^+$ est l'ensemble des matrices
  ${\rm diag}(p^{a_1},p^{a_2},...,p^{a_n})$, avec $a_1\geq a_2\geq...\geq a_n \in\mathbf{Z}$. Pour tout entier positif $m$, on note
  $T^+_m$ le sous-ensemble de $T^+$ form\'e des matrices ${\rm diag}(p^{a_1},p^{a_2},...,p^{a_n})$
 avec $a_1-a_n\leq m$. 

\begin{lemm} \label{conj}
  On a $d(g,1)\leq m$
 si et seulement si $g\in K T_m^+K$. En particulier, 
 $d(g,1)=0$ si et seulement si $g\in KZ$. 
\end{lemm}

\demo
 Soit $g=k_1tk_2$, avec $k_1,k_2\in K$ et $t={\rm diag}(p^{a_1},..., p^{a_n})\in T^+$. 
 Puisque $H$ et $H^{p^m}$ sont distingu\'es dans $K$, on a $g^{-1}H^{p^m}g\subset H$
 si et seulement si $ t^{-1}H^{p^m}t \subset H$. On conclut en utilisant les \'egalit\'es
 $t (x_{ij})_{i,j} t^{-1}=(p^{a_i-a_j} x_{ij})_{i,j}$ et $H^{p^m}=1+p^{m+\kappa} {\rm M}_n(\zp)$.

\begin{lemm}\label{dg1}
On a $d(g,1)\leq l+m$ si et seulement si l'on peut \'ecrire $g=g_1g_2$, avec $d(g_1,1)\leq l$ et $d(g_2,1)\leq m$.

\end{lemm}

\demo C'est une cons\'equence directe du lemme \ref{conj} et du fait que 
$T_l^+\cdot T_m^+=T_{l+m}^+$.

\medskip
  Si $W$ est un sous-$L[KZ]$-module d'un $L[G]$-module $\Pi$
  et si $h\in\N$, on note $$W^{[h]}=\sum_{d(g,1)\leq h} g\cdot W;$$
c'est un sous-$L[KZ]$-module de $\Pi$. Notons que $g\cdot W$ et $d(g,1)$ ne d\'ependent
que de l'image de $g$ dans $S:=G/KZ$. 

\begin{lemm}
 L'ensemble $\{s\in S,\  d(s,1)\leq h\}$ est fini.
\end{lemm}

\demo Soit $T_{h,0}^+$ l'ensemble des
$t={\rm diag}(p^{a_1},...,p^{a_n})\in T_h^+$ tels que $a_n=0$. 
Alors $T_{h,0}^+$ est un ensemble fini et $T_{h,0}^+Z\supset T_h^+$. 
Si $I_h$ est un syst\`eme de repr\'esentants de $K/K_h$ (avec $K_h=1+p^h {\rm M}_n(\zp)$), alors
$$KT_h^+K\subset \bigcup_{t\in T_{h,0}^+, k\in I_h} k K_h t KZ\subset \bigcup_{t\in T_{h,0}^+, k\in I_h} k t KZ,$$
la derni\`ere inclusion \'etant une cons\'equence du fait que $t^{-1}K_h t\subset K$ pour tout $t\in T_{h,0}^+$.
On conclut en utilisant le lemme \ref{conj}.

\medskip
   On suppose maintenant que $\Pi$ est une $L$-repr\'esentation de Banach admissible de~$G$, ayant un caract\`ere central.
   Le sous-espace $\Pi^{(h)}$ de $\Pi$ (d\'efini en consid\'erant $\Pi$ comme une repr\'esentation du 
   pro-$p$-groupe uniforme $H$) est stable par $KZ$, puisque $KZ$ normalise $H$. 

\begin{defi} On dit que la repr\'esentation
$\Pi$ est {\it coh\'erente} s'il existe
$m(\Pi)$ tel que $\Pi^{(h+k)}=(\Pi^{(h)})^{[k]}$, pour tous $h\geq m(\Pi)$ et $k\in\N$.

\end{defi}

\begin{rema}
(i) L'inclusion $(\Pi^{(h)})^{[k]}\subset \Pi^{(h+k)}$
est vraie pour n'importe quelle repr\'esentation de Banach $\Pi$ (cela d\'ecoule de la prop. \ref{dcalage}).

(ii) Le lemme \ref{dg1} montre que $W^{[k+1]}=(W^{[k]})^{[1]}$.
Pour montrer que $\Pi$ est coh\'erente, il suffit donc de v\'erifier que
$\Pi^{(h+1)}=(\Pi^{(h)})^{[1]}$
pour tout $h$ assez grand.

(iii) Si $\Pi$ est coh\'erente,
alors $\Pi^{(h)}$ engendre $\Pi^{\rm an}$, en tant que repr\'esentation
de~$G$, pour tout~$h$ assez grand (cela d\'ecoule de ce que $\Pi^{\rm an}$ est la r\'eunion des 
$\Pi^{(h)}$). 
On peut se demander si une propri\'et\'e de ce genre
est automatiquement v\'erifi\'ee
pour une repr\'esentation de Banach admissible de $G$
(au moins dans le cas de longueur finie). C'est le cas pour 
$G={\rm GL}_2(\qp)$ (cf. th.~\ref{COH}).
\end{rema}

\begin{prop} \label{extcoh}
   Si $0\to\Pi_1\to\Pi\to\Pi_2\to 0$ est une suite
exacte de repr\'esentations de Banach admissibles
de $G$ et si $\Pi_1$ est coh\'erente, alors $\Pi$ est coh\'erente
si et seulement si $\Pi_2$ l'est.
\end{prop}

\demo
 C'est une cons\'equence de l'exactitude du foncteur $\Pi\mapsto \Pi^{(h)}$.

\section{Vecteurs analytiques des repr\'esentations unitaires de ${\bf GL}_2(\Q_p)$} \label{vectan}

Ce chapitre \'etend (et raffine) \`a toutes les paires $G$-compatibles les r\'esultats de \cite[chap. V]{Cbigone}, concernant l'espace $\Pi_{\delta}(D)^{\rm an}$ des vecteurs localement analytiques de la repr\'esentation $\Pi_{\delta}(D)$. 
L'approche est assez diff\'erente de celle de~\cite{Cbigone} m\^eme si
le noyau technique (\`a savoir le \S~\ref{actiondagger} et, en particulier,
la prop.~\ref{crucial}) est le m\^eme.

   \subsection{Pr\'eliminaires}   \label{prelimin}

   \quad \quad On fixe dans la suite une paire $G$-compatible $(D,\delta)$, avec $D\in\fget({\cal E})$, et un r\'eseau $D_0$ de $D$, stable par $\varphi$ et $\Gamma$. On note
      $\Pi=\Pi_{\delta}(D)$ et $\Pi_0=\Pi_{\delta}(D_0)$. Alors
   $\Pi_0$ est un r\'eseau de $\Pi$, ouvert, born\'e et stable par $G$. On munit $\Pi$ de la valuation
   $v_{\Pi}$, \`a valeurs dans $v_p(L)$, faisant de $\Pi_0$ la boule unit\'e de $\Pi$.

   On renvoie au chap.~\ref{Fctan}
   pour les anneaux de fonctions analytiques utilis\'es dans la suite. Rappelons que $\Lambda(\Gamma)$ est l'anneau des mesures \`a valeurs dans $\O_L$ sur $\Gamma$. Si $R$ est un anneau de s\'eries de Laurent (comme $\oe $, ${\cal R}$, $\oe^{\dagger,b}$, etc.),
on peut remplacer la variable $T$ par
$\gamma-1$ o\`u $\gamma$ est un g\'en\'erateur topologique de $\Gamma$ (ou plut\^ot de
l'image inverse de $1+p\Z_p$ dans $\Gamma$) pour construire un anneau
$R(\Gamma)$ (on renvoie le lecteur au \no 3 du
\S~V.1 de \cite{Cbigone} pour les d\'etails).

   Rappelons que $m(D_0)$ est un entier assez grand, qui ne
d\'epend que de $D_0$. Il est en particulier choisi tel que la prop.~\ref{CCF} s'applique \`a $D_0$, et donc $D_0^{\natural}\boxtimes_{\delta}\p1
\subset D_0^{(0,r_{m(D_0)}]}\boxtimes_{\delta}\p1$ (cf. \cite[cor. II.7.2]{Cmirab}). Comme $D_0^{\natural}$ est compact,
il existe $l_1=l_1(D_0)$ tel que $D_0^{\natural}\subset T^{-l_1}D_0^{\dagger, m(D_0)}$, donc $D_0^{\natural}\boxtimes_{\delta}\p1
\subset (T^{-l_1}D_0^{\dagger, m(D_0)})\boxtimes_{\delta}\p1$, o\`{u} l'on note $$X\boxtimes_{\delta}\p1=(D\boxtimes_{\delta}\p1)\cap (X\times X)$$
pour $X\in \{D^{(0,r_b]}, D^{\dagger}\}$ (pareil avec $D_0$ si $X\in \{T^a D_0^{\dagger,b}, D_0^{(0,r_b]}\}$).
 Quitte \`a augmenter encore $l_1$ et $m(D_0)$, on peut supposer qu'ils
    sont aussi associ\'es \`a $(\check{D_0},\delta^{-1})$, et que l'involution $i_{\delta}$
    de $\O_L[\Gamma]$ qui envoie $\sigma_a$ sur $\delta(a)\sigma_{\frac{1}{a}}$ s'\'etend en une involution continue de
    $\oe^{\dagger,b}(\Gamma)$, $\oe^{(0,r_b]}(\Gamma)$, ${\cal R}(\Gamma)$ pour $b\geq m(D_0)$ (voir \cite[lemme V.2.3]{Cbigone}).

\subsection{Vecteurs analytiques et surconvergence}\label{coefM}

   \quad Pour $m\geq 2$ on note $K_m$ le sous-groupe $1+p^m{\rm M}_2(\zp)$ de $G$,
   $$a_m^+=\left(\begin{smallmatrix} 1+p^m & 0\\0 & 1\end{smallmatrix}\right), \quad a_m^-=\left(\begin{smallmatrix} 1 & 0\\0 & 1+p^m\end{smallmatrix}\right), \quad u_m^{+}=\left(\begin{smallmatrix} 1 & p^m\\0 & 1\end{smallmatrix}\right), \quad
   u_m^{-}=\left(\begin{smallmatrix} 1 & 0\\p^m & 1\end{smallmatrix}\right)$$ et, pour $\alpha\in\mathbf{N}^4$, on note
   $$b_m^{\alpha}=(a_m^+-1)^{\alpha_1}\cdot (a_m^{-}-1)^{\alpha_2}\cdot (u_m^{+}-1)^{\alpha_3}\cdot (u_m^{-}-1)^{\alpha_4}\in \zp[K_m].$$

   Alors $(a_m^+,a_m^{-}, u_m^+, u_m^-)$ est un syst\`eme minimal de g\'en\'erateurs topologiques du pro-$p$ groupe uniforme $K_m$.

  \begin{defi} \label{pimb} Si $b\geq m+1$, soit
    $$\Pi^{(b)}=\{v\in \Pi,\  \lim_{|\alpha|\to\infty} v_{\Pi}(b_m^{\alpha}\cdot v)-p^m r_b\cdot |\alpha|=\infty\},$$
que l'on munit de la valuation $$v^{(b)}(v)=\inf_{\alpha\in\mathbf{N}^4} (v_{\Pi}(b_m^{\alpha}\cdot v)-p^m r_b\cdot |\alpha|),$$
pour laquelle c'est un espace de Banach.
\end{defi}

\begin{rema} \label{Amice}
{\rm (i)}
 Le th\'eor\`eme d'Amice~\cite{Am64} montre que pour tout $m$ on a un isomorphisme d'espaces vectoriels topologiques
   $$\varinjlim_{b} \Pi^{(b)}\simeq \Pi^{\rm an}.$$

{\rm (ii)}
Il r\'esulte de la prop.~\ref{dcalage} que l'espace $\Pi^{(b)}$ ne d\'epend pas du choix de $m\leq b-1$
(la valuation $v^{(b)}$ en d\'epend, mais les valuations obtenues en faisant varier $m$ sont toutes \'equivalentes).
\end{rema}

Le th.~\ref{mainanal} ci-dessous d\'ecrit l'espace de Banach $\Pi^{(b)}$ en fonction de $D$.
Nous commen\c{c}ons par pr\'eciser un peu les topologies sur les espaces divers et vari\'es apparaissant dans ce th\'eor\`eme.
On dispose d'une pl\'eiade d'anneaux de s\'eries de Laurent (voir le chap.~\ref{Fctan}), chacun ayant une topologie naturelle.
  Cela induit une topologie naturelle sur les modules libres de type fini sur ces anneaux
  (et qui ne d\'epend pas des choix de bases). En appliquant cette discussion aux modules $D_0^{\dagger,b}$, $D_0^{(0,r_b]}$, $D^{(0,r_b]}$,
  $D^{\dagger}$, $D_{\rm rig}$, etc, on obtient des topologies sur ces modules.

  Si $X\in \{D_0^{(0,r_b]}, D^{(0,r_b]}, D^{\dagger}\}$,
  on munit $X\boxtimes_{\delta}\p1$ de la topologie induite par l'inclusion $X\boxtimes_{\delta}\p1\subset X\times X$.
  Cette topologie est plus forte que celle induite par l'inclusion $X\boxtimes_{\delta}\p1\subset D_0\boxtimes_{\delta}\p1$
  (ou $X\subset D\boxtimes_{\delta}\p1$). Comme $D_0^{\natural}\boxtimes_{\delta}\p1$ est ferm\'e dans $D_0\boxtimes_{\delta}\p1$, il
  est aussi ferm\'e dans $D_0^{(0,r_b]}\boxtimes_{\delta}\p1$ pour $b>m(D_0)$, donc $D^{\natural}\boxtimes_{\delta}\p1$
  est ferm\'e dans $D^{(0,r_b]}\boxtimes_{\delta}\p1$ et dans $D^{\dagger}\boxtimes_{\delta}\p1$.
  On munit alors $(D_0^{(0,r_b]}\boxtimes_{\delta}\p1)/(D_0^{\natural}\boxtimes_{\delta}\p1)$ et $(X\boxtimes_{\delta}\p1)/(D^{\natural}\boxtimes_{\delta}\p1)$ de la topologie quotient,
  pour $X\in \{D^{(0,r_b]}, D^{\dagger}\}$.
  
\medskip
     On fixe une paire $G$-compatible $(D,\delta)$, avec $D\in \fget({\cal E})$ et on note $\Pi=\Pi_{\delta}(D)$
     et $\check{\Pi}=\Pi_{\delta^{-1}}(\check{D})$. 

   \begin{theo} \label{mainanal}
       Il existe $c=c(D,\delta)$ tel que pour tout $b\geq c$:

    {\rm (i)} Le sous-module $D^{(0,r_b]}\boxtimes_{\delta}\p1$ de
    $D\boxtimes_{\delta}\p1$ est stable par ${\rm GL}_2(\zp)$, qui agit contin\^{u}ment.

   {\rm (ii)} On a un isomorphisme canonique de ${\rm GL}_2(\zp)$-modules de Banach
      $$(D^{(0,r_b]}\boxtimes_{\delta}\p1)/(D^{\natural}\boxtimes_{\delta}\p1)\simeq \Pi^{(b)},$$
      et donc une suite exacte de ${\rm GL}_2(\zp)$-modules topologiques
       $$0\to \check{\Pi}^\dual \to D^{(0,r_b]}\boxtimes_{\delta}\p1\to \Pi^{(b)}\to 0.$$
   \end{theo}

    Avant de passer \`a la preuve du th.~\ref{mainanal}, qui occupe les \S\S~V.3-V.6, donnons-en quelques cons\'equences.
    Le r\'esultat suivant d\'ecoule formellement du th.~\ref{mainanal} et de la rem.~\ref{Amice}. 

 \begin{coro}\label{anali} 
 
 {\rm (i)} Le sous-module $D^{\dagger}\boxtimes_{\delta}\p1$ de $D\boxtimes_{\delta}\p1$ est stable par $G$, qui agit contin\^{u}ment.

 {\rm (ii)} On a un isomorphisme canonique de
   $G$-modules topologiques $$ (D^{\dagger}\boxtimes_{\delta}\p1)/(D^{\natural}\boxtimes_{\delta}\p1)\simeq \Pi^{\rm an},$$
et donc une suite exacte de $G$-modules topologiques
 $$0\to \check{\Pi}^\dual \to D^{\dagger}\boxtimes_{\delta}\p1\to \Pi^{\rm an}\to 0.$$
  \end{coro}

Le th.~\ref{mainanal} fournit un raffinement du th\'eor\`eme de
  Schneider et Teitelbaum pour les objets de ${\rm Rep}(\delta)$.

\begin{coro}\label{stdense}
 Si $\Pi\in {\rm Rep}_L(G)$ il existe
  $m_0\geq 2$ tel que $\Pi^{(b)}$ soit dense dans~$\Pi^{\rm an}$ 
(et donc aussi dans $\Pi$) pour tout $b\geq m_0+1$. 
\end{coro}

\demo Disons que $\Pi$ est {\it bonne} 
si elle satisfait le corollaire. 
Trivialement, toute repr\'esentation de dimension finie est bonne.
 Si $\Pi=\Pi_{\delta}(D)$ pour une paire $G$-compatible
$(D,\delta)$, alors $\Pi$ est bonne: cela d\'ecoule
du th.~\ref{mainanal}, du cor.~\ref{anali} et de  
la densit\'e de $D^{(0,r_b]}\boxtimes_{\delta}\p1$ dans $D^{\dagger}\boxtimes_{\delta}\p1$ (qui d\'ecoule 
de celle de $D^{(0,r_b]}$ dans $D^{\dagger}$, elle-m\^eme cons\'equence de la densit\'e de ${\cal E}^{(0,r_b]}$ dans
  ${\cal E}^{\dagger}$). 
  
   En utilisant le th.~\ref{recoverPi}, le cor.~\ref{dimfinsl} et ce que l'on vient de d\'emontrer, nous pouvons conclure dans le cas g\'en\'eral
   gr\^ace au lemme suivant qui
est une cons\'equence de l'exactitude des foncteurs 
        $\Pi\mapsto \Pi^{(b)}$ et $\Pi\mapsto \Pi^{\rm an}$
(prop.~\ref{exac1}).
   
   \begin{lemm}
    Soit $0\to \Pi_1\to \Pi\to \Pi_2\to 0$ une suite exacte dans ${\rm Rep}_L(G)$, avec $\Pi_1$ bonne. 
    Si une des $\Pi$ et $\Pi_2$ est bonne, alors l'autre l'est aussi.
   \end{lemm}
        
   Le corollaire suivant est un sous-produit de la preuve du th.~\ref{mainanal}.

\begin{coro}\label{uwu}
   Soit $\Pi\in {\rm Rep}_L(G)$. Si 
  $x\mapsto \left(\begin{smallmatrix} 1 & x \\0 & 1\end{smallmatrix}\right)v$ et
  $x\mapsto \left(\begin{smallmatrix} 1 & 0 \\x & 1\end{smallmatrix}\right)v$ sont localement
  analytiquesi sur $\Q_p$, alors $v\in \Pi^{\rm an}$.  
\end{coro}

\demo  Commençons par le cas o\`u $\Pi=\Pi_{\delta}(D)$ pour une paire $G$-compatible $(D,\delta)$. La preuve de la prop.~\ref{relevanal} ci-dessous n'utilise que la croissance des coefficients
 de Mahler des applications $x\to \left(\begin{smallmatrix} 1 & x \\0 & 1\end{smallmatrix}\right)v$ et
  $x\to \left(\begin{smallmatrix} 1 & 0 \\x & 1\end{smallmatrix}\right)v$. Elle montre
  que $v$ admet un rel\`evement \`a $D^{\dagger}\boxtimes_{\delta}\p1$. Le corollaire
  \ref{anali} permet alors de conclure.
  
     Passons au cas g\'en\'eral. Disons que $v\in \Pi$ est 
{\it presqu'analytique} si les applications
  $x\mapsto \left(\begin{smallmatrix} 1 & x \\0 & 1\end{smallmatrix}\right)v$ et
  $x\mapsto \left(\begin{smallmatrix} 1 & 0 \\x & 1\end{smallmatrix}\right)v$ sont localement
  analytiques, et que $\Pi$ est {\it bonne} 
si tout vecteur presqu'analytique est localement analytique 
  (notons que tout vecteur localement analytique est trivialement presqu'analytique). On d\'eduit du th\'eor\`eme
  \ref{recoverPi}, du cor.~\ref{dimfinsl} et du premier paragraphe l'existence d'un morphisme $\beta_{\Pi}: \Pi_1\to \Pi/ \Pi^{{\rm SL}_2(\qp)}$, dont 
  le noyau et le conoyau sont de dimension finie, et tel que $\Pi_1$ soit bonne. On conclut en utilisant le lemme suivant: 
  
  \begin{lemm}
  
   Soit $0\to \Pi_1\to \Pi\to \Pi_2\to 0$ une suite exacte dans ${\rm Rep}_L(G)$. 
   
   {\rm (i)} Si $\Pi_1$ est de dimension finie et si $\Pi$ ou $\Pi_2$ est bonne, alors 
   l'autre l'est aussi.
   
   {\rm (ii)} Si $\Pi_2$ est de dimension finie et si $\Pi_1$ est bonne, alors $\Pi$ est bonne.

  \end{lemm}

   \demo {\rm (i)} Supposons d'abord que $\Pi$ est bonne et soit $v\in \Pi_2$ presqu'analytique. Soit 
   $\hat{v}\in \Pi$ un rel\`evement de $v$. Notons $T=\left(\begin{smallmatrix} 1 & 1 \\0 & 1\end{smallmatrix}\right)-1\in \O_L[G]$.
   Comme $v$ est presqu'analytique, il existe $r>0$ tel que $p^{-rn} T^n(v)$ tende vers $0$ dans $\Pi_2$, ce qui veut dire qu'il existe
   $x_n\in \Pi_1$ tels que $p^{-rn}T^n(\hat v)-x_n\to 0$ dans $\Pi$. Puisque $\Pi_1$ est de dimension finie, elle est tu\'ee par $T$ 
 (lemme \ref{unitdimfin}) et donc $p^{-rn}T^{n+1}(\hat v)\to 0$ dans $\Pi$. Le th\'eor\`eme d'Amice entra\^\i ne que 
 $x\mapsto \left(\begin{smallmatrix} 1 & x \\0 & 1\end{smallmatrix}\right)\hat v$ est localement analytique. On obtient de m\^eme que 
 $x\mapsto \left(\begin{smallmatrix} 1 & 0 \\x & 1\end{smallmatrix}\right)\hat v$ est localement analytique, et donc $\hat v$ est localement analytique
 (car $\Pi$ est bonne) et son image $v$ dans $\Pi_2$ l'est aussi. Cela montre que $\Pi_2$ est bonne.

Le reste de l'\'enonc\'e est une cons\'equence de l'exactitude du foncteur $\Pi\mapsto\Pi^{\rm an}$.

\subsection{Rel\`evement \`a $D^{(0,r_b]}\boxtimes_{\delta}\p1$}

  La proposition ci-dessous montre que tout $v\in \Pi^{(b)}$ se rel\`eve \`a $D^{(0,r_b]}\boxtimes_{\delta}\p1$.
   Rappelons que $l_1$ est choisi tel que $D_0^{\natural}\subset T^{-l_1}D_0^{\dagger, b}$ (cf.~\S~\ref{prelimin}).

\begin{prop}\label{relevanal}
 Soient $b>m> m(D_0)$ et $v\in \Pi^{(b)}$. Alors $v$
 admet un rel\`evement \`a
 $(p^{s}T^{s'}D_0^{\dagger,b})\boxtimes_{\delta}\p1\subset D^{(0,r_b]}\boxtimes_{\delta}\p1$, avec
 $s=v^{(b)}(v)$ et $s'=-(p^m n_b+l_1)$.
\end{prop}

\demo On peut supposer que $v^{(b)}(v)\geq 0$, quitte \`a multiplier $v$ par une puissance de $p$. 
 Notons $\xi=(u_m^+)^{n_b}$, $\mu=(u^{-}_m)^{n_b}$ et
     $$a_k=\min( v_{\Pi}(\xi^k v), v_{\Pi}(\mu^k v)),$$
     de telle sorte que $\lim_{k\to\infty} a_k-p^m k=\infty$ et $a_k\geq p^m k+v^{(b)}(v)$ pour tout $k$ (par d\'efinition de
     $\Pi^{(b)}$ et
     $v^{(b)}$).

      Soient $X=D_0\boxtimes_{\delta}\p1$ et
     $Y=D_0^{\natural}\boxtimes_{\delta}\p1$, de telle sorte que $\Pi_0=X/Y$ est la boule unit\'e de 
     $\Pi$ pour la valuation $v_{\Pi}$. Alors $p^{a_0}\Pi_0\subset \Pi_0$, car $a_0\geq 0$, et 
      $v\in p^{a_0} \Pi_0$, puisque 
     $v_{\Pi}(v)\geq a_0$. Ainsi, $v$ poss\`ede un 
     rel\`evement $z=(z_1,z_2)\in p^{a_0}X$. Nous aurons besoin du lemme suivant.

     \begin{lemm}\label{suitexiste}
      Il existe une suite $(y_n)_{n\geq 0}$ d'\'el\'ements de $Y$ telle que pour tout $n\geq 1$
      $$\xi^n z-\sum_{k=0}^{n-1} p^{a_k}\xi^{n-k-1} y_k\in p^{a_n}X.$$

     \end{lemm}

     \demo On construit la suite en question par r\'ecurrence. Noter que $X$ et $Y$ sont stables par $\xi$ et $\eta$ (car ils sont stables par $G$), ainsi que
     l'\'egalit\'e $p^{n}X\cap Y=p^{n}Y$ (pour $n\geq 0$), car $Y$ est un sous-$\O_L$-module satur\'e de $X$ par d\'efinition.

     Supposons d'abord que $n=1$. Si $a_1<a_0$, on prend $y_0=0$, supposons donc que $a_1\geq a_0$.
      Comme $v_{\Pi}(\xi v)\geq a_1\geq 0$, on a $\xi z\in (p^{a_1}X+Y)\cap p^{a_0}X\subset p^{a_1}X+p^{a_0}Y$, ce qui
      montre l'existence de $y_0$. Supposons avoir trouv\'e $y_0,...,y_{n-1}$ et \'ecrivons
      $$\xi^n z=\sum_{k=0}^{n-1} p^{a_k}\xi^{n-k-1} y_k+p^{a_n} u$$ pour un $u\in X$.
      Si $a_n>a_{n+1}$, on prend
      $y_{n}=0$. Sinon, en appliquant $\xi$ \`a l'\'egalit\'e pr\'ec\'edente on obtient
      $p^{a_n}\xi u\in \xi^{n+1}z+Y\subset p^{a_{n+1}}X+Y$ (la deuxi\`eme inclusion suit de
      $v_{\Pi}(\xi^{n+1} v)\geq a_{n+1}$). On en d\'eduit que $p^{a_n}\xi u\in p^{a_{n+1}}X+p^{a_n}Y$
      et on choisit $y_n\in Y$ tel que $p^{a_n}\xi u-p^{a_n}y_n\in p^{a_{n+1}}X$. Cela permet de conclure.

      Revenons \`a la preuve de la proposition. En appliquant ${\rm Res}_{\zp}$ \`a la relation
      du lemme \ref{suitexiste}, et en utilisant le fait que $\xi$ agit par multiplication
     par $\varphi^m(T)^{n_b}$, on obtient $z_1-\sum_{k=0}^{n-1} A_k\in p^{a_n} D_0$, avec
       $$A_k=\frac{p^{a_k}}{\varphi^{m}(T)^{(k+1)n_b}}{\rm Res}_{\zp}(y_k)=p^{a_k-p^mk}\frac{1}{\varphi^{m}(T)^{n_b}} \left(\frac{p^{p^m}}{\varphi^m(T)^{n_b}}\right)^{k}{\rm Res}_{\zp}(y_k)$$ 
       
    Le lemme \ref{hyper} montre que $\frac{1}{\varphi^m(T)^{n_b}}\in T^{-p^m n_b}\oe^{\dagger,b}$. Puisque
     $p\in T^{n_b}\oe^{\dagger,b}$, on en d\'eduit que 
      $\frac{p^{p^m}}{\varphi^m(T)^{n_b}}\in \oe^{\dagger,b}$. En utilisant aussi le fait que 
      ${\rm Res}_{\zp}(y_k)\in D_0^{\natural}\subset T^{-l_1}D_0^{\dagger,b}$, on obtient enfin 
      $$A_k\in p^{a_k-p^m k} T^{-l_1-p^mn_b}D_0^{\dagger,b}.$$
     
      Puisque $a_k-p^m k$ tend vers $\infty$ et est minor\'e par $v^{(b)}(v)$, et puisque 
    $D_0^{\dagger,b}$ est complet pour la topologie $p$-adique, le paragraphe pr\'ec\'edent montre que 
    $\sum_{k\geq 0} A_k$ converge dans $p^{v^{(b)}(v)}T^{-l_1-p^m n_b} D_0^{\dagger,b}$.
    La relation $z_1-\sum_{k=0}^{n-1} A_k\in p^{a_n} D_0$ montre que $\sum_{k\geq 0} A_k$ tend vers 
    $z_1$ dans $\oe $. On en d\'eduit que $\sum_{k\geq 0} A_k$ tend vers $z_1$ dans 
    $p^{v^{(b)}(v)}T^{-l_1-p^m n_b} D_0^{\dagger,b}$, et donc que $z_1\in p^{v^{(b)}(v)}T^{-l_1-p^m n_b} D_0^{\dagger,b}$. Les m\^emes arguments
       [remplacer dans ce qui pr\'ec\`ede $\xi$ par $\mu$ et ${\rm Res}_{\zp}$ par ${\rm Res}_{\zp}(w\cdot)$] donnent
        la m\^eme estim\'ee pour $z_2$, ce qui permet de conclure.

\begin{rema}

 Suppposons que $D^{\natural}=D^{\sharp}$ (c'est par exemple le cas si $D$ est 
irr\'eductible de dimension $\geq 2$), de telle sorte que l'inclusion de $\tilde{D}$ dans $D\boxtimes_\delta\p1$ induise un
 isomorphisme de $B$-modules de Banach $\Pi\simeq \tilde{D}/\tilde{D}^+$ (cor.~\ref{lienKir}).
 En posant $X=\tilde{D_0}$ et $Y=\tilde{D}_0^+$, on v\'erifie sans mal que le lemme \ref{suitexiste}
 s'applique encore (le point est que $X$ et $Y$ sont stables par $\xi$ et $\mu$, et $Y\cap p^{n}X=p^n Y$ pour $n\geq 0$).
 Le reste de la preuve s'applique et montre que tout $v\in \Pi^{(b)}$ admet un rel\`evement \`a
 $\tilde{D}^{(0,r_b]}$. Par contre, l'image dans $\Pi$ d'un \'el\'ement de 
$\tilde{D}^{(0,r_b]}$ n'est pas toujours localement analytique.

\end{rema}

  Il nous reste \`a montrer que l'image de $D^{(0,r_b]}\boxtimes_{\delta}\p1$ dans $\Pi$ est contenue
   dans~$\Pi^{(b)}$, si $m$ et $b$ sont assez grands. Cela va demander un certain nombre d'estim\'ees techniques, auxquelles sont
   d\'edi\'ees les parties \ref{vectpropres} et \ref{actiondagger} ci-dessous.

 \subsection{Vecteurs propres de $\psi$}\label{vectpropres}

\quad Si $\alpha\in \O_L^\dual $, on pose
  $${\cal C}^{\alpha}=(1-\alpha\cdot \varphi)D_0^{\psi=\alpha}\subset D_0\boxtimes_{\delta}\zpet.$$

   La proposition suivante est une version de la prop.~V.2.1\footnote{Cette proposition n'est vraie qu'apr\`es tensorisation par $L$, le probl\`eme \'etant que
    $D_0^{\psi=1}$ n'est pas toujours contenu dans $D_0^{\natural}$, m\^eme sous les hypoth\`eses de loc.cit. Comme le montre la suite, cela ne
    change rien aux arguments.} de \cite{Cbigone}. Voir \ref{prelimin} pour les objets $i_{\delta}$, $\Lambda(\Gamma)$, etc.

\begin{prop}\label{psi1}
 Soit $P\in \O_L[X]$ tel que $P(\psi)=0$ sur $D_0^{\sharp}/D_0^{\natural}$. Soit
$\alpha\in \O_L^\dual $ tel que $\alpha$ et $\beta:=(\delta(p)\alpha)^{-1}$ ne soient pas des racines de $P$ et $\alpha^{-1}$ et $\beta^{-1}$ ne soient pas des valeurs propres
de $\varphi$ sur $D^{\rm nr}$. Alors $w_{\delta}({\cal C}^{\alpha})\cap {\cal C}^{\beta}$
est d'indice fini dans ${\cal C}^{\beta}$.

\end{prop}

\demo  On commence par montrer que $w_{\delta}({\cal C}^{\alpha})\otimes_{\O_L} L={\cal C}^{\beta}\otimes_{\O_L} L$.
Soit $z\in D_0^{\psi=\alpha}$ et soit $z'=P(\alpha)z=P(\psi)z$. Comme $D_0^{\psi=\alpha}\subset D_0^{\sharp}$ 
et $P(\psi)D_0^{\sharp}\subset D_0^{\natural}$ par hypoth\`ese, 
 on a
$z'\in D_0^{\natural}$, donc $(\alpha^{-n}z')_{n\geq 0}\in D_0^{\natural}\boxtimes_{\delta}\qp$.
  Comme $D_0^{\natural}\boxtimes_{\delta}\p1$ se surjecte sur $D_0^{\natural}\boxtimes_{\delta}\qp$ (lemme \ref{resqp}), il existe
  $x=(x_1,x_2)\in D_0^{\natural}\boxtimes_{\delta}\p1$ tel que ${\rm Res}_{\zp} \left( \left(\begin{smallmatrix} p^n & 0\\0 & 1\end{smallmatrix}\right)x  \right)
=\alpha^{-n}z'$ pour tout $n\geq 0$. Alors (ibid.) $\left(\begin{smallmatrix} p & 0\\0 & 1\end{smallmatrix}\right)x-\alpha^{-1}x\in
{\rm Ker} ({\rm Res}_{\qp})=(0,D_0^{\rm nr})$ et un petit calcul montre que ceci entra\^{i}ne $\psi(x_2)-\beta x_2\in D_0^{\rm nr}$.
Comme $\beta^{-1}$ n'est pas valeur propre de $\varphi\in {\rm End}_{L} (D^{\rm nr})$, il existe
$u\in D^{\rm nr}$ tel que $\beta\varphi(u)-u=\psi(x_2)-\beta x_2$. Alors
$x_2+\varphi(u)\in D_0^{\psi=\beta}\otimes_{\O_L} L$ et donc ${\rm Res}_{\zpet}(x_2)={\rm Res}_{\zpet} (x_2+\varphi(u))\in
{\cal C}^{\beta}\otimes_{\O_L} L$. Comme $${\rm Res}_{\zpet} (x_2)=w_{\delta} ({\rm Res}_{\zpet} (x_1))=
P(\alpha) w_{\delta} ((1-\alpha\varphi)(z))$$ et $P(\alpha)\ne 0$, on conclut que
$w_{\delta}({\cal C}^{\alpha}\otimes_{\O_L} L)\subset {\cal C}^{\beta}\otimes_{\O_L} L$. Par sym\'etrie et puisque $w_{\delta}$ est une involution,
 cette inclusion est une \'egalit\'e.

   Comme ${\cal C}^{\alpha}$ et ${\cal C}^{\beta}$ sont des $\Lambda(\Gamma)$-modules de type fini \cite[cor. VI.1.3]{Cmirab}
 et comme $w_{\delta}$ est $i_{\delta}$-semi-lin\'eaire, le paragraphe pr\'ec\'edent montre l'existence d'une constante $c=c(P,\alpha,D_0)$ telle que $w_{\delta}({\cal C}^{\beta})
 \subset p^{-c}{\cal C}^{\alpha}$. Soit alors $x\in{\cal C}^{\beta}$. On vient de voir
 qu'il existe $y\in D_0^{\psi=\alpha}$ tel que $w_{\delta}(x)=p^{-c}(1-\alpha\varphi)y$. Si $y'$ est un autre
 choix, alors $y-y'\in D_0^{\varphi=\frac{1}{\alpha}}$. Comme $w_{\delta}(x)\in D_0$, on a $y\pmod {p^c}\in
 (D_0/p^c D_0)^{\varphi=\frac{1}{\alpha}}$. De plus, si $y\pmod {p^c}=0$, alors $x\in w_{\delta}({\cal C}^{\alpha})
 \cap {\cal C}^{\beta}$. Ainsi, l'application $x\to y\pmod {p^c}$ induit une injection de
 ${\cal C}^{\beta}/({\cal C}^{\beta}\cap w_{\delta}({\cal C}^{\alpha}))$ dans le quotient
 de $(D_0/p^c D_0)^{\varphi=\frac{1}{\alpha}}$ par l'image de $D_0^{\varphi=\frac{1}{\alpha}}$. Comme ce quotient est fini (car
 $(D_0/p^c D_0)^{\varphi=\frac{1}{\alpha}}$ est contenu dans $(D_0/p^c D_0)^{\rm nr}$, qui est fini),
 cela permet de conclure.

 \begin{rema} Puisque $\psi$ est un endomorphisme du $O_L$-module de type fini $D_0^{\sharp}/D_0^{\natural}$, on peut toujours trouver un polyn\^ome non nul $P$ comme dans la proposition pr\'ec\'edente. 

  \end{rema}

 \subsection{L'action de $G$ sur $T^a D_0^{\dagger,b}\boxtimes_{\delta}\p1$} \label{actiondagger}

   \quad Le but de cette partie est de contr\^{o}ler l'action de $G$ sur les modules
   $(T^a D_0^{\dagger,b})\boxtimes_{\delta}\p1$, plus pr\'ecis\'ement de d\'emontrer 
la prop.~\ref{action} ci-dessous. La plupart des arguments sont adapt\'es de \cite[chap. V]{Cbigone}. 
       On fixe une base $e_1,e_2,...,e_d$ de $D_0^{\dagger, m(D_0)}$ sur $\oe^{\dagger, m(D_0)}$ (c'est aussi une base de $D_0^{\dagger,b}$ sur $\oe^{\dagger,b}$ pour tout $b\geq m(D_0)$).
       Les constantes
$c,c_1,c_2,...,m_1,m_2,...$ ci-dessous ne d\'ependent que de $D_0$, $\delta$ et du choix de la base $e_1,e_2,...,e_d$.

\begin{prop}\label{engendre}
 Il existe $m_1\geq m(D_0)$ tel que $w_{\delta}$ laisse stable $D_0^{(0,r_b]}\boxtimes_{\delta}\zpet$ pour
 tout $b\geq m_1$.
\end{prop}

\demo
 On choisit $P,\alpha$ et $\beta$ comme dans la prop.~\ref{psi1} et on note $M$
le $\Lambda(\Gamma)$-module $w_{\delta}({\cal C}^{\alpha})
 \cap {\cal C}^{\beta}$. On choisit ensuite $m_1\geq m(D_0)$ tel que pour tout $b\geq m_1$:

 $\bullet$ 
${\cal C}^{?}$ est inclus dans
 $D_0^{(0,r_b]}\boxtimes_{\delta}\zpet$ et
 $\oe^{(0,r_b]}(\Gamma)\otimes_{\Lambda(\Gamma)} {\cal C}^{?}=D_0^{(0,r_b]}\boxtimes_{\delta}\zpet$, si $?\in \{\alpha,\beta\}$
 (un tel $m_1$ existe, c.f.~\cite[cor. V.1.13]{Cbigone}).

 $\bullet$ L'inclusion de $M$ dans ${\cal C}^{\beta}$ induit un isomorphisme
$$\oe^{(0,r_b]}(\Gamma)\otimes_{\Lambda(\Gamma)} M=
 \oe^{(0,r_b]}(\Gamma)\otimes_{\Lambda(\Gamma)} {\cal C}^{\beta}$$
  (cette condition est automatique
car ${\cal C}^{\beta}/M$ est tu\'e par une puissance de $\sigma_{1+p}-1$ puisqu'il
de longueur finie sur $\O_L$ d'apr\`es la prop.~\ref{psi1}, et $\sigma_{1+p}-1$ est inversible
dans $\oe^{(0,r_b]}(\Gamma)$).

 Alors $\oe^{(0,r_b]}(\Gamma)\otimes_{\Lambda(\Gamma)} M=D_0^{(0,r_b]}\boxtimes_{\delta}\zpet$ et $w_{\delta}(M)\subset  {\cal C}^{\alpha}\subset D_0^{(0,r_b]}\boxtimes_{\delta}\zpet$ (par d\'efinition de $M$), ce qui permet de conclure, en utilisant
 la $i_{\delta}$-semi-lin\'earit\'e de $w_{\delta}$.

\begin{coro}\label{d0rb}
 Si $b>m_1$, alors $D_0^{(0,r_b]}\boxtimes_{\delta}\p1$ est stable par ${\rm GL}_2(\zp)$ et 
\linebreak
$D^{\dagger}\boxtimes_{\delta}\p1$ est stable
 par $G$.
\end{coro}

\demo
 Cf.~\cite[lemme~II.1.10]{Cbigone}: c'est une cons\'equence formelle des formules donnant l'action de $G$ sur $D_0\boxtimes_{\delta}\p1$ (cf. rem.~\ref{squelette}), de la prop.~\ref{engendre}
et des inclusions $\psi(D_0^{(0,r_{b}]})\subset D_0^{(0,r_{b-1}]}$,
 $\varphi(D_0^{(0,r_{b-1}]})\subset D_0^{(0,r_{b}]}$ et $\sigma_a(D_0^{(0,r_b]})\subset D_0^{(0,r_b]}$
qui impliquent que $D_0^{(0,r_b]}$ est stable par ${\rm Res}_{p\Z_p}=\varphi\circ\psi$.

\medskip

  Si $b\geq m> m_1$, on note $\tau_m=\sigma_{1+p^m}-1$ et $$M_m^{\dagger,b}=(1+T)\varphi^m(D_0^{\dagger,b}) \quad \text{et} \quad M_m^{(0,r_b]}=(1+T)\varphi^m(D_0^{(0,r_b]}).$$ 
 Remarquons que $\left(\begin{smallmatrix} 1+p^m & 0\\0 & 1\end{smallmatrix}\right)-1$ agit comme $\tau_m$ sur
  $D_0$ et $\left(\begin{smallmatrix} 1 & p^m\\0 & 1\end{smallmatrix}\right)-1$ y 
  agit par multiplication par $\varphi^m(T)$.
  Le r\'esultat suivant (\cite[prop. V.1.14]{Cbigone} et sa preuve) compare les deux actions, ce qui est fondamental pour la suite.
  Rappelons que l'on a fix\'e une base $(e_i)_i$ du $\oe^{\dagger,b}$-module $D_0^{\dagger,b}$, et que
  $\Gamma_m=\chi^{-1}(1+p^{m}\zp)$. Si $m\geq 2$, on d\'efinit
  des anneaux $\oe^{\dagger,b}(\Gamma_m)$, etc, en rempla\c{c}ant simplement la variable $T$ par $\tau_m$.

\begin{prop}\label{mmdagger}
 Il existe $m_2>m_1$ tel que pour tous $b\geq m\geq m_2$ on ait

 {\rm (i)} $\tau_m$ est bijectif sur $M_m^{(0,r_b]}$ et $\tau_m^a$ induit une bijection de
 $M_m^{\dagger,b}$ sur $\varphi^m(T)^a\cdot M_m^{\dagger,b}$, pour tout $a\in\mathbf{Z}$.

 {\rm (ii)} $M_m^{\dagger,b}$ (resp. $M_m^{(0,r_b]}$) est un $\oe^{\dagger,b}(\Gamma_m)$ (resp. $\oe^{(0,r_b]}(\Gamma_m))$-module libre de base $((1+T)\varphi^m(e_i))_{i}$.

\end{prop}

  On fixe un tel $m_2$ et on le note simplement $m$. Voir le \S~\ref{coefM} pour les notations $K_m$, $a_m^{+}$, $a_m^{-}$, $b_m^{\alpha}$, etc.

\begin{lemm}\label{crucial}
 Il existe une constante $c\geq 1$ telle que:

 {\rm (i)}  $w_{\delta}(\tau_m^a M_m^{\dagger,b})\subset \tau_m^{a-c}M_m^{\dagger,b}$ pour tous $b\geq m$ et $a\in\mathbf{Z}$.

 {\rm (ii)} $(g-1)^n(\tau_m^a M_m^{\dagger,b})\subset \tau_m^{n+a-c}M_m^{\dagger,b}$ pour tous $b\geq m$, $a\in\mathbf{Z}$, $n\geq 0$ et $g\in K_m$.
\end{lemm}

\demo {\rm (i)} Notons $f_i:=(1+T)\varphi^m(e_i)\in D_0^{(0,r_{2m}]}$.
 Alors (prop. \ref{engendre}) $w_{\delta}(f_i)\in D_0^{(0, r_{2m}]}=D_0^{\dagger,2m}[\frac{1}{T}]$. On fixe
$c'$ tel que $$w_{\delta}(f_i)\in \varphi^m(T)^{-c'+l(D_0)}D_0^{\dagger, 2m}$$ pour
 $1\leq i\leq d$ (voir le lemme \ref{psi} pour $l(D_0)$). Comme $w_{\delta}$ commute \`a ${\rm Res}_{1+p^m\zp}$ (car $1+p^m\zp$ est stable
 par $w$), il existe $g_i\in D_0$ tels que
   $w_{\delta}(f_i)=(1+T)\varphi^m(g_i)$. Alors $\varphi^m(g_i)\in \varphi^m(T)^{-c'+l(D_0)} D_0^{\dagger,2m}$,
   donc $g_i\in T^{-c'} D_0^{\dagger,m}$ (utiliser le lemme \ref{psi} et l'identit\'e $g_i=\psi^m(\varphi^m(g_i))$) et donc finalement 
   $w_{\delta}(f_i)\in \varphi^{m}(T)^{-c'} M_m^{\dagger,m}=
   \tau_m^{-c'} M_m^{\dagger,m}$ (prop. \ref{mmdagger}).

   Soient enfin $b\geq m$, $a\in\mathbf{Z}$ et notons $X=M_m^{\dagger,b}$. Comme
   les $e_i$ forment une base de $D_0^{\dagger,b}$ sur $\oe^{\dagger,b}$, les $f_i$ forment
    une base de $X$ sur $\oe^{\dagger,b}(\Gamma_m)$ (prop. \ref{mmdagger}). En utilisant la $i_{\delta}$-semi-lin\'earit\'e
    de $w_{\delta}$ et le fait que $w_{\delta}(f_i)\in \tau_m^{-c'}X$, on obtient $w_{\delta}(\tau_m^a X)\subset \tau_m^{a-c'}X$, ce qui permet de conclure.

    {\rm (ii)} On va montrer que $c=8c'$ marche (avec $c'$ comme dans la preuve de {\rm (i)}, dont on garde les notations).
    Le {\rm (i)} de prop. \ref{mmdagger} montre que $(g-1)^n(\tau_m^a X)=\tau_m^{a+n}X$ si
 $a\in\mathbf{Z}$ et $g\in\{a_m^+, u_m^+\}$. En combinant cela avec le {\rm (i)}, on obtient pour $g\in \{u_m^+, a_m^+\}$
 $$(wgw-1)^n (\tau_m^a X)=w(g-1)^n w(\tau_m^a X)\subset w(\tau_m^{a+n-c'} X)
    \subset \tau_m^{a+n-2c'}X$$ et donc $b_m^{\alpha}(\tau_m^a X)\subset \tau_m^{|\alpha|+a-c}X$.

   Soit maintenant $g\in K_m$ quelconque et \'ecrivons
      $(g-1)^n=\sum_{\alpha\in\mathbf{N}^4} c_{\alpha} b^{\alpha}$ dans $\Lambda(K_m)$.
     Alors $c_{\alpha}\in \zp$ et $v_p(c_{\alpha})\geq n-|\alpha|$
      quand $|\alpha|<n$. Comme $p$ est multiple de $\tau_m^{n_b}$ (et donc de $\tau_m$)
      dans $\oe^{\dagger,b}(\Gamma_m)$, on obtient
      $c_{\alpha} b^{\alpha} (\tau_m^a X)\subset \tau_m^{\max(n,|\alpha|)+a-c}X$
      pour tout $\alpha\in\mathbf{N}^4$. Comme $X$ est complet pour la topologie
      $\tau_m$-adique (car $\oe^{\dagger,b}(\Gamma_m)$ l'est), cela permet
      de conclure que $(g-1)^n(\tau_m^a X)\subset \tau_m^{n+a-c}X$, ce qui finit la preuve.

\begin{prop}\label{action}
 Il existe $c_1>c$ tel que:

{\rm (i)} Pour tout $a\in\mathbf{Z}$ on a $w_{\delta}((T^aD_0^{\dagger,b})^{\psi=0})\subset (T^{a-c_1}D_0^{\dagger,b})^{\psi=0}$.

{\rm (ii)} Pour tous $b>2m$, $a\in\mathbf{Z}$, $n\geq 1$ et $g\in K_m$
  $$(g-1)^n (T^a D_0^{\dagger,b})\subset T^{a+p^m n-c_1} D_0^{\dagger,b}.$$
\end{prop}

\demo On va montrer que l'on peut prendre $c_1=p^m(1+c+l(D_0))$.
Fixons $b>2m$, $a\in\mathbf{Z}$, $n\geq 1$ et notons, pour simplifier, $q=[\frac{a}{p^m}]$
et $Y=\tau_m^{q-l(D_0)} M_m^{\dagger, b-m}$.

\begin{lemm}\label{decomp} Soit $A$ (resp. $B$) l'ensemble des $i\in \{0,1,...,p^m-1\}$ tels que $p$ ne divise pas $i$ (resp. $p$ divise $i$). Alors
 $T^a D_0^{\dagger,b}\subset \sum_{i\in A} \left(\begin{smallmatrix} i & 0\\0 & 1\end{smallmatrix}\right)Y+
 \sum_{i\in B} \left(\begin{smallmatrix} 1 & i-1\\0 & 1\end{smallmatrix}\right)Y$ et $(T^a D_0^{\dagger,b})^{\psi=0}
 \subset \sum_{i\in A} \left(\begin{smallmatrix} i & 0\\0 & 1\end{smallmatrix}\right)Y$.
\end{lemm}

\demo
 Soit $z\in T^a D_0^{\dagger,b}$ et posons $z_i=\psi^m((1+T)^{-i} z)$, de telle sorte que
 $z=\sum_{i=0}^{p^m-1} (1+T)^i\cdot \varphi^m(z_i)$ et $z_i\in T^{q-l(D_0)}D_0^{\dagger, b-m}$
 (lemme \ref{psi}). On en d\'eduit (prop. \ref{mmdagger}) que $x_i=(1+T)\varphi^m(\sigma_{\frac{1}{i}}(z_i))$
 (pour $i\in A$) et $y_i=(1+T)\varphi^m(z_i)$ (pour $i\in B$) sont des \'el\'ements de
 $Y$ et on conclut en remarquant que $z= \sum_{i\in A} \left(\begin{smallmatrix} i & 0\\0 & 1\end{smallmatrix}\right)x_i+
 \sum_{i\in B} \left(\begin{smallmatrix} 1 & i-1\\0 & 1\end{smallmatrix}\right)y_i$. La deuxi\`eme assertion s'en
 d\'eduit, car si $\psi(z)=0$, alors $z_i=0$ pour tout $i\in B$.

\medskip
 Revenons \`a la preuve de la prop.~\ref{action}. En appliquant le {\rm (ii)} du lemme \ref{crucial}, la prop. \ref{mmdagger} {\rm (i)} et le lemme \ref{hyper} (dans cet ordre) on obtient, pour $g\in K_m$,
\begin{align*}
   (g-1)^n(Y)&\subset \tau_m^{q-l(D_0)-c+n} M_m^{\dagger,b-m}=\varphi^m(T)^{q-c-l(D_0)+n} M_m^{\dagger, b-m}\\
   &\subset T^{p^m(q-c-l(D_0)+n)} D_0^{\dagger,b}\subset T^{a+p^m n-c_1}D_0^{\dagger,b}.
\end{align*}
   On conclut pour le {\rm (ii)} en utilisant le lemme \ref{decomp} et le fait que $K_m$ est distingu\'e dans ${\rm GL}_2(\zp)$.
 Le {\rm (i)} se d\'emontre de la m\^eme mani\`ere, en utilisant le {\rm (i)} du lemme \ref{mmdagger}.

\subsection{Fin de la preuve du th.~\ref{mainanal}}

On note
   $p_{\Pi}: D\boxtimes_{\delta}\p1\to \Pi$ la projection naturelle. Elle envoie $D_{0}\boxtimes_{\delta}\p1$ dans $\Pi_0$. 
   Coupl\'ee avec la prop.~\ref{relevanal}, la proposition ci-dessous permet de conclure quant \`a la preuve du th.~\ref{mainanal}.

\begin{prop}\label{presqueanal} Il existe une constante $c_2$ telle que si $a\in\mathbf{Z}$, $b> 2m+1$ et $z\in (T^a D_0^{\dagger,b})\boxtimes_{\delta}\p1$, alors $v:=p_{\Pi}(z)\in \Pi^{(b)}$ et $v^{(b)}(v)\geq ar_b-c_2$.

\end{prop}

\demo Si $z=(z_1,z_2)\in (T^a D_0^{\dagger,b})\boxtimes_{\delta}\p1$, alors 
$z=z_1+w\cdot \varphi(\psi(z_2))$ et 
$\varphi(\psi(z_2))\in T^{a-p(1+l(D))} D_0^{\dagger,b}$ (lemme 
   \ref{psi}). Il suffit donc de d\'emontrer la proposition pour $z\in T^aD_0^{\dagger,b}$. Nous aurons besoin du
lemme suivant:

  \begin{lemm}\label{rban}
 Si $z\in D_0\boxtimes_{\delta}\p1$, alors
 $v_{\Pi}(p_{\Pi}(z))\geq k$ si et seulement si
\linebreak
   $\{\check{z},z\}_{\p1}\in p^k\O_L$ pour tout $\check{z}\in \check{D}_0^{\natural}\boxtimes_{\delta^{-1}}\p1$.

\end{lemm}

\demo
 Par dualit\'e de Schikhof \cite{ST1}, le vecteur $v=p_{\Pi}(z)$ de $\Pi_0$ est dans $p^k\Pi_0$ si et seulement si
$l(v)\in p^k\O_L$ pour tout $l\in \Pi_0^\dual $. Le r\'esultat suit du fait que $\{\,\,,\,\}_{\p1}$ induit un isomorphisme
$\Pi_0^\dual =\check{D}_0^{\natural}\boxtimes_{\delta^{-1}}\p1$.

\medskip
Revenons \`a la d\'emonstration de la prop.~\ref{presqueanal}.
  Les lemmes \ref{rban} et \ref{weak}, la $K_m$-\'equivariance de $\{\,\,,\,\}$ et l'inclusion $D_0^{\natural}\subset T^{-l_1}D_0^{\dagger,b-1}$
  ram\`enent la preuve de la prop.~\ref{presqueanal} \`a celle de l'assertion suivante: il existe une constante $C$ telle que
    $$\lim_{|\alpha|\to\infty} v_p(\{b_m^{\alpha}\check{z}, z\})-p^m r_b|\alpha|=\infty, \quad
   \inf_{\alpha} \left( v_p(\{b_m^{\alpha}\check{z}, z\})-p^m r_b|\alpha|\right)\geq ar_b-C$$
 pour tous $a\in \Z$, $b>2m+1$, $z\in T^a D_0^{\dagger,b}$ et $\check{z}\in T^{-l_1}\check{D}_0^{\dagger,b-1}$.

Nous allons montrer que $C=l_1+4c_1$ convient. Comme
    $\check{z}\in T^{-l_1}\check{D}_0^{\dagger,b-1}$
et $T^{-l_1}\check{D}_0^{\dagger,b-1}\subset T^{-l_1}\check{D}_0^{\dagger,b}$, 
on d\'eduit de la prop.~\ref{action} que 
    $$ b_{m}^{\alpha}\check{z} \in
    T^{p^m|\alpha|-C}\check{D}_0^{\dagger,b-1}\subset T^{p^m|\alpha|-C}\check{D}_0^{\dagger,b}.$$ 
    L'in\'egalit\'e $ \inf_{\alpha} \left( v_p(\{b_m^{\alpha}\check{z}, z\})-p^m r_b|\alpha|\right)\geq ar_b-C$
d\'ecoule alors du lemme \ref{residupadic}. 

   Il nous reste \`a montrer que $\lim_{|\alpha|\to\infty} v_p(\{b_m^{\alpha}\check{z}, z\})-p^m r_b|\alpha|=\infty$.
   Tout
   \'el\'ement $f$ de $\oe^{\dagger,b}$ peut s'\'ecrire sous la forme $$f=\sum_{k\geq 0} f_k \left(\frac{p}{T^{n_b}}\right)^k,$$
   avec $f_k\in \oe^{\dagger,b-1}$ tendant vers $0$ pour la topologie $p$-adique, donc on peut \'ecrire
   $$z=\sum_{k\geq 0} p^{u_k} y_k \left(\frac{p}{T^{n_b}}\right)^k,$$ avec
   $y_k\in T^aD_0^{\dagger,b-1}$ et $u_k\in\mathbf{N}$ tendant vers $\infty$. Notons
    $$x_{k,\alpha}=\left\{b_m^{\alpha} \check{z}, p^{u_k} y_k \left(\frac{p}{T^{n_b}}\right)^k\right\}.$$
 On a $$p^{u_k} y_k \left(\frac{p}{T^{n_b}}\right)^k\in p^{k+u_k}T^{a-kn_b} D_0^{\dagger,b-1}\subset p^{k+u_k}T^{a-kn_b} D_0^{\dagger,b}$$
 et, comme on l'a d\'ej\`a vu,
     $$ b_{m}^{\alpha}\check{z} \in
    T^{p^m|\alpha|-C}\check{D}_0^{\dagger,b-1}\subset T^{p^m|\alpha|-C}\check{D}_0^{\dagger,b}.$$ 
Combin\'ees avec le lemme \ref{residupadic}, l'\'egalit\'e $n_b r_{b-1}=p$ et les in\'egalit\'es 
  $ar_b, ar_{b-1}\geq -|a|$ et $Cr_{b-1}, Cr_b\leq C$, les relations pr\'ec\'edentes donnent
  $$v_p(x_{k,\alpha})-p^m r_b |\alpha|\geq u_k-|a|-C+\max(0,  (p-1)(p^m |\alpha|r_b-k)).$$
 Un petit exercice d'analyse r\'eelle montre alors que $\inf_{k} (v_{p}(x_{k,\alpha})-p^mr_b|\alpha|)$ 
tend vers~$+\infty$ quand $|\alpha|\to\infty$, ce qui permet de conclure.

  \section{Le module $D_{\rm rig}\boxtimes_{\delta}\p1$ et l'espace $\Pi_{\delta}(D)^{\rm an}$}\label{drigp1}
  
    On fixe dans ce chapitre une paire $G$-compatible $(D,\delta)$, avec $D\in\fget({\cal E})$, et on note 
    $\Pi=\Pi_{\delta}(D)$ et $\check{\Pi}=\Pi_{\delta^{-1}}(\check{D})$. On 
  construit une extension non triviale $D_{\rm rig}\boxtimes_{\delta}\p1$ de
$\Pi^{\rm an}$ par $(\check{\Pi}^{\rm an})^\dual$.
Cette construction n'est pas utilis\'ee dans le chapitre suivant, consacr\'e
\`a la preuve du th.~\ref{maincompletion}, mais est tr\`es utile pour
une \'etude fine de $\Pi^{\rm an}$ (cf.~\cite{Cvectan,Dosp,DP}, par exemple).
            
      \subsection{Continuit\'e de l'action de $w_{\delta}$}
      
     Soit $D_0$ un $\oe$-r\'eseau de $D$ stable par $\varphi$ et $\Gamma$. 
       Soit $m$ comme apr\`es la prop. \ref{mmdagger} et soient $a\geq b> 2m$.
 Fixons une base $e_1,e_2,...,e_d$ de $D_0^{\dagger,b}$ sur
   $\oe^{\dagger,b}$. C'est aussi une base de $D^{]0,r_b]}$ sur ${\cal E}^{]0,r_b]}$, ce qui nous permet de
   poser $$v^{[r_a,r_b]}(z)=\min_{1\leq i\leq d} v^{[r_a,r_b]}(f_i)\quad \text{si}\quad z=\sum_{i=1}^{d} f_i e_i\in D^{]0,r_b]}.$$
Rappelons que $D\boxtimes\Z_p^\dual=D^{\psi=0}$
(et de m\^eme si on remplace $D$ par $D^{(0,r_b]}$ pour $b$ assez grand).

\begin{prop} \label{sixpointcinq} Il existe une constante $c$ telle que pour tous $a\geq b>2m$ et tout $z\in (D^{(0,r_b]})^{\psi=0}$ on ait 
$$ v^{[r_a,r_b]}(w_{\delta}(z))\geq v^{[r_a,r_b]}(z)-c.$$

\end{prop}

\demo On peut multiplier $z$ par une puissance de $p$ sans changer l'in\'egalit\'e, donc on peut supposer que 
$z\in D_0^{(0,r_b]}$ et $[v^{[r_a,r_b]}(z)]=N\geq 1$. Nous aurons besoin du lemme suivant:

\begin{lemm}\label{penible}
  Soient $a, b, N\in \mathbf{N}^\dual $ tels que $a\geq b$ et soit $f\in \oe^{(0,r_b]}$. Si
  $v^{[r_a,r_b]}(f)\geq N$, alors 
  $$f\in \sum_{i=0}^{N-1} p^{N-1-i} T^{in_a} \oe^{\dagger,b}.$$

\end{lemm}

\demo
 Ecrivons $$f=\sum_{n<0} a_n T^n+\sum_{n=0}^{n_a-1} a_n T^n+\sum_{n=n_a}^{2n_a-1} a_nT^n+...+\sum_{n\geq (N-1)n_a} a_n T^n.$$
 Par hypoth\`ese $v_p(a_n)+nr_a\geq N$ et $v_p(a_n)+nr_b\geq N$ pour tout $n$. En particulier $v_p(a_n)>N$ si $n<0$, donc 
 $\sum_{n<0} a_n T^n\in p^{N-1} \oe^{\dagger,b}$ (lemme \ref{oedagger}). Ensuite, si 
 $0\leq n<n_a$, on a $v_p(a_n)>N-1$, donc $\sum_{n=0}^{n_a-1} a_n T^n\in p^{N-1} \oe^{+}\subset 
 p^{N-1}\oe^{\dagger,b}$. Le m\^eme argument montre que 
 $\sum_{n=n_a}^{2n_a-1} a_n T^n\in p^{N-2} T^{n_a} \oe^{\dagger,b}$,..., $\sum_{n\geq (N-1)n_a} a_n T^n\in 
 T{(N-1)n_a} \oe^{\dagger,b}$.

\medskip
  Revenons \`a la preuve de la prop.~\ref{sixpointcinq}.
D'apr\`es le lemme \ref{penible} on a 
$z\in \sum_{i=0}^{N-1} p^{N-1-i} T^{in_a} D_0^{\dagger,b}$. Puisque $\psi(z)=0$, on a 
$$z={\rm Res}_{\zpet}(z)\in \sum_{i=0}^{N-1} p^{N-1-i} {\rm Res}_{\zpet} (T^{in_a} D_0^{\dagger,b}).$$
Le lemme \ref{psi} donne l'existence d'une constante $c_2$ telle que 
${\rm Res}_{\zpet} (T^{in_a} D_0^{\dagger,b})\subset (T^{in_a-c_2}D_0^{\dagger,b})^{\psi=0}$ pour tous
$a\geq b>2m$ et tout $i$. Le {\rm (i)} de la prop.~\ref{action} fournit une constante $c_1$ telle que 
$w_{\delta}((T^d D_0^{\dagger,b})^{\psi=0})\subset T^{d-c_1} D_0^{\dagger,b}$ pour tous $a\geq b>2m$ et 
$d\in \mathbf{Z}$. On a donc, avec $c=c_1+c_2$, 
$$w_{\delta}(z)\in \sum_{i=0}^{N-1} p^{N-1-i} T^{in_a-c} D_0^{\dagger,b},$$
et donc $$v^{[r_a,r_b]}(w_{\delta}(z))\geq \inf_{0\leq i<N} (N-1-i+(in_a-c)r_a)\geq N-1-c>v^{[r_a,r_b]}(z)-c-2,$$
d'o\`u le r\'esultat.

      \begin{coro}\label{penible2}
 L'involution $w_{\delta}$ de $(D^{(0,r_b]})^{\psi=0}$ s'\'etend de mani\`ere unique en une involution 
      continue de $(D^{]0,r_b]})^{\psi=0}$ pour tout $b>2m$.

      \end{coro}
      
      \demo Le module $(D^{(0,r_b]})^{\psi=0}=\oplus_{i=1}^{p-1}(1+T)^i \varphi( D^{(0,r_{b-1}]})$ est dense dans $(D^{]0,r_b]})^{\psi=0}=\oplus_{i=1}^{p-1}(1+T)^i \varphi( D^{]0,r_{b-1}]})$, puisque 
      $D^{(0,r_{b-1}]}$ l'est dans $D^{]0,r_{b-1}]}$. Cela d\'emontre l'unicit\'e de l'extension \'eventuelle de $w_{\delta}$. L'existence est une cons\'equence 
      de la proposition pr\'ec\'edente, de la densit\'e de $(D^{(0,r_b]})^{\psi=0}$ dans $(D^{]0,r_b]})^{\psi=0}$ et de la compl\'etude de $(D^{]0,r_b]})^{\psi=0}$.

\medskip
      
       Le cor.~\ref{penible2} fournit une involution continue $w_{\delta}$ sur le module $D_{\rm rig}^{\psi=0}=\cup_{b>2m} (D^{]0,r_b]})^{\psi=0}$, qui \'etend l'involution 
       $w_{\delta}$ sur $(D^{\dagger})^{\psi=0}$. On d\'efinit alors, de la mani\`ere usuelle
  $$D_{\rm rig}\boxtimes_{\delta}\p1=\{(z_1,z_2)\in D_{\rm rig}\times D_{\rm rig},\ 
  {\rm Res}_{\zpet}(z_2)=w_{\delta}({\rm Res}_{\zpet}(z_1))\},$$
 que l'on munit de la topologie induite par l'inclusion $D_{\rm rig}\boxtimes_{\delta}\p1\subset D_{\rm rig}\times D_{\rm rig}$. 
 Notons que l'application $z\mapsto ({\rm Res}_{\zp}(z), \psi({\rm Res}_{\zp}(wz)))$ induit un isomorphisme d'espaces vectoriels topologiques 
 $D_{\rm rig}\boxtimes_{\delta}\p1\simeq D_{\rm rig}\times D_{\rm rig}$, l'application inverse \'etant donn\'ee par
 $(z_1,z_2)\mapsto (z_1, \varphi(z_2)+w_{\delta}({\rm Res}_{\zpet}(z_1)))$. La densit\'e de $D^{\dagger}$ dans $D_{\rm rig}$ entra\^\i ne donc celle de 
 $D^{\dagger}\boxtimes_{\delta}\p1$ dans $D_{\rm rig}\boxtimes_{\delta}\p1$.

\Subsection{L'action de $G$ sur $D_{\rm rig}\boxtimes_{\delta}\p1$}
   \begin{prop}
    L'action de $G$ sur $D^{\dagger}\boxtimes_{\delta}\p1$ s'\'etend par continuit\'e en une action continue de $G$ sur 
    $D_{\rm rig}\boxtimes_{\delta}\p1$.
   \end{prop}
   \demo
    Les formules du squelette d'action (voir la rem.~\ref{squelette}) permettent de d\'efinir une action de 
    $G$ sur $D_{\rm rig}\boxtimes_{\delta}\p1$ (le fait qu'il s'agit bien d'une action d\'ecoule de la densit\'e de 
    $D^{\dagger}\boxtimes_{\delta}\p1$ dans $D_{\rm rig}\boxtimes_{\delta}\p1$ et du fait que ces formules d\'efinissent 
    une action de $G$ sur $D^{\dagger}\boxtimes_{\delta}\p1$). La continuit\'e de l'action se d\'emontre de la m\^eme mani\`ere que 
la prop.~\ref{sixpointcinq}, en utilisant le {\rm (ii)} de la prop.~\ref{action}.
      
\medskip
    On renvoie 
      au \S~\ref{dh} pour les alg\`ebres ${\cal D}(K_m)$ et ${\cal D}_h(K_m)$, et au $\S$ \ref{coefM} pour les~$b_m^{\alpha}$.

\begin{prop}
   Il existe une constante $c$ telle que pour tous $a\geq b>2m$, $z\in D^{]0,r_b]}$ et 
   $\alpha\in\mathbf{N}^4$ on ait 
   $$v^{[r_a,r_b]}(b_m^{\alpha}z)\geq v^{[r_a,r_b]}(z)+p^{m}|\alpha|r_a-c.$$
\end{prop}

\demo La preuve est enti\`erement analogue \`a celle de la prop.~\ref{sixpointcinq}, en utilisant le lemme 
\ref{penible} et le {\rm (ii)} de la prop.~\ref{action}.

\begin{coro}\label{estimrarb}
Pour tous $a\geq b>2m$, $z\in D^{]0,r_b]}$ et $\lambda=\sum_{\alpha\in\mathbf{N}^4} c_{\alpha} b_m^{\alpha}\in {\cal D}_{a-m}(K_m)$, la s\'erie $\sum_{\alpha} c_{\alpha} b_m^{\alpha}z$ converge dans $D^{]0,r_b]}$ et
$$v^{[r_a,r_b]}(\sum_{\alpha} c_{\alpha} b_m^{\alpha}z)\geq v^{[r_a,r_b]}(z)+v^{(a-m)}(\lambda)-c.$$

\end{coro}

\demo Une suite de $D^{]0,r_b]}$ converge dans $D^{]0,r_b]}$ si et seulement si elle converge pour la valuation 
$v^{[r_a,r_b]}$ pour tous $a\geq b$. Le r\'esultat suit donc de 
 la proposition pr\'ec\'edente et de la d\'efinition de $v^{(a-m)}$.

\smallskip

\begin{prop} \label{distrlocal}
  Soit $H$ un sous-groupe ouvert compact de $G$.

{\rm (i)}  Si $b$ est assez grand, l'action de $H$ sur $D^{]0,r_b]}\boxtimes_{\delta}\p1$
 s'\'etend en une structure de ${\cal D}(H)$-module topologique.

 {\rm (ii)}  $D_{\rm rig}\boxtimes_{\delta}\p1$ est un ${\cal D}(H)$-module topologique.

\end{prop}

\demo
  Comme $H$ est commensurable \`a $K_m=1+p^m {\rm M}_2(\zp)$ (avec $m$ comme ci-dessus), on peut supposer que $H=K_m$.
Si $z=(z_1,z_2)\in D^{]0,r_b]}\boxtimes_{\delta}\p1$, on peut \'ecrire $z$ sous
la forme $z=z_1+w\cdot {\rm Res}_{p\zp}(z_2)$, avec $z_1, {\rm Res}_{p\zp}(z_2)\in D^{]0,r_b]}$. En passant \`a la limite projective (sur $a$) dans
   le cor.~\ref{estimrarb}, on obtient une application continue
   ${\cal D}(K_m)\times (D^{]0,r_b]}\boxtimes_{\delta}\p1)\to D^{]0,r_b]}\boxtimes_{\delta}\p1$,
   d\'efinie par $(\lambda,z)\mapsto \sum_{\alpha\in\mathbf{N}^4}
   c_{\alpha} b_m^{\alpha}z$ si $\lambda=\sum_{\alpha\in \mathbf{N}^4} c_{\alpha}b_m^{\alpha}$.
   Cette application \'etend la structure de $L[K_m]$-module 
\linebreak
de
   $D^{]0,r_b]}\boxtimes_{\delta}\p1$, et comme $L[K_m]$ est dense
dans ${\cal D}(K_m)$, cela prouve que $D^{]0,r_b]}\boxtimes_{\delta}\p1$
   est un 
${\cal D}(K_m)$-module (topologique d'apr\`es ce qui pr\'ec\`ede).
Ceci d\'emontre le (i) et, le (ii) \'etant une cons\'equence imm\'ediate
du (i), cela permet de conclure.

\begin{prop}\label{reslocal}
Soit $H$ un sous-groupe ouvert compact de $G$, qui stabilise l'ouvert compact $U\subset \p1(\qp)$. Alors
$D_{\rm rig}\boxtimes_{\delta} U$ est un sous-${\cal D}(H)$-module de $D_{\rm rig}\boxtimes_{\delta}\p1$
 et ${\rm Res}_U(\lambda\cdot z)=\lambda\cdot {\rm Res}_U(z)$
 pour tous $z\in D_{\rm rig}\boxtimes_{\delta}\p1$ et $\lambda\in {\cal D}(H)$.
 \end{prop}

 \demo
   Cela d\'ecoule de la continuit\'e de l'action de ${\cal D}(H)$, de la densit\'e de $L[H]$ dans ${\cal D}(H)$
   et de la $H$-\'equivariance de l'application ${\rm Res}_U$.

\subsection{Description de $\Pi^{\rm an}$ via $D_{\rm rig}\boxtimes_{\delta}\p1$}

   L'accouplement $\{\,\,,\,\}_{\p1}$
sur 
\linebreak
$(\check{D}^{\dagger}\boxtimes_{\delta^{-1}}\p1)\times (D^{\dagger}\boxtimes_{\delta}\p1)$
 s'\'etend en un accouplement $G$-\'equivariant parfait (voir la discussion qui pr\'ec\`ede \cite[prop. V.2.10]{Cbigone})
 $$ \{\,\,,\,\}_{\p1}:(\check{D}_{\rm rig}\boxtimes_{\delta^{-1}}\p1)\times (D_{\rm rig}\boxtimes_{\delta}\p1)\to L.$$

\begin{theo}\label{mainan}  $(\Pi^{\rm an})^\dual $ est isomorphe comme $G$-module topologique \`a l'orthogonal
 de $D^{\natural}\boxtimes_{\delta}\p1$ dans $\check{D}_{\rm rig}\boxtimes_{\delta^{-1}}\p1$.

\end{theo}

\demo
 Soit $M$ l'orthogonal de $D^{\natural}\boxtimes_{\delta}\p1$ dans $\check{D}_{\rm rig}\boxtimes_{\delta^{-1}}\p1$.
Nous aurons besoin du lemme suivant:

\begin{lemm}\label{contpairing}
 Pour tout $\check{z}\in \check{D}_{\rm rig}\boxtimes_{\delta^{-1}}\p1$ 
 l'application $ D^{\dagger}\boxtimes_{\delta}\p1\to L$, 
donn\'ee par $z\mapsto \{\check{z}, z\}_{\p1}$, est continue. 
De plus, l'application 
 $\check{D}_{\rm rig}\boxtimes_{\delta^{-1}}\p1\to (D^{\dagger}\boxtimes_{\delta}\p1)^\dual $ ainsi obtenue est continue. 
\end{lemm}

\demo
 Il suffit de v\'erifier que pour tout $\check{z}\in \check{D}_{\rm rig}$ l'application 
 $z\mapsto \{\check{z}, z\}$ est une forme lin\'eaire continue sur $D^{\dagger}$ et que 
 l'application $\check{D}_{\rm rig}\to (D^{\dagger})^\dual $ ainsi obtenue est continue. En revenant aux d\'efinitions
 des topologies de $D^{\dagger}$ et $\check{D}_{\rm rig}$, la continuit\'e de $\check{D}_{\rm rig}\to (D^{\dagger})^\dual $ d\'ecoule de
l'in\'egalit\'e 
 $$v_p(\{\check{z}, z\})\geq kr_b+v^{[r_a,r_b]}(\check{z})-2$$
 pour $a\geq b>2m$, $k\geq 1$, $z\in T^k D_0^{\dagger,b}$ et $\check{z}\in D^{]0,r_b]}$.
 Pour d\'emontrer cette in\'egalit\'e on se ram\`ene par densit\'e et 
 $L$-lin\'earit\'e (et en utilisant le {\rm (ii)} du lemme \ref{valrarb}) \`a $\check{z}\in \frac{1}{p} T^{Nn_b}D_0^{\dagger,b}$, avec $N$ la partie enti\`ere de 
 $v^{[r_a,r_b]}(\check{z})\geq 0$. L'in\'egalit\'e suit alors du lemme~\ref{residupadic}. 

\medskip 
Revenons \`a la preuve du th.~\ref{mainan}. On d\'eduit du lemme \ref{contpairing}, de l'isomorphisme
$\Pi^{\rm an}\simeq (D^{\dagger}\boxtimes_{\delta}\p1)/(D^{\natural}\boxtimes_{\delta}\p1)$ (cor. \ref{anali})
et de la d\'efinition de $M$ une application lin\'eaire continue $\phi: M\to (\Pi^{\rm an})^\dual $, induite par 
$\{\,\,,\,\}_{\p1}$. Explicitement, on a $\langle \phi(\check{z}), v\rangle=\{\check{z}, z_v\}_{\p1}$ pour tout rel\`evement 
$z_v\in D^{\dagger}\boxtimes_{\delta}\p1$ de $v\in \Pi^{\rm an}$ et tout $\check{z}\in M$.
L'application 
$\phi$ est $G$-\'equivariante, puisque 
$\{\,\,,\,\}_{\p1}$ l'est. Nous allons montrer que $\phi$ est un hom\'eomorphisme, en construisant son inverse.

Commen\c cons par constater que $\phi$ est injective car un \'el\'ement
de ${\rm Ker}(\phi)$ est orthogonal \`a $D^{\dagger}\boxtimes_{\delta}\p1$
et donc \`a $D_{\rm rig}\boxtimes_{\delta}\p1$ par
densit\'e de $D^{\dagger}\boxtimes_{\delta}\p1$, et donc est nul puisque
$\{\ ,\ \}_{\p1}$ est un accouplement parfait.
     
Le (ii) de la prop.~\ref{distrlocal} montre que l'inclusion 
       $\Pi^\dual \simeq\check{D}^{\natural}\boxtimes_{\delta^{-1}}\p1\subset \check{D}_{\rm rig}\boxtimes_{\delta^{-1}}\p1$ induit une application ${\cal D}(H)$-lin\'eaire continue 
       $\xi:{\cal D}(H)\otimes_{\Lambda(H)} \Pi^\dual  \to \check{D}_{\rm rig}\boxtimes_{\delta^{-1}}\p1$. Puisque 
       $\check{\Pi}^\dual $ et $\Pi^\dual $ sont orthogonaux, l'image de $\xi$ est contenue dans $M$. 

 La prop.~\ref{dualdepih} fournit un isomorphisme d'espaces vectoriels topologiques 
$\iota:(\Pi^{\rm an})^\dual \simeq {\cal D}(H)\otimes_{\Lambda(H)} \Pi^\dual $ ($H$ \'etant par exemple 
${\rm GL}_2(\zp)$), et la compos\'ee $\iota\circ\phi\circ\xi$ est l'identit\'e
car c'est l'identit\'e sur $\Pi^\dual$. 

Comme $\phi$ est injective, cela implique que son inverse est $\xi\circ\iota$,
ce qui permet de conclure.

\begin{coro}  $(\Pi^{\rm an})^\dual $ et $(\check{\Pi}^{\rm an})^\dual $ sont exactement orthogonaux
 pour l'accouplement $\{\,\,,\,\}_{\p1}$.

\end{coro}

\demo Le th.~\ref{mainan} et le cor.~\ref{checkstab} montrent que $(\Pi^{\rm an})^\dual $ est l'orthogonal de $\check{\Pi}^\dual $ dans $\check{D}_{\rm rig}\boxtimes_{\delta^{-1}}\p1$.
 Or $\check{\Pi}^\dual $ est un sous-espace dense de $(\check{\Pi}^{\rm an})^\dual $ (par densit\'e de $\check{\Pi}^{\rm an}$ dans $\check{\Pi}$ combin\'ee
   \`a la r\'efl\'exivit\'e de $\check{\Pi}^{\rm an}$ et au th\'eor\`eme de Hahn-Banach), donc $(\Pi^{\rm an})^\dual $
   est en fait l'orthogonal de  $(\check{\Pi}^{\rm an})^\dual $, ce qui permet de conclure.

\begin{coro} L'injection $(\check{\Pi}^{\rm an})^\dual \to D_{\rm rig}\boxtimes_{\delta}\p1$ fournie par le th.~\ref{mainan} induit une suite exacte de $G$-modules topologiques
    $$0\to (\check{\Pi}^{\rm an})^\dual \to D_{\rm rig}\boxtimes_{\delta}\p1\to \Pi^{\rm an}\to 0.$$

\end{coro}

\demo
  D'apr\`es le th.~\ref{mainan}, $(\check{\Pi}^{\rm an})^\dual $ est un sous-espace ferm\'e de $D_{\rm rig}\boxtimes_{\delta}\p1$.
   Soit $Y$ le quotient. Puisque $\{\,\,,\,\}_{\p1}$ induit une dualit\'e parfaite entre
   $\check{D}_{\rm rig}\boxtimes_{\delta^{-1}}\p1$ et $D_{\rm rig}\boxtimes_{\delta}\p1$, on
   obtient un isomorphisme topologique de $Y^\dual $ sur l'orthogonal de $(\check{\Pi}^{\rm an})^\dual $ dans
   $\check{D}_{\rm rig}\boxtimes_{\delta^{-1}}\p1$, donc sur $(\Pi^{\rm an})^\dual $ (corollaire pr\'ec\'edent).
   On a donc un isomorphisme de $G$-modules topologiques $Y^\dual \simeq (\Pi^{\rm an})^\dual $ et on conclut en observant que 
   $\Pi^{\rm an}$ et $Y$ sont r\'eflexifs (pour le dernier, cela d\'ecoule de ce que 
  $\check{D}_{\rm rig}\boxtimes_{\delta^{-1}}\p1$
    satisfait Hahn-Banach).

\begin{coro}
   Il existe $m=m(D)$ tel que $(\Pi^{\rm an})^\dual \subset \check{D}^{]0, r_m]}\boxtimes_{\delta^{-1}}\p1$. 

\end{coro}

\demo
$(\Pi^{\rm an})^\dual $ est un espace de Fr\'echet et le corollaire pr\'ec\'edent fournit une injection continue dans 
$\check{D}_{\rm rig}\boxtimes_{\delta^{-1}}\p1$, qui est la r\'eunion croissante des espaces de Fr\'echet 
$\check{D}^{]0, r_m]}\boxtimes_{\delta^{-1}}\p1$. Le r\'esultat s'ensuit.

\section{Compl\'et\'es unitaires universels}

\subsection{R\'eseaux invariants minimaux}

  Soit $G$ un groupe de Lie $p$-adique
 et soit $\Pi$ une repr\'esentation continue de $G$ sur un $L$-espace vectoriel localement convexe. 
  Rappelons qu'une $L$-repr\'esentation de Banach $B$ de $G$ est dite unitaire si $G$ pr\'eserve une valuation d\'efinissant la topologie de $B$. 
   Le {\it compl\'et\'e unitaire universel} $\widehat\Pi$ de $\Pi$ est
(s'il existe) une $L$-repr\'esentation de Banach unitaire de $G$, 
munie d'une application $L$-lin\'eaire continue, $G$-\'equivariante $\iota: \Pi\to \widehat{\Pi}$, qui est universelle au sens suivant: pour toute 
$L$-repr\'esentation de Banach unitaire $B$ de $G$, l'application 
$${\rm Hom}_{L[G]}^{\rm cont}(\widehat{\Pi}, B)\to {\rm Hom}_{L[G]}^{\rm cont}(\Pi, B), \quad f\mapsto f\circ \iota$$
est une bijection. Autrement dit, tout morphisme continu $\Pi\to B$ se factorise de mani\`ere unique \`a travers $\iota: \Pi\to \widehat{\Pi}$.

\begin{rema}
(i) 
Il d\'ecoule facilement de la d\'efinition que si $\widehat{\Pi}$ existe, alors l'image de $\iota$
est dense dans $\widehat \Pi$, et que $\widehat\Pi$ est unique \`a isomorphisme 
unique pr\`es.

(ii) M\^eme si $\widehat\Pi$ existe, il n'y a aucune raison
a priori pour que $\widehat\Pi\neq 0$, et classifier les repr\'esentations
de $G$ ayant un compl\'et\'e universel non nul est un probl\`eme
difficile et fondamental.

(iii) Si $\Pi$ est topologiquement irr\'eductible et si $\widehat\Pi$ existe et est non nul, alors
$\Pi$ admet une valuation invariante par $G$. En effet, dans ce cas l'application naturelle $\Pi\to\widehat\Pi$ est injective,
ce qui permet de consid\'erer la restriction de la valuation sur $\widehat\Pi$ \`a $\Pi$.

(iv) L'espace ${\rm LA}(\Z_p)$
des fonctions localement analytiques sur $\Z_p$,
vu comme repr\'esentation de $\Z_p$ (par $(b\cdot\phi)(x)=\phi(x-b)$),
n'a pas de compl\'et\'e unitaire universel; par contre, le m\^eme espace,
vu comme une repr\'esentation du semi-groupe $P^+=\matrice{\Z_p-\{0\}}{\Z_p}{0}{1}$,
admet comme compl\'et\'e unitaire universel l'espace ${\cal C}^0(\Z_p)$ des
fonctions continues sur $\Z_p$ [\,l'action de
$P^+$ est donn\'ee par $\big(\matrice{a}{b}{0}{1}\cdot\phi)\big)(x)=0$ si
$x\notin b+a\Z_p$, et $\big(\matrice{a}{b}{0}{1}\cdot\phi)\big)(x)=
\phi\big(\frac{x-b}a\big)$ si $x\in b+a\Z_p$\,].
\end{rema}

Si $\Pi$ est un $L$-espace vectoriel, un {\it r\'eseau} de $\Pi$
est un sous-$\O_L$-module de $\Pi$ qui engendre le $L$-espace
vectoriel~$\Pi$ (on ne demande pas \`a un r\'eseau d'\^etre s\'epar\'e
pour la topologie $p$-adique; un r\'eseau peut donc contenir des
sous-$L$-espaces vectoriels). La remarque suivante d'Emerton \cite[lemma 1.3]{Emcompl} sera utile pour la suite:

\begin{lemm}\label{Em}
  $\Pi$ admet un compl\'et\'e universel $\widehat\Pi$ si et seulement si
  $\Pi$ contient un $\O_L$-r\'eseau $M$ avec les propri\'et\'es suivantes:

  $\bullet$ $M$ est ouvert dans $\Pi$ et stable sous l'action de $G$.

  $\bullet$ $M$ est minimal pour ces propri\'et\'es, i.e. $M$ est contenu dans un homoth\'etique de tout sous-$\O_L$-r\'eseau
ouvert
de $\Pi$ stable par $G$.

   De plus, dans ce cas $$\widehat\Pi= L\otimes_{\O_L}\linv(M/p^nM).$$
\end{lemm}

\begin{prop}\label{coherent}
Si $\Pi$ est coh\'erente, alors $\Pi^{\rm an}=\cup_{h} \Pi^{(h)}$ admet un compl\'et\'e
universel $\widehat{\Pi^{\rm an}}$.  Plus pr\'ecis\'ement, si $h\geq m(\Pi)$ et si
$\Pi^{(h)}_0$  est la boule unit\'e de $\Pi^{(h)}$ pour la valuation $v^{(h)}$, alors 
$$M:=\sum_{g\in G} g\cdot\Pi^{(h)}_0$$
est un r\'eseau ouvert de $\Pi^{\rm an}$ et
$\widehat{\Pi^{\rm an}}$ est le compl\'et\'e de $\Pi$ relativement
\`a ce r\'eseau.
\end{prop}

\demo
  Soient $k\geq h\geq m(\Pi)$. Pour chaque $s$ tel que $d(s,1)\leq k-h$ on choisit
  un rel\`evement $\hat{s}$ \`a $G$.
   L'application continue 
  $$\bigoplus_{d(s,1)\leq k-h}\Pi^{(h)} \to \Pi^{(k)}, \quad 
   (x_s)_{d(s,1)\leq k-h}\mapsto \sum_{d(s,1)\leq k-h} \hat{s}\cdot x_s$$
   est surjective par hypoth\`ese. On en d\'eduit que 
    $\sum_{d(s,1)\leq k-h}s\cdot \Pi^{(h)}_0$ est un r\'eseau ouvert
de $\Pi^{(k)}$. En passant \`a la limite inductive, il s'ensuit que $M=\sum_{s\in S}s\cdot \Pi^{(h)}_0$
est un r\'eseau ouvert de $\Pi$, invariant par $G$ par
construction. Par ailleurs, si ${\cal L}$ est un r\'eseau ouvert de $\Pi^{\rm an}$, alors
${\cal L}\cap \Pi^{(h)}$ est un r\'eseau ouvert de $\Pi^{(h)}$, donc il existe 
$k$ tel que ${\cal L}\supset p^k\Pi^{(h)}_0$. Si de plus 
${\cal L}$ est invariant par $G$, alors ${\cal L}$ contient 
$p^kM$. Ainsi, $M$ est un r\'eseau ouvert $G$-invariant et minimal \`a homoth\'etie pr\`es.
Le lemme~\ref{Em} permet de conclure. 

\begin{rema}
Il n'est pas difficile de voir, en reprenant
les arguments ci-dessus, que la valuation
$v^{(h)}$ sur $\Pi^{(h)}$ est
g\'en\'eratrice au sens
d'Emerton~\cite[def.~1.13]{Emcompl}.
Autrement dit, la coh\'erence implique l'existence
d'une valuation g\'en\'eratrice et donc la prop.~\ref{coherent}
peut aussi se d\'eduire de
la prop.~1.14 de~\cite{Emcompl}.
\end{rema}

\subsection{Fonctorialit\'e}

 Si $u\in {\rm Hom}_{L[G]}^{\rm cont}(\Pi_1,\Pi_2)$,
et si $\Pi_1$ et $\Pi_2$ admettent des compl\'et\'es universels,
alors il existe un unique morphisme
$\hat u\in {\rm Hom}_{L[G]}^{\rm cont}(\widehat \Pi_1, \widehat \Pi_2)$,
tel que $\hat u\circ\iota_1=\iota_2\circ u$, o\`u $\iota_i:\Pi_i\to\widehat\Pi_i$ est
l'application canonique.

\begin{prop} \label{exactness}
Soit $0\to\Pi_1\to\Pi\to\Pi_2\to 0$ une suite exacte stricte de repr\'esentations de $G$ 
sur des $L$-espaces vectoriels localement convexes. Si 
$\widehat{\Pi}$ existe, alors:

{\rm (i)}  $\widehat{\Pi}_2$ existe aussi, et le morphisme $\widehat{\Pi}\to\widehat{\Pi}_2$ induit par 
$\Pi\to \Pi_2$ est surjectif. 

{\rm (ii)} Si de plus $\widehat{\Pi}_1$ existe, alors 
${\rm Im}(\widehat \Pi_1\to \widehat \Pi)$ est dense dans
${\rm Ker}(\widehat{\Pi}\to\widehat{\Pi}_2)$. 

\end{prop}

\demo {\rm (i)} Soit $M$ un r\'eseau ouvert, $G$-stable et minimal (\`a homoth\'etie pr\`es) de $\Pi$, comme dans le lemme \ref{Em}.
Comme $\Pi\to\Pi_2$ est stricte et surjective,
l'image $M_2$ de $M$ dans $\Pi_2$
est
un r\'eseau ouvert de $\Pi_2$, stable par $G$.
Si $M'_2$ est un autre r\'eseau ouvert de $\Pi_2$, stable par $G$,
et si $M'$ est l'image inverse de $M'_2$ dans $\Pi$, alors
$M'\cap M$ est un r\'eseau ouvert de $\Pi$ qui est stable par
$G$; comme $M$ est minimal, il existe $k\in\N$ tel que
$M'\cap M$ contienne $p^k M$, et donc $M'$ contient
$p^kM_2$.  Il s'ensuit que $M_2$ est minimal
(\`a homoth\'etie pr\`es) et le r\'esultat suit du lemme~\ref{Em}, qui montre aussi que 
$\linv\, M_2/p^nM_2$ est un r\'eseau ouvert born\'e de
$\widehat \Pi_2$.

{\rm (ii)} $M_1=\Pi_1\cap M$, est
un r\'eseau ouvert de $\Pi_1$, stable par $G$,
et on a une suite exacte $0\to M_1\to M\to M_2\to 0$.
En passant aux s\'epar\'es compl\'et\'es pour la topologie
$p$-adique, puis en inversant $p$, on en d\'eduit la suite
exacte
$$0\to L\otimes_{\O_L}(\linv\, M_1/p^nM_1)\to\widehat\Pi\to\widehat\Pi_2\to 0.$$
L'application naturelle $M_1\to \linv\, M_1/p^nM_1$ induit un morphisme continu
\linebreak
$f: \Pi_1\to  L\otimes_{\O_L}(\linv M_1/p^nM_1)$, d'image dense. 
Par d\'efinition de 
$\hat{\Pi}_1$, $f$ est induit par un unique morphisme
$\varphi: \widehat\Pi_1\to L\otimes_{\O_L}(\linv M_1/p^nM_1)$. Puisque 
$f$ est \`a image dense, il en est de m\^eme de $\varphi$, ce qui permet de conclure.

\begin{coro} \label{chapeau} On garde les notations et hypoth\`eses de la prop.~\ref{exactness}, et on suppose de plus que 
 $\widehat{\Pi}_1$ est admissible.
Alors on a une suite exacte de $L$-espaces vectoriels $\widehat{\Pi}_1\to \widehat{\Pi}\to \widehat{\Pi}_2\to 0$. 
\end{coro}

\demo Soient $X={\rm Ker}(\widehat{\Pi}\to \widehat{\Pi}_2)$ et $H$ un sous-groupe
ouvert compact de $G$. La prop.~\ref{exactness} 
montre l'existence d'un morphisme d'image dense $f: \widehat{\Pi}_1\to X$.
Ainsi, $X^\dual $ est un sous-$\Lambda(H)$-module de~$(\widehat{\Pi}_1)^\dual $,
 qui est de type fini par admissibilit\'e de~$\widehat{\Pi}_1$. 
Comme $\Lambda(H)$ est noeth\'erien, 
$X$ est admissible, et 
l'image de $f$ est ferm\'ee \cite{ST1}. 
Donc $f$ est surjectif, ce qui permet de conclure.

\begin{rema}
 $\widehat\Pi_1\to\widehat\Pi$ n'est pas
toujours injective.

\end{rema}

\Subsection{Le compl\'et\'e universel de $\Pi^{\rm an}$}
Supposons dor\'enavant que $G={\rm GL}_2(\Q_p)$.

\begin{theo}\label{COH}
Pour tout $\Pi\in {\rm Rep}_L(G)$, $\Pi^{\rm an}$ est coh\'erente et son compl\'et\'e universel
est $\Pi$. 
\end{theo}

\demo
 La preuve va demander quelques pr\'eliminaires. On commence par supposer que $\Pi=\Pi_{\delta}(D)$ pour un $D\in \fget({\cal E})$. 
Soit $D_0$ un r\'eseau de $D$ et soit $\Pi_0=\Pi_{\delta}(D_0)$, un r\'eseau de $\Pi$, ouvert, born\'e
et $G$-stable. Pour tout $b>m(D)$ on note $X_b$
le sous-$\O_L$-module $(D_0^{\dagger,b}\boxtimes_{\delta}\p1)/(D_0^{\natural}\boxtimes_{\delta}\p1)$
  de $\Pi_0$. Soit $\Pi_0^{(b)}$ la boule unit\'e de~$\Pi^{(b)}$ pour la valuation 
  $v^{(b)}$. Les prop.~\ref{relevanal} et~\ref{presqueanal}, et le fait que $p\in T^{n_b} \oe^{\dagger,b}$, montrent
qu'il existe $b_0>m(D)$ et une constante $c$ tels que $p^{c} \Pi_0^{(b)} \subset X_b\subset p^{-c} \Pi_0^{(b)}$ pour tout 
$b\geq b_0$.

\begin{lemm} \label{coh} Il existe une constante $c_1$ telle que 
pour tout $b\geq b_0$ 
$$p^{c_1} X_b\subset \sum_{d(g,1)\leq b-b_0} g\cdot X_{b_0}.$$
\end{lemm}

\demo Il existe $c_1$ tel que $w_{\delta}({\rm Res}_{\zpet}( D_0^{\dagger,b}))\subset p^{-c_1} D_0^{\dagger,b}$
pour tout $b\geq b_0$ (utiliser la prop.~\ref{action} et le fait que $p\in T^{n_b}\oe\oe $). On a donc $p^{c_1} D_0^{\dagger,b}\subset D_0^{\dagger,b}\boxtimes_{\delta}\p1$.
Ensuite, tout $z\in D_0^{\dagger,b}$ s'\'ecrit sous la forme 
$$z=\sum_{i=0}^{p^{b-b_0}-1} \matrice{1}{i}{0}{1} \matrice{p^{b-b_0}}{0}{0}{1} u_i,$$
avec $u_i=\psi^{b-b_0}((1+T)^{-i}z)\in D_0^{\dagger,b_0}$. Puisque tout 
 $x\in D_0^{\dagger,b}\boxtimes_{\delta}\p1$ s'\'ecrit 
\linebreak
$x={\rm Res}_{\zp}(x)+w\cdot {\rm Res}_{p\zp}(w\cdot x)$,
 on en d\'eduit que 
 $$ p^{c_1} (D_0^{\dagger,b}\boxtimes_{\delta}\p1)\subset \sum_{d(g,1)\leq b-b_0} g\cdot (D_0^{\dagger,b_0}\boxtimes_{\delta}\p1),$$
ce qui permet de conclure.

\begin{lemm} \label{commens}
 On a $\Pi^{\rm an}\cap \Pi_0\subset p^{-(c+c_1+1)}  \sum_{g\in G} g\cdot \Pi_0^{(b_0)}$. 
\end{lemm}

\demo
 Soit $v\in p\cdot (\Pi_0\cap \Pi^{\rm an})$. Alors 
 $v$ a un rel\`evement $z$ \`a $p  D_0\boxtimes_{\delta}\p1$ et, puisque 
 $v\in \Pi^{\rm an}$, le cor.~\ref{anali} montre que $z\in D^{\dagger}\boxtimes_{\delta}\p1$.
 Puisque $p D_0\cap D^{\dagger}\subset \cup_{b\geq b_0}  D_0^{\dagger,b}$
 (cela se d\'eduit du lemme \ref{intersect}), on conclut que $v\in X_b$ pour un certain 
 $b\geq b_0$. Or, le lemme pr\'ec\'edent et l'inclusion $X_{b_0}\subset p^{-c} \Pi_0^{(b_0)}$
 entrainent $$X_b\subset p^{-c-c_1} \sum_{d(g,1)\leq b-b_0} g\cdot \Pi_0^{(b_0)}\subset
 p^{-c-c_1} \sum_{g\in G} g\cdot \Pi_0^{(b_0)},$$
 ce qui permet de conclure.

\medskip
  Revenons \`a la preuve du th.~\ref{COH}. En tensorisant par $L$ l'inclusion du lemme \ref{coh}
 on obtient $\Pi^{(b)}\subset \sum_{d(g,1)\leq b-b_0} g\cdot \Pi^{(b_0)}$ pour tout 
 $b\geq b_0$, ce qui montre que $\Pi^{\rm an}$ est coh\'erente. La prop.~\ref{coherent} montre
 que $\widehat{\Pi^{\rm an}}$ existe, et c'est le compl\'et\'e de $\Pi^{\rm an}$ par rapport au r\'eseau
 minimal (\`a homoth\'etie pr\`es) $\sum_{g\in G} g\cdot \Pi_0^{(b_0)}$. Le lemme 
 \ref{commens} montre que $M=\Pi^{\rm an}\cap \Pi_0$ 
est commensurable \`a  $\sum_{g\in G} g\cdot \Pi_0^{(b_0)}$,
 donc $\widehat{\Pi^{\rm an}}$ est le compl\'et\'e de $\Pi^{\rm an}$ par rapport au r\'eseau 
 $M$. Puisque $\Pi^{\rm an}$ est dense dans $\Pi$, 
  l'injection naturelle $M/p^nM \to \Pi_0/p^n\Pi_0$ est un isomorphisme. En passant \`a la limite et en inversant $p$, on obtient $\widehat{\Pi^{\rm an}}=\Pi$.
  
   Jusque l\`a nous avons suppos\'e que $\Pi=\Pi_{\delta}(D)$ pour une paire $G$-compatible $(D,\delta)$. Supposons maintenant que 
   $\Pi\in{\rm Rep}_L(G)$ est quelconque. Disons que $\Pi$ est {\it bonne} si elle
est coh\'erente et \'egale au compl\'et\'e unitaire universel de ses vecteurs localement analytiques. 
Notons que si $\Pi$ est une repr\'esentation coh\'erente, alors $\widehat{\Pi^{\rm an}}$ existe (prop.~\ref{coherent})
et l'injection $\Pi^{\rm an}\to \Pi$ induit un morphisme $\widehat{\Pi^{\rm an}}\to \Pi$. Notons aussi qu'une repr\'esentation de dimension finie est bonne.

Le th.~\ref{recoverPi} nous fournit une paire $G$-compatible $(D,\delta)$ et une application 
\linebreak
   $\beta: \Pi_{\delta}(D)\to \Pi/\Pi^{{\rm SL}_2(\qp)}$, dont le noyau et le conoyau sont de dimension finie sur~$L$. 
Comme $\Pi^{{\rm SL}_2(\qp)}$ est de dimension finie (cor. \ref{dimfinsl}), et comme
   $\Pi_{\delta}(D)$ est bonne d'apr\`es ce qui pr\'ec\`ede, le r\'esultat
s'obtient en appliquant plusieurs fois le lemme suivant.

\begin{lemm}
Soit $0\to\Pi_1\to\Pi\to\Pi_2\to 0$ une suite exacte dans ${\rm Rep}_L(G)$. Si $\Pi_1$ est bonne et si 
$\Pi$ ou $\Pi_2$ est bonne, l'autre l'est aussi. 
\end{lemm}

\demo L'exactitude du foncteur $\Pi\mapsto \Pi^{\rm an}$ nous fournit une suite exacte 
$0\to \Pi_1^{\rm an}\to \Pi^{\rm an}\to \Pi_2^{\rm an}\to 0$. Puisque 
la coh\'erence est stable par quotient et extension (prop.~\ref{extcoh}), les hypoth\`eses faites entra\^\i nent la coh\'erence de 
$\Pi_1, \Pi$ et $\Pi_2$ et donc (prop. \ref{coherent}) l'existence de 
$\widehat{\Pi_1^{\rm an}}$, $\widehat{\Pi^{\rm an}}$ et 
$\widehat{\Pi_2^{\rm an}}$. De plus, comme $\Pi_1$ est bonne, 
$\widehat{\Pi_1^{\rm an}}\simeq \Pi_1$ est admissible. Le cor.~\ref{chapeau} fournit donc une suite exacte 
$\widehat{\Pi_1^{\rm an}}\to \widehat{\Pi^{\rm an}}\to \widehat{\Pi_2^{\rm an}}\to 0$, s'ins\'erant dans un diagramme commutatif
$$\xymatrix{
&\widehat{\Pi_1^{\rm an}}\ar[r]\ar[d]&\widehat{\Pi^{\rm an}}\ar[r]\ar[d]&\widehat{\Pi_2^{\rm an}}\ar[r]\ar[d]&0\\
0\ar[r]&\Pi_1\ar[r]&\Pi\ar[r]&\Pi_2\ar[r]&0}.$$
Par hypoth\`ese la fl\`eche verticale de gauche est un isomorphisme. 
Une chasse au diagramme montre que si une des fl\`eches verticales 
restantes est un isomorphisme,
l'autre l'est aussi, ce qui permet de conclure.


\begin{thebibliography}{94}

\bibitem{Am64}
{\sc Y. Amice}, Interpolation $p$-adique, Bull. Soc. France
{\bf 92} (1964), 117--180.


\bibitem{BL} {\sc L.~Barthel} et {\sc R.~Livn\'e}, Irreducible modular representations of ${\rm GL}_2$ of a local field, 
Duke Math.~J.~75 (1994), 261-292.


\bibitem{Br1} {\sc C.~Breuil}, Sur quelques repr\'esentations modulaires et p-adiques de ${\rm GL}_2(\qp)$ I, 
Compositio Math. 138 (2003), 165-188.

\bibitem{Br} {\sc C.~Breuil}, Invariant \textit{L}
 et s\'erie sp\'eciale $p$-adique, Ann. Sci. Ecole Norm. Sup
 37 (2004),  559-610.
 
 
\bibitem{BC} {\sc L.~Berger} et {\sc P.~Colmez}, Familles de repr\'esentations de de Rham et monodromie p-adique,
Ast\'erisque 319 (2008),  303-337.




\bibitem{CCsurconv} {\sc F.~Cherbonnier} et {\sc P.~Colmez}, Repr\'esentations
$p$-adiques surconvergentes, Invent. Math. 133 (1998),  581-611.

\bibitem{Cmirab} {\sc P.~Colmez}, $(\varphi,\Gamma)$-modules et repr\'esentations
du mirabolique de $G$, Ast\'erisque 330 (2010),  61-153.

\bibitem{Cbigone} {\sc P.~Colmez}, Repr\'esentations de ${\rm GL}_2(\qp)$ et
$(\varphi,\Gamma)$-modules, Ast\'erisque 330 (2010),  281-509.

\bibitem{Cserieunit} {\sc P.~Colmez}, La s\'erie principale unitaire de
${\rm GL}_2(\qp)$, Ast\'erisque 330, 213-262.

\bibitem{Cvectan} {\sc P.~Colmez}, La s\'erie principale unitaire de
${\rm GL}_2(\qp)$: vecteurs localement analytiques, preprint.


\bibitem{Dosp}
{\sc G.~Dospinescu}, Actions infinit\'esimales dans la correspondance
de Langlands locale $p$-adique pour ${\bf GL}_2(\Q_p)$,
Math.~Ann.~{\bf 354} (2012), 627--657.

\bibitem{DBenjamin} {\sc G.~Dospinescu} et {\sc B.~Schraen},  Endomorphism algebras of admissible
$p$-adic representations of $p$-adic Lie groups, Representation Theory (\`a para\^\i tre).

\bibitem{DP} {\sc G.~Dospinescu}, en pr\'eparation.



\bibitem{Sautoy} {\sc K.~Dixon}, {\sc M.~ Du Sautoy}, {\sc A.~Mann} et {\sc D.~Segal}, 
Analytic pro-$p$-groups, $2$nd edition, Cambridge University Press, 1999.

\bibitem{Emcompl} {\sc M.~Emerton}, $p$-adic $L$-functions and unitary completions of representations of $p$-adic reductive groups,
 Duke Math. J. 130 (2005), no. 2, 353-392.

\bibitem{Emlocan} {\sc M.~Emerton}, Locally analytic vectors in representations of locally p-adic analytic groups,
to appear in Memoirs of the AMS, disponible \`a \url{http://www.math.northwestern.edu/~emerton/preprints.html}.

\bibitem{FoGrot} {\sc J.-M.~Fontaine}, Repr\'esentations
$p$-adiques des corps locaux. I, in {\it The
Grothendieck Festschrift, Vol II}, Progr. Math. 87, Birkhauser, 1990,  249-309.


\bibitem{LXZ} {\sc R. Liu}, {\sc B. Xie}, et {\sc Y. Zhang},
 Locally Analytic Vectors of Unitary Principal Series of $\rm{GL}_2(\qp)$, 
Ann. E.N.S.~{\bf 45} (2012), 167--190.


\bibitem{KK1} {\sc K.~Kedlaya},
A $p$-adic monodromy theorem, Ann. of Math.~{\bf 160} (2004),
93--184.

\bibitem{KiAst} {\sc M.~Kisin}, Deformations of $G_{\qp}$
and $G$ representations, Ast\'erisque 330,
2010,  511-528.

\bibitem{Pa} {\sc V.~Pa\v{s}k\={u}nas}, The image
of Colmez's Montr\'eal functor, Publ.~IHES (\`a para\^\i tre).

\bibitem{Garcia} {\sc G.~Perez-Garcia} et {\sc W.~Schikhof}, {\it Locally convex spaces over non-archimedean 
valued fields}, Cambridge studies in advanced mathematics, 2010. 


\bibitem{ST1} {\sc P.~Schneider} et {\sc J.~Teitelbaum}, Banach space representations
and Iwasawa theory, Israel J. Math. 127 (2002), 359-380.


\bibitem{STInv} {\sc P.~Schneider} et {\sc J.~Teitelbaum}, 
Algebras of $p$-adic distributions and admissible
representations, Invent. Math. 153 (2003),  145-196.


\end{thebibliography}
\end{document}